%% file: BR08b.tex
\documentclass[envcountsect,envcountsame]{llncs}

\def\isPCversion{1} 

\RequirePackage[OT1]{fontenc}

\usepackage{a4wide}

\usepackage{amsmath, amssymb}

\usepackage{amsfonts}

\usepackage{graphicx,pstricks}

\usepackage{comment}

\usepackage[round]{natbib}

\include{macrosbis}

\def\blfootnote{\xdef\@thefnmark{}\@footnotetext}

\begin{document}

\begin{frontmatter}

\title{
Adaptive FDR control under independence and dependence}
\titlerunning{
Adaptive FDR control under independence and dependence}

\author{Gilles Blanchard\inst{1} \and \'Etienne Roquain\inst{2}}

\institute{
Fraunhofer FIRST.IDA, Berlin, Germany,\\
\email{blanchar@first.fraunhofer.de}
\and
University of Paris 6, LPMA, 4, Place Jussieu, 75252 Paris cedex 05, France\\
\email{etienne.roquain@upmc.fr}
}

\authorrunning{Blanchard, G., and Roquain, E.}

\maketitle

\begin{abstract}

In the context of multiple hypotheses testing, the proportion $\pi_0$
of true null hypotheses in the pool of hypotheses to test 
often plays a crucial role, although it is generally
unknown {\em a priori}. A testing procedure using an implicit or explicit estimate
of this quantity in order to improve its efficency is called {\em adaptive}.
In this paper, we focus on
the issue of False Discovery Rate (FDR) control and we present new
adaptive multiple testing procedures with control of the FDR.
First, in the context of assuming independent $p$-values, we present two new procedures and
give a unified review of other existing adaptive procedures that have provably
controlled FDR. We report
extensive simulation results comparing these procedures
and testing their robustness when the independence
assumption is violated. The new proposed procedures appear
competitive with existing ones. The overall best, though, is reported to
be Storey's estimator, but for a parameter setting that does not
appear to have been considered before.
  Second, we propose adaptive versions of step-up procedures that have provably
  controlled FDR under positive dependences and unspecified dependences of the $p$-values, respectively. 
  While simulations only show an improvement over non-adaptive
  procedures in limited situations, these are to our knowledge among the
  first theoretically founded adaptive multiple testing procedures that control the FDR when the $p$-values
are not independent.
\end{abstract}

\end{frontmatter}

\pagestyle{plain}    
\pagenumbering{arabic}

\section{Introduction}

\subsection{Adaptive multiple testing procedures}

Spurred by an increasing number of application fields, in partilar bioinformatics,
the topic of multiple testing (which enjoys a long history in the
statistics literature) has generated a renewed, growing attention in the recent years.
For example, using microarray data, the goal is to detect which genes (among several ten
of thousands) exhibit a significantly different level of expression in
two different experimental conditions.  Each gene represents a
``hypothesis'' to be tested in the statistical sense.  The genes'
expression levels fluctuate naturally (not to speak of other sources
of fluctuation introduced by the experimental protocol), and, because
they are so many genes to choose from, it is important to control
precisely what can be deemed a significant observed difference.
Generally it is assumed that the natural fluctuation distribution of a
{\em single} gene is known and the problem is to take into account the
number of genes involved (for more details, see for instance
\citealp{DSP2003}).

In this work, we focus on building multiple testing procedures with a
control of the false discovery rate (FDR). This quantity is defined as
the expected proportion of type I errors, that is, the proportion of
true null hypotheses among all the null hypotheses that have been
rejected (i.e. declared as false) by the procedure.  In their seminal
work on this topic, \citet{BH1995} proposed the
celebrated \textit{linear step-up} (LSU) procedure, that is proved to
control the FDR under independence between the $p$-values.  Later, it
was proved \citep{BY2001} that the LSU procedure still controls the FDR
when the $p$-values have positive dependences (or more precisely, a
specific form of positive dependence called PRDS). Under unspecified
dependences, the same authors have shown that the FDR control still
holds if the threshold collection of the LSU procedure is divided by a
factor $1+1/2+\dots+1/m$, where $m$ is the total number of null
hypotheses to test. 

More recently, the latter result has been
generalized \citep{BF2007,BR2008EJS,Sar2006}, by showing that there is a family of
step-up procedures (depending on the choice of a kind of prior distribution)
that still control the FDR under unspecified dependences between the
$p$-values.

However, all of these procedures, which are built in order to control the FDR at a
level $\alpha$, can be showed to have actually their FDR upper bounded by
$\pi_0\alpha$, where $\pi_0$ is the proportion of true null
hypotheses in the initial pool. Therefore, when most of the hypotheses are false (i.e.,
$\pi_0$ is small), these procedures are inevitably
conservative, since their FDR is in actuality much lower than the fixed target $\alpha$. 
In this context, the challenge of \textit{adaptive control} of
the FDR (e.g., \citealp{BH2000,Black2004}) is to
integrate an estimation of the unknown proportion $\pi_0$ in the
threshold of the previous procedures and to prove that the FDR is
still rigorously controlled by $\alpha$.

Adaptive procedures are therefore of practical interest if it is expected that
$\pi_0$ is, or can be, significantly smaller than 1. An example of such a
situation occurs when using hierarchical procedures
(e.g., \citealp{BH2006}) which first selects some clusters of hypotheses that are
likely to contain false nulls, and then apply a multiple testing procedure on the
selected hypotheses. Since  a large part of the true null
hypotheses is expected to be false in the second step, an adaptive procedure
 is needed in order to keep the FDR close to the target level.

A number of adaptive procedures have been proposed in the recent
literature and can loosely be divided into the following categories:
\begin{itemize}
\item {\em plug-in} procedures, where some initial estimator of $\pi_0$ is
directly plugged in as a multiplicative level correction to the usual
procedures. In some cases (e.g. Storey's estimator, see \citealp{Storey2002}), the resulting plug-in adaptive 
procedure (or a variation thereof) has been proved to control the FDR
at the desired level \citep{BKY2006,STS2004}. A variety of other
estimators of $\pi_0$ have been proposed
(e.g. \citealp{MR2006,JC2007,Jin2008}); while their asymptotic consistency
(as the number of hypotheses tends to infinity) is generally established, their use in plug-in
adaptive procedures has not always been studied.
\item {\em two-stage} procedures: in this approach, a first round
of multiple hypothesis testing is performed using some fixed algorithm, then 
the results of this first round are used  in order to tune the parameters of
a second round in an adaptive way. This can generally be interpreted
as using the output of the first stage to estimate $\pi_0$. Different 
procedures following this general approach have been proposed
\citep{BKY2006,Sar2006,Far2007}; more generally, multiple-stage
procedures can be considered.
\item {\em one-stage} procedures, which perform a single round of
  multiple testing (generally step-up or step-down), based on a particular
threshold collection that renders it adaptive \citep{FDR2008,GBS2008}.
\end{itemize}





In addition, some other works \citep{GW2004,STS2004,FDR2008} have studied the 
question of adaptivity to the parameter $\pi_0$\, from an {\em
  asymptotic} viewpoint.
In this framework, the more specific {\em random effects} model is --
most generally, though not always -- considered,
in which $p$-values are assumed independent, each hypothesis has a probability
$\pi_0$ of being true, and all false null hypotheses share the same alternate
distribution. 
The behavior of different procedures is then studied under the
limit where the number of tested hypotheses grows to infinity. 
One advantage of this approach and specific model is that it allows to derive
quite precise results (see \citealp{Neu08}, for a precise study of limiting behaviors
of central limit type under this model, including for some of the new
procedures introduced in the present paper). However, we emphasize that in the present
work our focus is decidedly on the nonasymptotic side, using finite samples and arbitrary alternate
hypotheses.

 To complete this overview, let us also mention another interesting and different direction opened up recently, that of
adaptivity to the alternate distribution. If the alternate distribution
is known {\em a priori}, it is well-known that the optimal testing statistics are
likelihood ratios between the null and the alternate (which can then be
transformed into $p$-values). When the alternate is unknown though,
one can hope to estimate, implicitly or explicitly, the alternate
distribution from the observed data, and consequently approximate the
optimal test statistics (\citealp{SC2007} proposed an asymptotically
consistent approach to this end; see also \citealp{Spjo1972},
\citealp{Storey2007})\,. Interestingly, this point of view is also intimately linked to
some traditional topics in statistical learning 
such as classification and of optimal novelty detection 
(see, e.g., \citealp{SB2009}). However, in the present paper we will
focus on adaptivity to the parameter $\pi_0$ only. 


\subsection{Overview of this paper}

The contributions of the present paper 
are the following.
A first goal of the paper is to introduce a number of novel adaptive procedures:
\begin{enumerate}
\item We introduce a new {\em one-stage} step-up procedure that is more 
powerful than the standard LSU procedure in a large range of
situations, and provably controls the FDR under independence. 
This procedure is called one-stage adaptive, because the estimation of $\pi_0$ is performed implicitly.
\item Based on this, we then build a new {\em two-stage} adaptive procedure, which 
is more powerful in general 
than the procedure proposed by \citet{BKY2006},
while provably controlling the FDR under independence.
\item  Under the assumption of positive or arbitrary dependence of the $p$-values,
we introduce new two-stage adaptive versions of known step-up procedures (namely, of the LSU under positive
dependences, and of the family of procedures introduced by \citealp{BF2007}, under unspecified dependences). 
These adaptive versions provably control the FDR and result in an improvement of power over the 
non-adaptive versions in some situations (namely, when the number of hypotheses rejected in the first
stage is large, typically more than $60\%$). 
\end{enumerate}

A second goal of this work is to present a review  of several existing
adaptive step-up procedures with provable FDR control under independence.
For this, we present the theoretical FDR control as a consequence of a single general theorem for plug-in procedures, 
which was first established by \citet{BKY2006}. Here, we give 
a short self-contained proof of this result 
that is of independent interest. The latter is based on some tools
introduced earlier \citep{BR2008EJS,Roq2007}, that aim to unify FDR control proofs.
Related results and tools also appear independently in \citet{FDR2008,Sar2008}.

A third goal is to compare both the existing and our new adaptive
procedures in an extensive simulation study under both independence
and dependence, following the simulation model and methodology used by
\citet{BKY2006}. 
Concerning the new one- and two- step procedures with theoretical FDR
control under independence, these are generally quite competitive in comparison to existing
ones. However we also report that the best procedure overall (in terms of
power, among procedures that are robust enough to the dependent case)
appears to be
the plug-in procedure based on the well-known Storey estimator \citep{Storey2002}
used with the somewhat nonstandard parameter $\lambda=\alpha$\,.
This outcome was in part unexpected since to
the best of our knowledge, this fact had never been
pointed out so far (the usual default recommended choice is
$\lambda=\frac{1}{2}$ and turns out to be very unstable in dependent
situations); this is therefore an important conclusion of this paper
regarding practical use of these procedures. 

Concerning the new
two-step procedure with theoretical FDR control under dependence,
simulations show an (admittedly limited) improvement over their non-adaptive
counterpart in favorable situations which correspond to what 
was expected from the theoretical study (large proportion of false
hypotheses). The observed improvement is unfortunately not striking
enough to be able to recommend using these procedures in practice;
their interest is therefore at this point mainly theoretical, in that 
 these are to our knowledge among the
  first theoretically founded adaptive multiple testing procedures that control the FDR when the $p$-values
are not independent.

The paper is organized as follows: in Section 2, we introduce the
 mathematical framework, and we recall the existing
 non-adaptive results in FDR control. In Section 3 we deal with the
setup of independent $p$-values. We expose our new procedures and review
the existing ones, and finally compare them in a simulation study.
The case of positive dependent and arbitrarily dependent $p$-values 
is examined in Section 4 where we introduce our new adaptive procedures
in this context. A conclusion is given in Section~5. Section 6 and 7 contains proofs 
of the results and lemmas, respectively. Some technical
remarks and discussions of links to other work are gathered at the end
of each relevant subsection, and can be skipped by the non-specialist reader.

\section{Preliminaries}

\subsection{Multiple testing framework}

In this paper, we will stick to the traditional statistical 
framework for multiple testing.
Let $(\mathcal{X},\mathfrak{X},\Proba)$ be a probability space;  we
want to infer a decision on $\Proba$ from an observation $x$ in
$\mathcal{X}$ drawn from $\Proba$\,. Let $\mathcal{H}$ be a finite set of null hypotheses
for $\Proba$, that is, each null hypothesis $h\in\cH$ corresponds to
some subset of distributions on $(\mathcal{X},\mathfrak{X})$ and
"$\Proba$ satisfies $h$" means that $\Proba$ belongs to this subset of
distributions. The number of null hypotheses $|\cH|$ is denoted by
$m$, where $|.|$ is the cardinality function. The underlying probability $\Proba$ being fixed, we denote
$\mathcal{H}_0=\{h\in\mathcal{H}| \Proba \mbox{ satisfies } h\}$ the
set of the true null hypotheses and $m_0=|\mathcal{H}_0|$ the number
of true null hypotheses. We let also $\pi_0:=m_0/m$ the proportion of
true null
hypotheses. 
We stress that $\cH_0$, $m_0$, and
$\pi_0$ are unknown and implicitly depend on the unknown $\mbp$\,. All
the results to come are always implicitly meant to hold for any
generating distribution $\mbp$\,.

We suppose further that there exists a set of \textit{$p$-value} functions
$\mbf{p}=(p_h,h\in\cH)$, meaning that each
$p_h:(\mathcal{X},\mathfrak{X})\mapsto [0,1]$ is a measurable function
and that for each $h\in\mathcal{H}_0$, $p_h$ is bounded stochastically
by a uniform distribution, that is,
\begin{equation}
\forall h \in \cH_0\, \qquad \forall t\in[0,1] ,\:\:\Proba(p_h \leq t)\leq t.\label{equ_stochunif}
\end{equation}
Typically, $p$-values are obtained from statistics that have a known distribution $P_0$
under the corresponding null hypothesis. In this case, if $F_0$ denotes the corresponding
cumulative distribution function, applying $1-F_0$ to the observed statistic results in a random variable
satisfying \eqref{equ_stochunif} in general. Here, we are however not concerned how these $p$-values
are constructed and only assume that they exist and are known (this is the
standard setting in multiple testing).

\subsection{Multiple testing procedure and errors}

A \textit{multiple testing procedure} is a measurable
function $$R :x\in \mathcal{X} \mapsto R(x)\in  \mathcal{P}(\cH),$$ 
which takes as
input an observation $x$ 
and returns a subset of $\cH$,
corresponding to the rejected hypotheses. As it is commonly the case, we will focus here on
multiple testing procedure based on $p$-values, that is, we will implicitly
assume that $R$ is of the form $R(\mbf{p})$.

A multiple testing procedure $R$ can make two kinds of errors:
a \textit{type I error} occurs for $h$ when $h$ is true and is
rejected by $R$\,, that is, $h\in\cH_0\cap R$. Following the Neyman-Pearson
general philosophy for hypothesis testing, the primary concern
is to control the number of type I errors of a testing procedure.
Conversely, \textit{a type II error} occurs for $h$ when $h$ is 
false and is not rejected by $R$, that is $h\in\cH_0^c\cap R^c$.

The most traditional way to control type I error is to upper bound the
``Family-wise error rate'' (FWER), which is the probability that one or
more true null hypotheses are wrongly rejected. However, procedures with
a controlled FWER are 
very ``cautious'' not to make a single
error, and thus reject only few hypotheses. This conservative way of measuring
the type I error for multiple hypothesis testing can be a serious hindrance
in practice, since it requires 
to collect large enough
datasets so that significant evidence can be found 
under this strict error control criterion.
More recently, a more liberal measure of type I errors
has been introduced in multiple testing \citep{BH1995}:
the {\em false discovery rate} (FDR), which is the averaged
proportion of true null hypotheses in the set of all the rejected
hypotheses:
\begin{definition}[False discovery rate]
The {\em false discovery rate} of a multiple testing procedure $R$ for
a generating distribution $\mbp$ is given by
\begin{equation}
 \FDR(R):=\e{\frac{|R\cap \cH_0|}{|R|} \ind{|R|>0}}\,.\label{equ_FDR}
\end{equation}
\end{definition}

\begin{remark}
Throughout this paper we will use the following convention: whenever there
is an indicator function inside an expectation, this has logical priority over any
other factor appearing in the expectation. What we mean is that if 
other factors include expressions that may not be defined (such as 
the ratio $\frac{0}{0}$) outside
of the set defined by the indicator, this is safely ignored. This results in
more compact notations, such as in the above definition. Note also
again that the dependence of the FDR on the unknown $\mbp$ is implicit.
\end{remark}

A classical aim, then, is to build procedures $R$ with FDR
upper bounded at a given, fixed level $\alpha$. Of course, if we choose
$R=\emptyset$, meaning that $R$ rejects no hypotheses, trivially $\FDR(R)=0\leq
\alpha$\,. Therefore, it is desirable to build procedures $R$ 
satisfying $\FDR(R)\leq
\alpha$ while at the same time having as few type II errors as possible. As a general rule,
provided that $\FDR(R)\leq \alpha$, we
want to build procedures that reject as many false hypotheses as
possible. The absolute power of a multiple testing procedure is
defined as the average proportion of false hypotheses correctly
rejected, $\e{\abs{R\cap \cH_0^c}}/\abs{\cH_0^c}$\,.
Given two procedures $R$ and $R'$\,, a particularly simple sufficient
condition for $R$ to be more powerful than $R'$ is when 
$R'$ if $R'  \subset R$  holds pointwise.
We will say in this case that $R$ is \textit{(uniformly) less conservative} than
$R'$\,.

\subsection{Self-consistency, step-up procedures, FDR control and adaptivity}\label{sec_step-upFDRcontrol}

We first define a general class of multiple testing procedures called
{\em self-consistent procedures} \citep{BR2008EJS}.

\begin{definition}[Self-consistency, non-increasing procedure]
\label{def-sc}
Let $\Delta: \set{0,1,\ldots,m} \rightarrow \mbr$\,, $\Delta(0)=0$\,,
be a function called {\em threshold collection}; a multiple testing procedure $R$ is said
to satisfy the self-consistency condition with respect to $\Delta$ if
\begin{equation}
R\subset\{h\in\cH\telque p_h\leq \Delta(|R|) \}\label{gen_selfcond} 
\end{equation}
holds almost surely.
Furthermore, we say that $R$ is non-increasing if for all $h\in\cH$ the function $p_h\mapsto |R(\mathbf{p})|$ is non-increasing, that is if $|R|$ is non-increasing in each $p$-value.   
\end{definition}

The class of self-consistent procedures includes well-known types of procedures, notably step-up
and step-down.  The assumption that a procedure is non-increasing,
which is required in addition to self-consistency in some of the
results to come, is relatively natural as a lower $p$-value means we
have more evidence to reject the corresponding hypothesis.
We will mainly focus on the {\em step-up} procedure,
which we define now. For this, we sort the $p$-values in increasing order using the notation 
$p_{(1)}\leq \dots \leq p_{(m)}$ and put $p_{(0)}=0$\,. This order is of course itself
random since it depends on the observation. 
\begin{definition}[Step-up procedure]\label{def_stepup}
The \textit{step-up procedure} with threshold collection $\Delta$ 
 is defined as
\[
R=\{h\in\cH\telque p_h \leq p_{(k)}\},\:\mbox{ where
}\:{k}=\max\{ 0\leq i \leq m \telque p_{(i)}\leq \Delta(i)\}.\]
\end{definition}

A trivial but important property of a step-up procedure is the following. 
\begin{lemma}
The step-up procedure with threshold collection $\Delta$ is non-increasing and self-consistent with respect to $\Delta$\,.
\end{lemma}
Therefore, a result valid for any non-increasing self-consistent procedure also holds for the corresponding step-up procedure. This will be used extensively through the paper and thus should be kept in mind by the reader. 

Among all procedures that are self-consistent with
respect to $\Delta$\,, the step-up is uniformly less conservative than
any other \citep{BR2008EJS} and is therefore of primary interest.
However, to recover procedures of a more general form (including step-down for instance), the statements of this paper will be preferably expressed in terms of self-consistent procedures when it is possible.

Threshold collections are generally scaled by the target FDR level
$\alpha$\,. Once correspondingly rewritten under the normalized form $\Delta(i) = \alpha
\beta(i)/m$\,, we will call $\beta$ the {\em shape function} for
threshold collection $\Delta$\,.
In the particular case where
the shape function $\beta$ is the identity function,
the procedure is called
the \textit{linear step-up} (LSU) {\em procedure} (at level $\alpha$).


The LSU plays a prominent role in multiple testing under FDR control;
it was the first procedure for which FDR control was proved and it is probably the most widely used
procedure in this context. More precisely, when the $p$-values are assumed to be 
independent, 
 the following theorem holds. 

\begin{theorem}\label{th_BH}
  Suppose that the $p$-values of $\mathbf{p}=(p_h,h\in\cH)$ are independent. 
  Then any non-increasing self-consistent procedure with respect to threshold collection $\Delta(i) = \alpha i/m$
has FDR upper bounded by $\pi_0\alpha$\,, where $\pi_0=m_0/m$ is the proportion 
of true null hypotheses. (In particular, this is the case for the linear
  step-up procedure). 
  Moreover, if the $p$-values associated to true null hypotheses are exactly distributed like a uniform
  distribution, 
  the linear  step-up procedure has FDR equal to $\pi_0\alpha$\,. 
\end{theorem}


The first part of this result, in the case of the LSU,
 was proved in the landmark paper of \citet{BH1995};
the second part (also for the LSU) was proved by
\citet{BY2001} and \citet{FR2001}.

\citet{BY2001} extended the previous result about FDR
control of the linear step-up procedure to the case of $p$-values with
a certain form of positive dependence called {\em positive regressive
dependence from a subset} (PRDS).  We skip a formal definition for now
(we will get back to this topic in Section \ref{newsu_df}).   The
extension of this result to self-consistent
procedures (in the independent as well as PRDS cases) 
was established by \citet{BR2008EJS} and \citet{FDR2008}.




However, when no particular assumptions are made on the dependences
between the $p$-values, it can be shown that the above FDR control is
not generally true.  This situation is called {\em unspecified} or
{\em arbitrary} dependence.  A modification of the LSU was first
proposed in \cite{BY2001} which was proved to have a controlled FDR
under arbitrary dependence. This result was extended by \cite{BF2007}
and \cite{BR2008EJS}
(see also a related result of \citealp{Sar2006}): it can be shown
that self-consistent procedures (not necessarily non-increasing) based
on a particular class of shape functions have controlled FDR:
\begin{theorem}\label{th_BF}
Under unspecified dependences between the $p$-values of
$\mathbf{p}=(p_h,h\in\cH)$, consider $\beta$ a shape function of the
form:
\begin{equation}
\beta(r)=\int_{0}^{r} u d\nu(u),\label{equ_beta}
\end{equation}
where $\nu$ is some fixed a priori probability distribution on $(0,\infty)$. Then 
any self-consistent procedure with respect to threshold collection $\Delta(i) = \alpha \beta(i)/m$\, 
has FDR upper bounded by $\alpha\pi_0$\,.
\end{theorem}


To recap, in all of the above cases, 
the FDR is actually controlled at the level $\pi_0 \alpha$ instead of the
target $\alpha$.
Hence, a direct corollary  of both of the above theorems 
is that the step-up procedure 
with shape function $\beta^*(x) = \pi_0^{-1} \beta(x)$
has FDR upper bounded 
$\alpha$ in either of the following situations:
\begin{itemize}
\item[-] $\beta(x)=x$ when the $p$-values are independent or PRDS,
\item[-] the shape function $\beta$ is of the form (\ref{equ_beta}) when the $p$-values have unspecified dependences.
\end{itemize}
Note that, since $\pi_0\leq 1$, using $\beta^*$ always gives rise to a
less conservative procedure than using $\beta$ 
(especially when $\pi_0$ is small). However, since
$\pi_0$ is unknown, the shape function $\beta^*$ is not directly accessible.
We therefore will call the step-up procedure using $\beta^*$ the
\textit{Oracle step-up procedure} based on shape function $\beta$
(corresponding to one of the above cases).

Simply put, the role of adaptive step-up procedures is to mimic the latter oracle in order
to obtain more powerful procedures. 
Adaptive procedures are often step-up procedures using
the modified shape function $G\beta $\,,
where $G$ is some estimator of $\pi_0^{-1}$:

\begin{definition}[Plug-in adaptive step-up procedure]\label{defadap}
Given a level $\alpha\in(0,1)$, a \textit{shape function} $\beta$ and
an estimator $G:[0,1]^{\cH}\rightarrow (0,\infty)$ of the quantity $\pi_0^{-1}$\,,
the \textit{plug-in adaptive step-up procedure} of shape function
$\beta$ and using estimator $G$ (at level $\alpha$) 
is defined as
\[
R=\{h\in\cH\telque p_h \leq p_{(k)}\},\:\mbox{ where
}\:k=\max\{i\telque p_{(i)}\leq \alpha \beta(i)G(\mathbf{p})/m\}.\]
The (data-dependent) function $\Delta(\bp,i)=\alpha
\beta(i)G(\mathbf{p})/m$ is called the \textit{adaptive threshold collection}
corresponding to the procedure. In the particular case where the shape
function $\beta$ is the identity function on $\mathbb{R}^+$, the
procedure is called an \textit{adaptive linear step-up procedure}
using estimator $G$ (and at level $\alpha$). 
\end{definition}

Following the previous definition, an adaptive plug-in procedure is composed of two different steps:
\begin{enumerate}
\item Estimate $\pi_0^{-1}$ with an estimator $G$\,.
\item Take the step-up procedure of shape function $G \beta$\,.
\end{enumerate}
A subclass of plug-in adaptive procedures is formed by so-called 
{\em two-stage procedures}, when the estimator $G$ is actually based on a first,
non-adaptive, multiple testing procedure. This can obviously be possibly iterated and
lead to multi-stage procedures. The distinction between generic plug-in procedures
and two-stage procedures is somewhat informal and generally meant only to provide
some kind of nomenclature between different possible approaches.

The main theoretical task is to ensure that an adaptive procedure of this type still correctly
controls the FDR. 
The mathematical difficulty obviously comes from the additional random variations
of the estimator $G$ in the procedure.

\section{Adaptive procedures with provable FDR control under independence}\label{section_indepadapt}

In this section, we introduce two new adaptive procedures that provably
control the FDR under independence. The first one is one-stage and does not include an
explicit estimator of $\pi_0^{-1}$\,, hence it is not a plug-in
procedure. We then propose to use this as the first stage in a new
two-stage procedure, which constitutes the second proposed method.

For clarity, we first introduce the new one-stage procedure; we then discuss
several possible plug-in procedures, including our new proposition and
several procedures proposed by other authors. FDR control for these
various plug-in procedures can be studied using a general theoretical device
introduced by \cite{BKY2006} which  we reproduce here with 
a self-contained and somewhat simplified proof. Finally,  to compare these different approaches, we close this
section with extensive simulations which both examined  the performance under independence and 
the robustness under (possibly strong) positive correlations.

\subsection{New adaptive one-stage step-up procedure}\label{sec:onestageadapt}

We present here our first main contribution,
a one-stage adaptive step-up procedure. This means that the
estimation step is implicitly included in the shape function $\beta$\,.

\begin{theorem}\label{th_main1}
Suppose that the $p$-values of $\mathbf{p}=(p_h,h\in\cH)$ are
independent and let $\lambda\in (0,1)$ be fixed. Define the adaptive 
threshold collection 
\begin{equation} \label{equ_onestage_lambda}
\Delta(i)=\min\bigg((1-\lambda)\frac{\alpha i}{m-i+1},\lambda\bigg)\,.
\end{equation}
Then any non-increasing self-consistent procedure with respect to $\Delta$
has FDR upper bounded by $\alpha$\,. In particular,
this is the case of the corresponding step-up procedure,
denoted by \textit{BR-1S-$\lambda$}\,.
\end{theorem}

The above result will be proved in
Section~\ref{sec_proofadaptproc}. Our proof is in part based on
Lemma~1 of \cite{BKY2006}. Note that an alternate proof of Theorem
\ref{th_main1} has been established in \cite{Sar2008} without using this lemma while
nicely connecting the FDR upper-bound to the false non-discovery rate.

Below, we will focus on the choice $\lambda=\alpha$\,, leading to the threshold collection
\begin{equation}
\Delta(i)={\alpha}\min\bigg((1-\alpha)\frac{i}{m-i+1}\,,1\bigg) \label{equ_onestage_alpha}\,.
\end{equation}
For $i\leq (m+1)/2$, the threshold
\eqref{equ_onestage_alpha} is
\[
\alpha\frac{(1-\alpha)i}{m-i+1}\,,
\]
and thus our approach differs from the threshold collection of 
the standard LSU procedure threshold by the factor $\frac{(1-\alpha)m}{m-i+1}$. 

It is interesting to note that the correction factor
$\frac{m}{m-i+1}$ appears in Holm's step-down procedure \citep{Holm1979}
for FWER control. 
The latter is a well-known improvement of Bonferroni's procedure
(which corresponds to the fixed threshold $\alpha/m$),
taking into account the proportion of true nulls,
and defined as the step-down procedure\footnote{The step-down procedure
with threshold collection $\Delta$ rejects the hypotheses corresponding to the $k$ smallest $p$-values,
where $k=\max\{ 0\leq i \leq m \telque \forall j \leq i\,\;\; p_{(j)}\leq
\Delta(j) \}$. It is self-consistent with respect to $\Delta$ but 
uniformly more conservative than the step-up
procedure with the same threshold collection, compare with Definition \ref{def_stepup}.} with
threshold collection $\alpha/(m-i+1)$\,. Here we therefore prove that
this correction is suitable as well for the linear step-up procedure, 
in the framework of FDR control.

If $r$ denotes the final number of rejections of the new one-stage
procedure, 
we can interpret the ratio $\frac{(1-\lambda)m}{m-r+1}$
between the adaptive threshold and the LSU threshold at the same point
as an {\em a posteriori} estimate for $\pi_0^{-1}$\,.
In the next section we propose to use this quantity in a plug-in,
2-stage adaptive procedure.

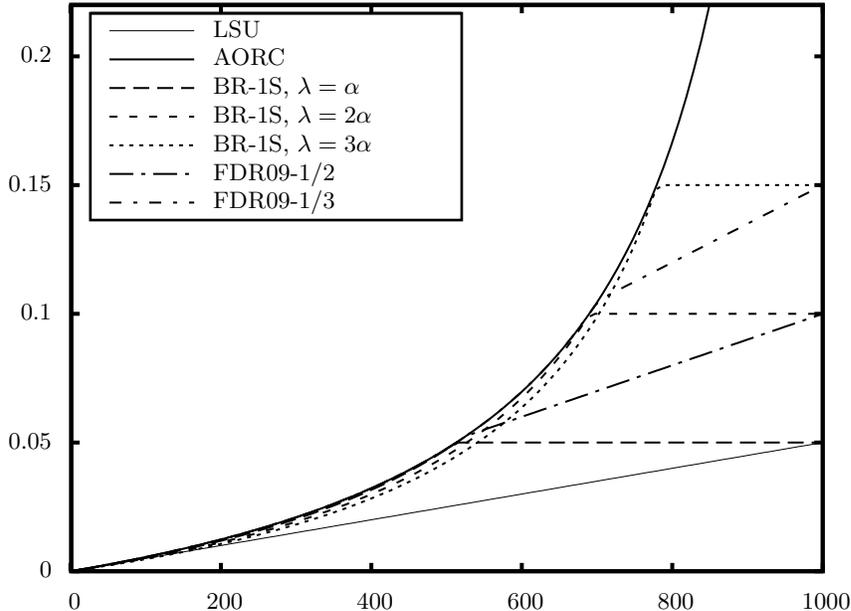
\begin{figure}[h!]
\begin{center}
\input{compar_glob.tex}
\end{center}
\vspace{-.5cm}
\caption{\label{fig_stepupadapt}For $m=1000$ null hypotheses and $\alpha=5\%$: comparison of 
 the new threshold collection \textit{BR-1S-$\lambda$} given by \eqref{equ_onestage_lambda} to that of
the LSU, the AORC and \textit{FDR09-$\eta$}\,.}
\end{figure}

As Figure \ref{fig_stepupadapt} illustrates, 
our procedure is generally less conservative than the
(non-adaptive) linear step-up procedure (LSU). 
Precisely, the new procedure can only be more conservative than the LSU
procedure in the marginal case where the factor 
$\frac{(1-\alpha)m}{m-i+1}$ is smaller than one. This happens only
when the proportion of null
hypotheses rejected by the LSU procedure is positive but less than
$\alpha+1/m$\, (and even in this region the ratio of the two threshold
collections is never less than $(1-\alpha)$\,). Roughly speaking, this
situation with only few rejections can only happen if there are few false hypotheses to
begin with ($\pi_0$ close to 1) or if the false hypotheses are very
difficult to detect (the distribution of false $p$-values is close
to being uniform).

In the interest of being more specific, we briefly investigate this issue in the next
lemma, considering the particular {\em Gaussian random effects} model (which is relatively standard in 
the multiple testing literature, see e.g. \citealp{GW2004}) in order to
give a quantitative answer from an asymptotical point of view (when the number
of tested hypotheses grows to infinity). In the random effect model, hypotheses
are assumed to be randomly true or false with probability $\pi_0$\,, and the false
null hypotheses share a common distribution $P_1$\,. Globally, the $p$-values
then are i.i.d. drawn according to the mixture distribution $\pi_0 U[0,1] + (1-\pi_0)P_1$\,.

\begin{lemma}\label{lemma:probamarg}
Consider the random effects model where the
$p$-values are i.i.d. with common cumulative distribution function
$t\mapsto\pi_0 t+(1-\pi_0)F(t)$. Assume the true null hypotheses are
standard Gaussian with zero mean and the alternative hypotheses are
standard Gaussian with mean $\mu>0$\,. In this case
$F(t)=\overline{\Phi}(\overline{\Phi}^{-1}(t)-\mu)$\,,
where $\overline{\Phi}$ is the standard Gaussian upper tail
function. Assuming $\pi_0 < (1+\alpha)^{-1}$\,, define 
\begin{equation}
 \mu^{\star}=\overline{\Phi}^{-1}(\alpha^2)-\overline{\Phi}^{-1}\bigg(\frac{\alpha^{-1}-\pi_0}{1-\pi_0}
 \alpha^2\bigg).\nonumber 
\end{equation}
Then if $\mu>\mu^*$\,, the probability that the LSU rejects a proportion of null
hypotheses less than $1/m+\alpha$ tends to 0 as $m$ tends to infinity.
On the other hand, if $\pi_0 > (1+\alpha)^{-1}$\,, or $\mu<\mu^*$\,, then this probability tends to one.
\end{lemma}
For instance, taking in the above lemma the values $\pi_0=0.5$ and $\alpha=0.05$, 
results in the critical value $\mu^{\star}\simeq 1.51$\,. This lemma delineates
clearly in a particular case in which situation we can expect an improvement from the
adaptive procedure over the standard LSU.


\subsubsection{Comparison to other adaptive one-stage procedures.}
\label{FiDiRo}

Very recently, other adaptive one-stage procedures with important similarities 
to {\em BR-1S-$\lambda$} have been proposed by other authors. (The
present work was developed independently.)

Starting with some 
heuristic motivations,  \cite{FDR2008}  proposed the threshold collection
$t(i) = \frac{\alpha i}{m-(1-\alpha)i}$, which they dubbed the {\em
  asymptotically optimal rejection curve} (AORC). However, the step-up
procedure using this threshold collection as is
does not have controlled FDR (since $t(m)=1$\,, the
corresponding step-up procedure would always reject all the
hypotheses), and several suitable modifications were proposed by \cite{FDR2008},
the simplest one being $$t'_{\eta}(i) = \min\big(t(i), \eta^{-1} \alpha i/m\big),$$
which is denoted by \textit{FDR09-$\eta$} in the following.

  The theoretical FDR control proved in \cite{FDR2008} is studied 
asymptotically as the number of hypotheses grows to infinity.
 In that framework, asymptotical control at level $\alpha$ is shown to
 hold for any $\eta<1$.
On Figure \ref{fig_stepupadapt}, we represented the thresholds 
\textit{BR-1S-$\lambda$} and \textit{FDR09-$\eta$} for comparison, for
several choices of the parameters. The two families
appear quite similar, initially following the AORC curve then
branching out or capping at a point depending on the parameter. One noticeable
difference in the initial part of the curve is that while \textit{FDR09-$\eta$} exactly coincides with
the AORC, \textit{BR-1S-$\lambda$} is arguably
sligthly more conservative. This reflects the nature of the
corresponding theoretical result -- non-asymptotic control of the FDR
requires a somewhat more conservative threshold as compared to the
%
only asymptotic control of \textit{FDR-$\eta$}\,. Moreover,
 we can use \textit{BR-1S-$\lambda$} as a first step 
in a 2-step procedure, as will be argued in the next section.

The ratio between \textit{BR-1S-$\lambda$} and the AORC (before
the capping point) is a factor which, assuming $\alpha\geq
(m+1)^{-1}$\,, is lower bounded by $(1-\lambda)(1-\frac{1}{m+1})$\,.
This suggests that the value for $\lambda$ should be kept small,
this is why we propose $\lambda=\alpha$\, as a default choice.

Finally, the {\em step-down} procedure based on the same threshold
collection $t(i)$ (without modification) is proposed and studied by \cite{GBS2008}.
Using specific properties of step-down procedures, 
these authors proved the nonasymptotic FDR control of this procedure.




\subsection{Adaptive plug-in methods
}\label{sec:adaptplugin}

In this section, we consider different adaptive step-up procedures of
the plug-in type, i.e. based on an explicit estimator of $\pi_0^{-1}$\,. 
We first review a general method proposed by \cite{BKY2006}
in order to derive FDR control for such plug-in procedures (see also
Theorem~4.3 of \citealp{FDR2008}). 
We propose here a self-contained
proof of this result, which is somewhat more compact
than the original one (and also extends the original result from
step-up procedures to more generally self-consistent procedures).
Based on this result, we review the different
plug-in estimators considered by \cite{BKY2006} and add a new one to the lot,
based on the one-stage adaptive procedure introduced in the previous section.


Let us first introduce the following notations: for
each $h\in\cH$, we denote by $\mbf{p}_{-h}$ the collection of $p$-values $\mbf{p}$ restricted
to $\cH\setminus\set{h}$\,, that is, $\mbf{p}_{-h} = (p_{h'},h'\neq h)$\,.
We also denote $\mbf{p}_{0,h}=(\mbf{p}_{-h},0)$ the
collection $\mbf{p}$ where $p_h$ has been replaced
by $0$. 

\begin{theorem}[Benjamini, Krieger, Yekutieli 2006]\label{th_stepupadapt}
Suppose that the family $p$-values $\mathbf{p}=(p_h,h\in\cH)$ is
independent. Let $G:[0,1]^\cH\rightarrow
(0,\infty)$ be a measurable, coordinate-wise non-increasing function.
Consider a non-increasing multiple testing procedure $R$ which is self-consistent with
respect to the adaptive linear threshold collection $\Delta(\bp,i) =
\alpha G(\bp) i/m$\,. Then the following holds:
\begin{equation}
\label{adapt-control}
\FDR(R) \leq \frac{\alpha}{m} \sum_{h \in \cH_0} \e{G(\bp_{0,h})}\,.
\end{equation}
In particular, if for any $h \in \cH_0$\,, it holds that
$\e{G(\bp_{0,h})} \leq \pi_0^{-1}$\,, then $\FDR(R)\leq \alpha$\,.
\end{theorem}

We will apply the above result to the following estimators, depending
on a fixed parameter $\lambda \in (0,1)$ or $k_0 \in \set{1,\ldots,m}$:
\begin{align*}
\text{[Storey-$\lambda$]} && G_1(\mbf{p})
&=\frac{(1-\lambda)m}{\sum_{h\in\cH}\ind{p_h>\lambda}+1}\, ;\\ 
\text{[Quant-$\frac{k_0}{m}$]} && G_2(\mbf{p})
&=\frac{(1-p_{(k_0)})m}{m-k_0+1}\, ;\\ 
\text{[BKY06-$\lambda$]} && G_3(\mbf{p}) &=\frac{(1-\lambda)m}{m-|R_0(\mbf{p})|+1}\,, \text{ where } R_0
\text{ is the standard LSU at level } \lambda\,;\\
\text{[BR-2S-$\lambda$]} && G_4(\mbf{p}) &=\frac{(1-\lambda)m}{m-|R'_0(\mbf{p})|+1}\,, \text{ where
} R'_0 \text{ is BR-1S-$\lambda$ (see Theorem~\ref{th_main1})}.
\end{align*}
Above, the notations ``Storey-$\lambda$", ``Quant-$\frac{k_0}{m}$",
``BKY06-$\lambda$" and ``BR-2S-$\lambda$" refer to the plug-in
adaptive linear step-up procedures associated to $G_1$, $G_2$, $G_3$
and $G_4$, respectively.

Estimator $G_1$  is usally called \textit{modified Storey's estimator} and was
initially introduced by \cite{Storey2002} from an heuristics on the
$p$-values histogram (originally without the ``+1'' in the numerator, hence the
name ``modified''). Its intuitive justification is as follows: the set
$S_\lambda$ of $p$-values larger than the threshold $\lambda$ contains on average at
least a proportion $(1-\lambda)$ of the true null hypotheses. Hence, a
natural estimator of $\pi_0^{-1}$ is $(1-\lambda) m/|S_\lambda \cap
\cH_0| \leq (1-\lambda) m/|S_\lambda| \simeq G_1(\mbf{p})$\,. Therefore, we expect that
Storey's estimator is generally an overestimate of $\pi_0^{-1}$\,.
A standard choice is $\lambda=1/2$ (as in the SAM software of \citealp{ST2003}). 
FDR control for the corresponding plug-in step-up procedure was proved
by  \cite{STS2004} 
(actually, for the modification $\wt{\Delta}(\mbf{p},i) =
  \min(\alpha G_1(\mbf{p}) i/m,\lambda)$\,) and by
\cite{BKY2006}. 

Estimator $G_2$ was introduced by \cite{BH2000} and \cite{ETST2001}, from
a slope heuristics on the $p$-values c.d.f. Roughly speaking, $G_2$
appears as a Storey's estimator with the data-dependent choice
$\lambda=p_{(k_0)}$\,, and can therefore be interpreted as
the quantile version of the Storey estimator. A standard value for
$k_0$ is $\lfloor m/2 \rfloor$, resulting in the so-called median
adaptive LSU (see \cite{BKY2006} and the references therein).

Estimator $G_3$ was introduced by \cite{BKY2006} 
for the particular choice $\lambda=\alpha/(1+\alpha)$. More precisely, 
a slightly less conservative version, without the ``+1'' in the
denominator, was used in \cite{BKY2006}. We forget about this refinement here, noting that it
results only in a very slight improvement.

Finally, the estimator $G_4$ is new and follows exactly the same
philosophy as $G_3$, that is, uses a step-up procedure as a first
stage in order to estimate $\pi_0^{-1}$\,, but this time based on our
adaptive one-stage step-up procedure introduced in the previous section, rather than
the standard LSU. Note that since $R'_0$ is less conservative than
$R_0$ (except in marginal cases), we generally have  $G_2\leq G_3$ pointwise
and our estimator improves over the one of \cite{BKY2006}.

These different estimators all satisfy the sufficient condition
mentioned in Theorem \ref{th_stepupadapt}, and
we thus obtain the following corollary: 
\begin{corollary} 
\label{cor_stepupadapt} Assume
that the $p$-values of $\mbf{p}=(p_h,h\in\cH)$ are independent.
For $i=1,2,3,4$\,, and any $h \in \cH_0$\,, it holds that 
$\e{G_i(\bp_{0,h}}\leq \pi_0^{-1}$\,. Therefore, the
plug-in adaptive linear step-up procedure at level $\alpha$ using
estimator $G_i$ has FDR smaller than $\alpha$\,.
\end{corollary}

The above result for $G_1$, $G_2$ and $G_3$ (for
 $\lambda=\alpha/(1+\alpha)$)
was proved by \cite{BKY2006}. For completeness, we reproduce
shortly the corresponding arguments in the appendix.

In other words, Corollary~\ref{cor_stepupadapt} states that under
independence, for any $\lambda$ and $k_0$, the 
plug-in adaptive procedures  Storey-$\lambda$, Quant-$\frac{k_0}{m}$, BKY06-$\lambda$ and
BR-2S-$\lambda$ all control the FDR at level $\alpha$.

\begin{remark}\label{remark_CCadapt_improv} The result proved by
 \cite{BKY2006} is actually slightly sharper than Theorem \ref{th_stepupadapt}.
Namely, if $G(\cdot)$ is moreover supposed to be coordinate-wise
left-continuous, it is possible to prove that Theorem
  \ref{th_stepupadapt} still holds when $\bp_{0,h}$ in
the RHS of \eqref{adapt-control} is 
replaced by the slightly better 
$\wt{\mbf{p}}_h =
(\mbf{p}_{-h},\wt{p}_h(\mbf{p}_{-h}))$\,, defined as the collection of $p$-values
$\mbf{p}$ where $p_h$ has been replaced by
$\wt{p}_h(\mbf{p}_{-h})=\max\big\{p\in[0,1]\;\big|\;
p\leq\alpha\pi(h)|R(\mbf{p}_{-h},p)|G(\mbf{p}_{-h},p)\big\}.$ This
improvement then permits to get rid of the ``+1'' in the denominator of
$G_3$\,. Here, we opted for simplicity and a more straightforward
statement, noting that this improvement is not crucial.
\end{remark}

\begin{remark} \label{FiDiRo2} The one-stage step-up procedure of \cite{FDR2008} (see 
previously the discussion at the end of Section \ref{FiDiRo}) --- for which there is no result proving
non-asymptotic FDR control up to our knowledge --- can also be interpreted 
intuitively as an adaptive version of the LSU using estimator $G_2$\,,
where the choice of parameter $k_0$ is data-dependent.
Namely, assume we reject at least $i$ null hypotheses 
whenever $p_{(i)}$ is lower than the standard LSU threshold 
times the estimator $G_2$\, wherein parameter $k_0=i$\, is used. This corresponds
to the inequality $p_{(i)} \leq \frac{(1-p_{(i)}) k }{m-i+1}$\,, which,
solved in $p_{(i)}$\,, gives the threshold collection of
\cite{FDR2008}. Remember from Section \ref{FiDiRo} that this threshold 
collection must actually be modified in order to be useful, since it
otherwise always leads to reject all hypotheses. The modification
leading to \textit{FDR09-$\eta$}
consists in capping the estimated $\pi_0^{-1}$ at a level $\eta$\,,
i.e. using $\min(\eta,G_2)$ instead of $G_2$ in the above reasoning.
In fact, the proof of \cite{FDR2008} relies on a result which is
essentially a reformulation of Theorem~\ref{th_stepupadapt}
for a specific form of estimator.
\end{remark}

\begin{remark} \label{unfavorablecases} 
The estimators $G_i$, $i=1,2,3,4$ are not necessarily larger
than 1, and to this extent can in some unfavorable cases result in
the final procedure being actually more conservative than the standard
LSU. This can only happen in the situation where
either $\pi_0$ is close to 1 (``sparse signal'') or the alternative hypotheses are
difficult to detect (``weak signal''); 
if such a situation is anticipated, it is more appropriate to use the regular
non-adaptive LSU.

For the Storey-$\lambda$ estimator, we can control precisely the probability that
such an unfavorable case arises by using Hoeffding's inequality \citep{Hoeff1963}:
assuming the true nulls are i.i.d. uniform on $(0,1)$ and the false nulls i.i.d. of c.d.f. $F(\cdot)$, we write by definition of $G_1$
\begin{align}
\prob{G_1(\mathbf{p})<1} &= \prob{\frac{1}{m} \sum_{h \in \cH}^m 
(\ind{p_h>\lambda}-\prob{p_h>\lambda})>(1-\pi_0)(F(\lambda)-\lambda)-m^{-1}}\nonumber\\
& \leq  \exp(-2(mc^2+1)) \label{equ_hoeffding},
\end{align}
\end{remark}
where we denoted $c=(1-\pi_0)(F(\lambda)-\lambda)$\,, 
and assumed additionally $c>m^{-1}$\,.
The behavior of the bound mainly depends on $c$\,, which can get small
only if $\pi_0$ is close to 1 (sparse signal) or $F(\lambda)$ is close
to $\lambda$ (weak signal), illustrating the above point. In general,
provided $c>0$ does not depend on $m$\,,  the probability that the Storey procedure fails to outperform the LSU
vanishes exponentially as $m\rightarrow \infty$\,.

\subsection{Theoretical robustness of the adaptive procedures under maximal dependence}\label{sec:theexdep}

For the different procedures proposed above, the theory only provides
the correct FDR control under independence between the $p$-values. An
important issue is to know how robust this control is
when dependences are present (as it is often the case in
practice). However, the analytic computation of the FDR under
dependence is generally a difficult task, and this issue is often
tackled empirically through simulations in a pre-specified model 
(we will do so in Section~\ref{sec_simu}).

In this short section, we present theoretical computations of the FDR
for the previously introduced adaptive step-up procedures,
under the maximally dependent model where all the $p$-values
are in fact equal, that is $p_h\equiv p_1$ for all $h\in\cH$\, (and
$m_0=m$). It corresponds to the case where we perform $m$ times the
same test, with the same $p$-value. Albeit relatively trivial and
limited, this case leads to very simple FDR computations and provides
at least some hints concerning the robustness under dependence of the different
procedures studied above.



\begin{proposition}\label{prop-rho1}
Suppose that we observe $m$ identical $p$-values
$\mathbf{p}=(p_1,...,p_m)=(p_1,...,p_1)$ with $p_1\sim U([0,1])$ and
assume $m=m_0$. Then, the following holds: 
 \begin{align*}
\FDR(\textit{BR-1S-$\lambda$})&= \min(\lambda,(1-\lambda) \alpha m),\\
\FDR(\textit{FDR09-$\eta$})&= \alpha \eta^{-1}, \\
\FDR(\textit{Storey-$\lambda$})&=\min\bigg(\lambda,\alpha(1-\lambda)m\bigg)
+ \left(\alpha(1-\lambda)(1+m^{-1})-\lambda\right)_+,\\
\FDR(\textit{Quant-$k_0/m$})&= \frac{\alpha }{(1+\alpha) - (k_0 -1)m^{-1}},\\
\FDR(\textit{BKY06-$\lambda$})&=\FDR(\textit{BR-2S-$\lambda$})=\FDR(\textit{Storey-$\lambda$}).
\end{align*}
\end{proposition}

Interestingly, the above proposition suggests specific choices of the
 parameters $\lambda$, $\eta$ and $k_0$ to ensure control of the FDR
at level $\alpha$ under maximal dependence:

\begin{itemize}
\item[\textbullet] For BR-1S-$\lambda$, putting $\lambda_2=\alpha/(\alpha +  m^{-1})$, Proposition~\ref{prop-rho1} gives that
$\FDR(\textit{BR-1S-$\lambda$})=\lambda$ whenever $\lambda\leq
\lambda_2$. 
This suggests to take $\lambda=\alpha$\,, and is thus in
accordance with the default choice proposed in Section~\ref{sec:onestageadapt}.
\item[\textbullet] For {FDR09-$\eta$}, no choice of $\eta<1$ will lead
to the correct FDR control under maximal dependence. However, the larger $\eta$\,, the smaller the
FDR in this situation. Note that $\FDR(\textit{FDR09-$\frac 1 2$})=2\alpha$.
\item[\textbullet] For {Storey-$\lambda$}, {BKY06-$\lambda$} and
{BR-2S-$\lambda$}, putting $\lambda_1=\alpha / (1+\alpha + m^{-1})$,
we have
$\FDR=\lambda$ for $\lambda_1\leq \lambda\leq
\lambda_2.$ 
This suggests to choose $\lambda=\alpha$ within
these three procedures.  Furthermore, note that the standard choice
$\lambda=1/2$ for {Storey-$\lambda$} leads to a very poor control
under maximal dependence:
$\FDR(\textit{Storey-$1/2$})= {\min(\alpha m , 1 )}/{2}$.
\item[\textbullet] For {Quant-$k_0/m$}, we see that the value
of $k_0$ maximizing the FDR while maintaining it below $\alpha$ is
$k_0=\lfloor\alpha m\rfloor +1$. Remark also that the standard
choice $k_0=\lfloor m/2 \rfloor$ leads to
$\FDR(\textit{Quant-$k_0/m$})={2\alpha}/{(1+2\alpha +2m^{-1})}\simeq
2\alpha$.
\end{itemize}


Nevertheless, we would like to underline that the above computations
should be interpreted with caution, as the maximal dependence
case is very specific and cannot possibly give an accurate idea of the behavior of
the different procedures when the correlation between the $p$-values
are strong but not equal to $1$\,. For instance, it is well-known that
the LSU procedure has FDR far below $\alpha$ for strong positive
correlations, but its FDR is equal
to $\alpha$ in the above extreme model
(see \citealp{FDR2007}, for a comprehensive study of the LSU under
positive dependence). Conversely, the FDR of some adaptive
procedures can be higher under moderate dependence than under maximal
dependence. This behavior appears in the simulations of the next section,
illustrating the complexity of the issue.




\subsection{Simulation study}\label{sec_simu}

How can we compare the different adaptive procedures defined above?
For a fixed $\lambda$, we have pointwise $G_1\geq G_4
\geq G_3$\,, which shows that the adaptive procedure [Storey-$\lambda$]
is always less conservative than [BR-2S-$\lambda$], itself
less conservative than [BKY06-$\lambda$] (except in the marginal
cases where the one-stage adaptive procedure is more conservative than
the standard step-up procedure, as delineated earlier for example in
Lemma \ref{lemma:probamarg}).  It would
therefore appear that one should always choose [Storey-$\lambda$] and disregard the
other ones.  However, a important point made by \cite{BKY2006} for
introducing $G_3$ as a better alternative to the (already known
earlier) $G_1$ is that, on simulations with positively dependent test
statistics, 
the plug-in procedure using $G_1$ with $\lambda = 1/2$ had very poor
control of the FDR, while the FDR was still controlled for the plug-in
procedure based on $G_3$. While the
positively dependent case is not covered by the theory, it is of course
very important to ensure that a multiple testing procedure is sufficiently
robust in practice so that the FDR does not vary too much in this
situation.

In order to assess the quality of our new procedures, we
compare here the different methods on a simulation study
following the setting used by \cite{BKY2006}. Let
$X_i=\mu_i+\varepsilon_i,$ for $i,1\leq i\leq m$, where $\varepsilon$
is a $\mathbb{R}^m$-valued centred Gaussian random vector such that
$\E(\varepsilon^2_i)=1$ and for $i\neq j$,
$\E(\varepsilon_i\varepsilon_j)=\rho$, where $\rho\in[0,1]$ is a
correlation parameter. Thus, when $\rho=0$ the $X_i$'s are
independent, whereas when $\rho>0$ the $X_i$'s are positively
correlated (with a constant pairwise correlation).
For instance, the
$\varepsilon_i$'s can be constructed by taking
$\varepsilon_i:=\sqrt{\rho}\:\:U+\sqrt{1-\rho}\:\:Z_i,$ where $Z_i$, $1\leq
i\leq m$ and $U$ are all i.i.d $\sim \mathcal{N}(0,1)$.

Considering the one-sided null hypotheses $h_i$ : ``$\mu_i\leq 0$"
against the alternatives ``$\mu_i>0$" for $1\leq i\leq m$, we define
the $p$-values $p_i=\overline{\Phi}(X_i)$, for $1\leq i\leq m$, where
$\overline{\Phi}$ is the standard Gaussian distribution tail. We choose
a common mean $\bar{\mu}$ for all false hypotheses, that is, 
for $i,1\leq i\leq m_0$, $\mu_i=0$ and for $i,m_0+1\leq i\leq m$,
$\mu_i=\bar{\mu}$\,;  the $p$-values corresponding to the null
means follow exactly a uniform distribution. 

Note that the case $\rho=1$ and $m=m_0$ (i.e. $\pi_0=1$) corresponds to the maximally dependent case studied in Section \ref{sec:theexdep}.

We compare the following step-up multiple testing procedures: first, the one-stage step-up procedures defined in Section~\ref{sec:onestageadapt}:
\begin{itemize}
\item[-][\textit{BR08-1S}-$\alpha$] The new procedure of Theorem \ref{th_main1}, 
with parameter $\lambda=\alpha$\,,
\item[-][\textit{FDR09-$\frac{1}{2}$}] The procedure proposed in \cite{FDR2008} and
described in Section \ref{FiDiRo}, with $\eta=\frac{1}{2}$\,.
\end{itemize}
Second, the adaptive plug-in step-up procedure defined in Section~\ref{sec:adaptplugin}:
\begin{itemize}
\item[-][\textit{Median LSU}] The procedure [Quant-$\frac{k_0}{m}$] with the choice $\frac{k_0}{m} = \frac{1}{2}$\,,
\item[-][\textit{BKY06-$\alpha$}] The procedure [\textit{BKY06-$\lambda$}] with 
the parameter choice $\lambda=\alpha$\,,
 \item[-][\textit{BR08-2S}-$\alpha$]  The procedure [\textit{BR08-2S}-$\lambda$] with 
the parameter choice $\lambda=\alpha$\,,
\item[-][\textit{Storey-$\lambda$}] With the choices $\lambda=1/2$ and $\lambda=\alpha$\,.
\end{itemize}
Finally, we used as oracle reference [\textit{LSU Oracle}], the step-up procedure with the threshold collection $\Delta(i)=\alpha i/m_0$, using a perfect estimation of $\pi_0$. 

The parameter choice $\lambda=\alpha$ for [\textit{Storey-$\lambda$}] comes from the relationship of $G_3,G_4$
to $G_1$ in Section~\ref{sec:onestageadapt}, and form the discussion of the maximally dependent
case in Section \ref{sec:theexdep}.
Note that the procedure studied by \cite{BKY2006} is actually
[BKY06-$\alpha/(1+\alpha)$] in our notation (up to the very slight modification explained in Remark~\ref{remark_CCadapt_improv}). This means that the procedure [\textit{BKY06-$\alpha$}] used in our
simulations is not exactly the same as in \cite{BKY2006}, but it is very close.

The three most important parameters in the simulation are the correlation coefficient $\rho$, the
proportion of true null hypotheses $\pi_0$, and the alternative mean $\bar{\mu}$ which represents
the signal-to-noise ration, or how easy it is to distinguish alternative hypotheses.
We present in Figures \ref{fig_pi0}, \ref{fig_moy}, and \ref{fig_rho} results of the simulations
for one varying parameter ($\pi_0$, $\bar{\mu}$ and $\rho$, respectively), the others being kept fixed.
Reported are, for the different methods: the  average FDR, and the average power relative to the reference 
[LSU-Oracle]. The absolute power is defined as the average proportion of false null hypotheses rejected,
and the  relative power as the mean of the number of true rejections of the
procedure divided by the number of true rejections of [LSU-Oracle]. 
Each point is an average of $10^5$ simulations, with fixed parameters $m=100$ and $\alpha=5\%$\,.

\subsubsection{Under independence ($\rho=0$)}

Remember that under independence of the $p$-values, 
the \textit{LSU}
procedure has FDR equal to $\alpha\pi_0$ and that the \textit{LSU
Oracle} procedure has  FDR equal to $\alpha$ (provided that
$\alpha\leq \pi_0$). The other procedures have their FDR upper bounded by
$\alpha$\, (in an asymptotical sense only for [FDR09-$\frac{1}{2}$]).

\begin{figure}[p]
\begin{center}
\begin{tabular}{cc}
 FDR & Average relative power to [LSU-oracle] \\
\hspace{-1cm} \scalebox{0.8}{\input{fig_rho_1_fdr.tex}} & \hspace{-1cm}
\scalebox{0.8}{\input{fig_rho_1_puis.tex}} \\
\hspace{-1cm} \scalebox{0.8}{\input{fig_rho_2_fdr.tex}} &\hspace{-1cm}
\scalebox{0.8}{\input{fig_rho_2_puis.tex}} \\
\hspace{-1cm} \scalebox{0.8}{\input{fig_rho_3_fdr.tex}} & \hspace{-1cm}
\scalebox{0.8}{\input{fig_rho_3_puis.tex}} \\
$\rho$ & $\rho$ \\
\end{tabular}
\end{center}
\caption{ \label{fig_rho} FDR and power relative to oracle
as a function of the pairwise correlation coefficient $\rho$\,.Target FDR is $\alpha=5\%$\,, 
total number of hypotheses $m=100$\,.  
The mean for the alternatives is $\bar{\mu}=3$.
From top to bottom: proportion of true null hypotheses $\pi_0 \in \set{0.2, 0.5, 0.8}$.}
\end{figure}
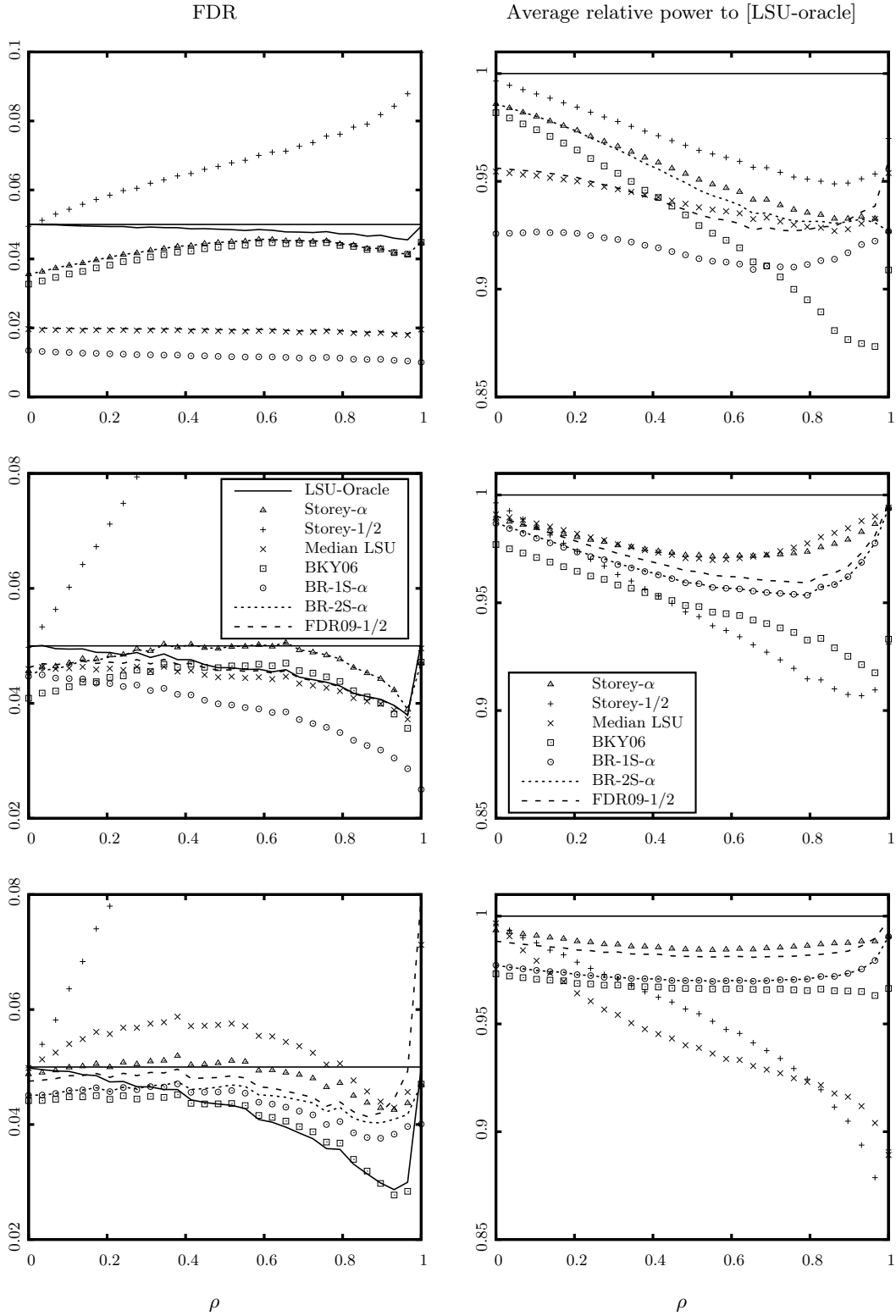

The situation where the $p$-values are independent corresponds to the first
row of Figures \ref{fig_pi0} and \ref{fig_moy} and the leftmost point of
each graph in Figure \ref{fig_rho}.
It appears that in the independent case, 
the following procedures can be consistently ordered in terms of (relative)
power over the range of parameters studied here:
\begin{center}
[\textit{Storey-$1/2$}] $\succ$ [\textit{Storey-$\alpha$}] $\succ$
[\textit{BR08-2S-$\alpha$}] $\succ$ [\textit{BKY06-$\alpha$}], 
\end{center}
the symbol ``$\succ$" meaning ``is (uniformly over our experiments) more powerful
than". 

Next, the procedures [median-LSU] and [FDR09-$\frac{1}{2}$] appear both
consistently less powerful than [Storey-$\frac{1}{2}$], and
[FDR09-$\frac{1}{2}$] is additionally also consistently less powerful than [Storey-$\alpha$].
Their relation to the remaining procedures depends on the parameters;
both [median-LSU] and [FDR09-$\frac{1}{2}$] appear to be more powerful
than the remaining procedures when $\pi_0>\frac{1}{2}$, and less
efficient otherwise. We note that [median-LSU] also appears to perform better
when $\bar{\mu}$ is low (i.e. the alternative hypotheses are harder to
distinguish).

Concerning our one-stage procedure [BR08-1S-$\alpha$], we note that it
appears to be indistinguishable from its two-stage counterpart
[BR08-2S-$\alpha$] when $\pi_0>\frac{1}{2}$\,, and significantly less
powerful otherwise. This also corresponds to our expectations, since
in the situation $\pi_0<\frac{1}{2}$\,, there is a much higher
likelihood that more than 50\% hypotheses are rejected, in which case
our one-stage threshold family hits its ``cap'' at level $\alpha$\,
(see e.g. Fig. \ref{fig_stepupadapt}; a similar qualitative explanation
applies to understand the behavior of $FDR09-1/2$). This is precisely to improve on this situation
that we introduced the 2-stage procedure, and we see that does 
in fact improve substantially the 1-stage version in that specific region.

The fact that [Storey-$\frac{1}{2}$] is uniformly more powerful than
the other procedures in the independent case corroborates the
simulations reported in \cite{BKY2006}. Generally speaking, under
independence we obtain a less biased estimate for $\pi_0^{-1}$
when considering Storey's estimator based on a ``high'' threshold like
$\lambda=\frac{1}{2}$\,. Namely, higher $p$-values are less likely to be
``contaminated'' by false null hypotheses; conversely, if we take a lower
threshold $\lambda$, there will be more false null hypotheses included in the
set of $p$-values larger than $\lambda$\,, leading to a pessimistic
bias in the estimation of $\pi_0^{-1}$\,. This qualitative reasoning
is also consistent with the observed behavior of [median-LSU], since
the set of $p$-values larger than the median is much more likely to be
``contaminated'' when $\pi_0<\frac{1}{2}$\,.

However, the problem with [Storey-$\frac{1}{2}$] is that the
corresponding estimation of $\pi_0^{-1}$ exhibits much more
variability than its competitors when there is a substantial
correlation between the $p$-values. As a consequence it is
a very fragile procedure. This phenomenon was already pinpointed
in \cite{BKY2006} and we study it next.

\begin{figure}[p]
\begin{center}
\begin{tabular}{cc}
 FDR & Average relative power to [LSU-oracle] \\
\hspace{-1cm} \scalebox{0.8}{\input{fig_pi0_1_fdr.tex}} & \hspace{-1cm}
\scalebox{0.8}{\input{fig_pi0_1_puis.tex}} \\
\hspace{-1cm} \scalebox{0.8}{\input{fig_pi0_2_fdr.tex}} & \hspace{-1cm}
\scalebox{0.8}{\input{fig_pi0_2_puis.tex}} \\
\hspace{-1cm} \scalebox{0.8}{\input{fig_pi0_3_fdr.tex}} & \hspace{-1cm}
\scalebox{0.8}{\input{fig_pi0_3_puis.tex}} \\
$\pi_0$ & $\pi_0$
\end{tabular}
\end{center}
\caption{ \label{fig_pi0} FDR and power relative to oracle 
as a function of the true proportion $\pi_0$ of null hypotheses\,. 
Target FDR is $\alpha=5\%$\,, total number of hypotheses $m=100$\,. 
The mean for the alternatives is $\bar{\mu}=3$. From top to bottom:
pairwise correlation coefficient $\rho\in\set{0, 0.2, 0.5}$.}
\end{figure}
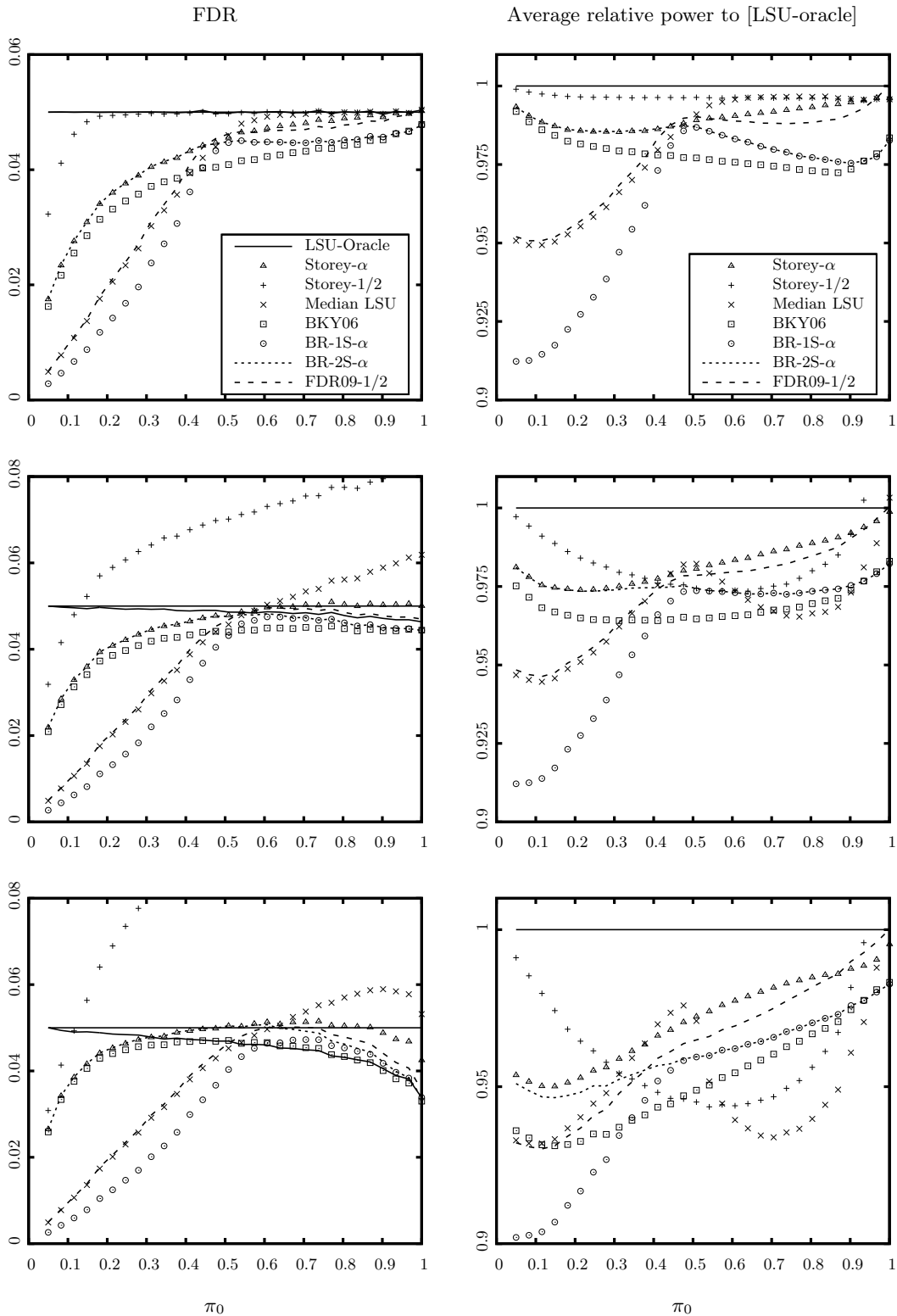

\begin{figure}[p]
\begin{center}
\begin{tabular}{cc}
 FDR & Average relative power to [LSU-oracle] \\
\hspace{-1cm} \scalebox{0.8}{\input{fig_moy_1_fdr.tex}}&  \hspace{-1cm}
\scalebox{0.8}{\input{fig_moy_1_puis.tex}} \\
\hspace{-1cm} \scalebox{0.8}{\input{fig_moy_2_fdr.tex}} &\hspace{-1cm}
\scalebox{0.8}{\input{fig_moy_2_puis.tex}} \\
\hspace{-1cm} \scalebox{0.8}{\input{fig_moy_3_fdr.tex}} &\hspace{-1cm}
\scalebox{0.8}{\input{fig_moy_3_puis.tex}} \\
$\bar{\mu}$ & $\bar{\mu}$
\end{tabular}
\end{center}
\caption{ \label{fig_moy} FDR and power relative to oracle
as a function of the common alternative hypothesis mean $\bar{\mu}$\,. Target FDR is $\alpha=5\%$\,,
total number of hypotheses $m=100$\,. 
The proportion of true null hypotheses
is $\pi_0=0.5$. From top to bottom:
pairwise correlation coefficient $\rho\in\set{0, 0.2, 0.5}$.}
\end{figure}
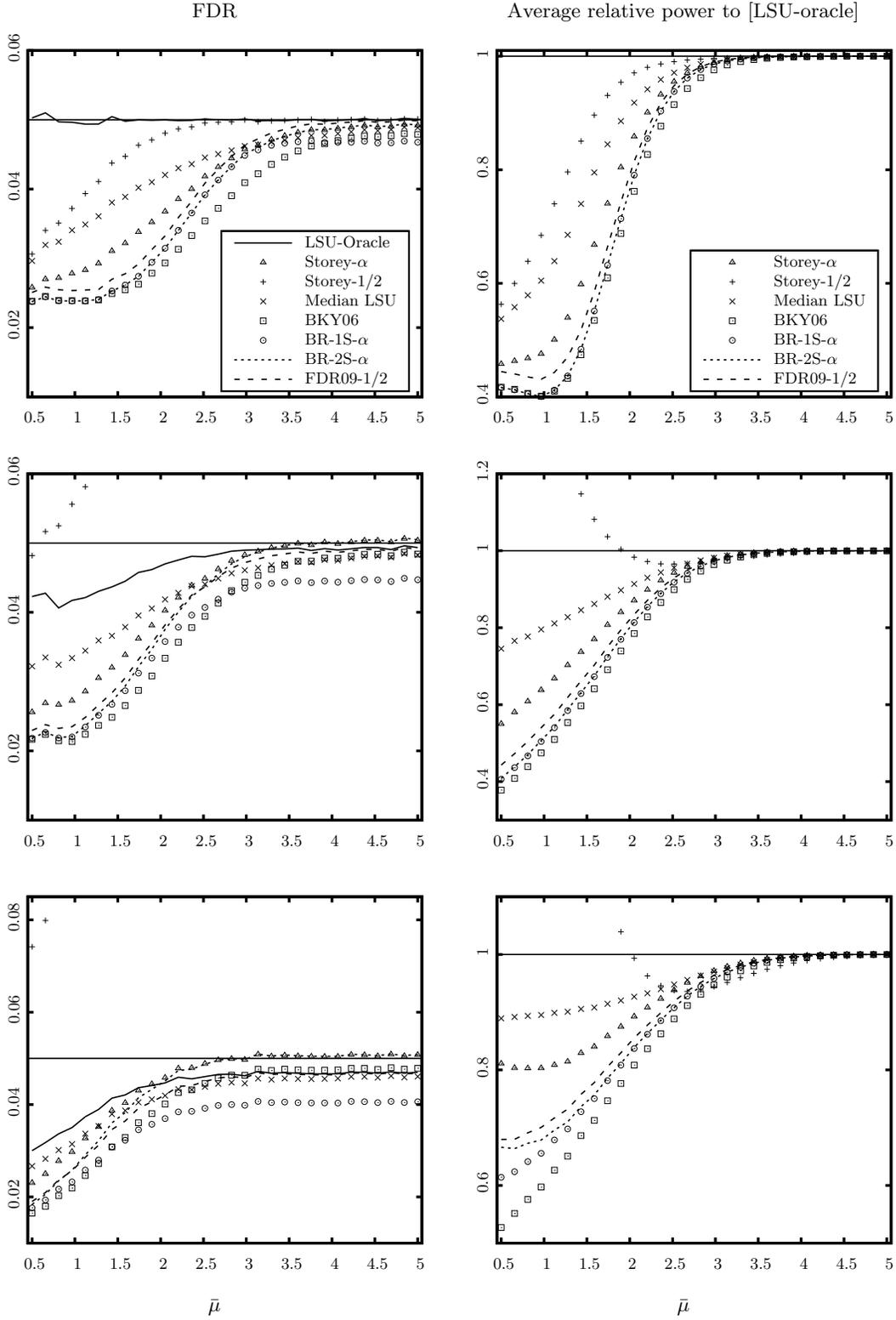

\subsubsection{Under positive dependences ($\rho>0$)}

Under positive dependences, remember that it is known theoretically from \cite{BY2001} 
that the FDR of the procedure
\textit{LSU} (resp. \textit{LSU Oracle}) is still bounded by
$\alpha\pi_0$ (resp. $\alpha$), but without equality. 
However, we do not
know from a theoretical point of view 
if the adaptive procedures have their FDR upper bounded by $\alpha$.
In fact, it was pointed out by \citet{Far2007}, in another work reporting simulations on 
adaptive procedures, that one crucial point for these seems
to be the variability of estimate of $\pi_0^{-1}$. Estimates of this quantity
that are not robust with respect to positive dependence will result in
failures for the corresponding multiple testing procedure.

The situation where the $p$-values are positively dependent corresponds to the
second and third 
rows ($\rho=0.2, 0.5$\,, respectively) of Figures \ref{fig_pi0} and \ref{fig_moy} 
and to all the graphs of Figure \ref{fig_rho}.

The most striking fact is that [Storey-$\frac{1}{2}$] does not control the
FDR at the desired level any longer under positive dependences, an can even
be off by quite a large factor. This is in accordance with the experimental findings of 
\cite{BKY2006}. Therefore, although this procedure was the favorite in the
independent case, it turns out to be not robust, which is very undesirable
for practical use where it is generally impossible to guarantee that the
$p$-values are independent. The procedure [median-LSU] appears to have
higher power than the remaining ones in the situations studied in Figure
\ref{fig_moy}, especially with a low signal-to-noise ratio. Unfortunately,
other situations appearing in Figures \ref{fig_pi0} and \ref{fig_rho} show
that [median-LSU] can exhibit a poor FDR control in some
parameter regions, most notably when $\pi_0$ is close to 1 and
positive dependences are present (see, e.g., Fig. \ref{fig_rho}, bottom row).
In a majority of practical situations, it is difficult to rule out
{\em a priori} that $\pi_0$ could be close to 1 (i.e., there is only a
small proporion of false hypotheses), or that there are no dependences.
For these reasons, our conclusion is that [median-LSU]
is also not robust enough in general to be reliable.

The other remaining procedures seem to still have a controlled FDR, or at least
to be very close to the FDR target level (except for [FDR09-$\frac{1}{2}$] when $\rho$ and $\pi_0$ are close to $1$). For these it seems that
the qualitative conclusions concerning power
comparison found in the independent case remain true.
To sum up:
\begin{itemize}
\item the best overall procedure seems to be [Storey-$\alpha$]: its FDR seems to be
under or only slightly over the target level in all situations, and it exhibits 
globally a power superior to other procedures.
\item then come in order of power, our 2-stage procedure [BR08-2S-$\alpha$], then
[BKY06-$\alpha$].
\item like in the dependent case, [FDR09-$\frac{1}{2}$] ranks second when $\pi_0>\frac{1}{2}$
but tends to perform noticeably poorer if $\pi_0$ gets smaller. Its
FDR is also not controlled if very strong correlations are present.
\end{itemize}

To conclude, the recommendation that we draw from these experiments is
that for practical use, we recommend in priority [Storey-$\alpha$], 
then as close seconds [BR08-2S-$\alpha$]
or [FDR09-$\frac{1}{2}$] (the latter when it is expected that
$\pi_0>1/2$\,, and that there are no very strong correlations
present). The procedudre [BKY06-$\alpha$] is also competitive but
appears to be in most cases noticeably outperformed by the above ones.
These procedures all exhibit good robustness to dependence for
FDR control as well as comparatively good power. The fact that [Storey-$\alpha$]
performs so well and seems to hold the favorite position
has up to our knowledge not been reported before (it was not
included in the simulations of \citealp{BKY2006}) and came somewhat as a surprise to us.

\begin{remark} As pointed out earlier, the fact that [FDR09-$\frac{1}{2}$] performs
  sub-optimally for $\pi_0<\frac{1}{2}$ appears to be strongly linked
  to the choice of parameter $\eta=\frac{1}{2}$\,. Namely, 
the implicit estimator of $\pi_0^{-1}$ in the procedure is
capped at $\eta$\, (see Remark \ref{FiDiRo2}). Choosing a higher value for
$\eta$ will reduce the sub-optimality region but increase the
variability of the estimate and thus decrease the overall
robustness of the procedure (if dependences are present; and also
under independence if only a small number of
hypotheses $m$ are tested, as for this procedure the convergence of the FDR towards
its asymptotically  controlled value becomes slower as $\eta$ grows towards 1).
\end{remark}

\begin{remark} Another 2-stage adaptive procedure was introduced in \cite{Sar2006},
which is very similar to a plug-in procedure using [Storey-$\lambda$]. In fact,
in the experiments presented in \cite{Sar2006}, the two procedures are almost equivalent,
corresponding to $\lambda=0.995$\,. We decided not to include this additional
procedure in our simulations to avoid overloading the plots. Qualitatively,
we observed that the procedures of \cite{Sar2006} or [Storey-0.995] are very similar in behavior
to [Storey-$\frac{1}{2}$]: very performant in the independent case but very fragile
with respect to deviations from independence.
\end{remark}

\begin{remark}
One could formulate the concern that the observed FDR control for
[Storey-$\alpha$] could possibly fail with other parameters settings,
for example when $\pi_0$ and/or $\rho$ are
close to one. We performed additional simulations to this regard,
which we summarize briefly here. We considered the following cases:
$\pi_0=0.95$
and varying $\rho\in[0,1]$\,; $\rho=0.95$ and varying
$\pi\in[0,1]$\,; finally $(\pi_0,\rho)$ varying both in
$[0.8,1]^2$\,, using a finer discretization grid to cover this region
in more detail.
In all the above cases Storey-$\alpha$ still had its FDR very close to (or below)
$\alpha$.  
Note also that the case $\rho \simeq 1$ and $\pi_0\simeq 1$ is in accordance with
the result of Section~\ref{sec:theexdep}, stating that
$\FDR(\textit{Storey-$\alpha$})= \alpha$ when $\rho=1$ and
$\pi_0=1$\,. Finally, we also performed additional experiments for
different choices of the number of hypotheses to test ($m=20$ and
$m=10^4$) and different choices of the target level
($\alpha=10\%,1\%$). In all of these cases were the results
qualitatively in accordance with the ones already presented here.
\end{remark}



\section{New adaptive procedures with provable FDR control under arbitrary dependence}\label{newsu_df}

In this section, we consider from a theoretical point of view the problem of constructing multiple testing procedures
that are adaptive to $\pi_0$ under arbitrary dependence conditions of the $p$-values.
The derivation of adaptive procedures that have provably controlled FDR under 
dependences appears to have been only studied scarcely (see
\citealp{Sar2006}, and \citealp{Far2007}). Here, we propose to use
a 2-stage procedure where the first stage is a multiple testing
with either controlled FWER or controlled FDR. The first option is relatively
straightfoward and is intended as a reference.
In the second case, we use Markov's inequality to estimate $\pi_0^{-1}$. Since Markov's
inequality is general but not extremely precise, the resulting
procedures are obviously quite conservative and are arguably of a
limited practical interest. However, we will show that they still
provide an improvement, in a certain regime, with respect to
the (non-adaptive) LSU procedure in the PRDS case and with respect to the
family of (non-adaptive) procedures proposed in Theorem \ref{th_BF} in
the arbitrary dependences case.


For the purposes of this section, we first recall the formal definition for PRDS dependence of \cite{BY2001}:
\begin{definition}[PRDS condition]
Remember that a set $D\subset [0,1]^{\cH}$ is said to be
\textit{non-decreasing} if for all $x,y\in [0,1]^{\cH}$, $x\leq y$
coordinate-wise and $x\in D$ implies $y\in D$. Then, the $p$-values
$\mathbf{p}=(p_h,h\in\cH)$ are said \textit{positively regressively
dependent on each one from} $\cH_0$ (PRDS on $\cH_0$ in short) if for
any non-decreasing measurable set $D\subset [0,1]^{\cH}$ and for all $h\in\cH_0$,
$u\in[0,1]\mapsto \Proba(\mathbf{p}\in D | p_h = u)$ is
non-decreasing. 
\end{definition}
On the one hand, it was proved by \cite{BY2001} that the LSU still has
controlled FDR at level $\pi_0 \alpha$ (i.e., Theorem \ref{th_BH}
still holds) under the PRDS assumption.  On the other hand, under
totally arbitrary dependences this result does not hold,
and  Theorem \ref{th_BF} provides a family of threshold collection 
resulting in controlled FDR at the same level in this case.

Our first result concerns a two-stage procedure where the first stage $R_0$
is any multiple testing procedure with controlled FWER, and where we
(over-) estimate $m_0$ via the straightforward estimator $(m-|R_0|)$\,.
This should be considered as a form of baseline reference for this type of two-stage procedure.
\begin{theorem}
\label{th_gen_holm}
Let $R_0$ be a non-increasing multiple testing procedure and assume that its FWER is
controlled at level $\alpha_0$\,, that is, $\prob{R_0 \cap \cH_0 \neq
\emptyset} \leq \alpha_0$\,. 
 Then the adaptive step-up procedure $R$
with data-dependent threshold collection
$\Delta(i)=\alpha_1 (m-|R_0|)^{-1}\beta(i)$ has FDR controlled at level
$\alpha_0 + \alpha_1$\, in
either of the
following dependence situations:
\begin{itemize}
\item the $p$-values $(p_h,h\in\cH)$ are PRDS on $\cH_0$ and the shape function is the identity function. 
\item the $p$-values have unspecified dependences and $\beta$ is a shape function of the form (\ref{equ_beta}).
\end{itemize}
\end{theorem}
Here it is clear that the price for adaptivity is a certain loss in FDR control for being able to
use the information of the first stage. If we choose $\alpha_0=\alpha_1=\alpha/2$\,, then this procedure will outperform
its non-adaptive counterpart (using the same shape function) only if there are more than $50\%$\,, rejected hypotheses in the first stage. Only
if it is expected that this situation will occur does it make sense to employ this procedure, since it will otherwise
perform worse than the non-adaptive procedure.

Our second result is a two-stage procedure where the first stage has controlled FDR.
First introduce, for a fixed constant
$\kappa\geq 2$\,, the following function: for $x\in[0,1]$, 
\begin{equation}\label{equ_Fc}
F_\kappa(x) =
\begin{cases}
1 & \text{ if } x\leq \kappa^{-1}\;\\
\frac{2\kappa^{-1}}{1-\sqrt{1-4(1-x)\kappa^{-1}}} & \text{ otherwise}
\end{cases}\,\,.
\end{equation}
If $R_0$ denotes the first stage, we propose using $F_\kappa(|R_0|)$ as an (under-)estimation 
of $\pi_0^{-1}$ at the second stage. We obtain the following result:

\begin{theorem}\label{th_gen_df}
Let $\beta$ be a fixed shape function, and $\alpha_0,\alpha_1 \in (0,1)$ such that $\alpha_0\leq \alpha_1$. 
Denote by $R_0$ the step-up procedure with threshold collection
$\Delta_0(i)=\alpha_0 \beta(i)/m$. Then the adaptive step-up procedure $R$
with data-dependent threshold collection $\Delta_1(i)=\alpha_1 \beta(i) F_\kappa(|R_0|/m) /m$ 
has FDR upper bounded by $\alpha_1+\kappa\alpha_0$ in either of the
following dependence situations:
\begin{itemize}
\item the $p$-values $(p_h,h\in\cH)$ are PRDS on $\cH_0$ and the shape function is the identity function. 
\item the $p$-values have unspecified dependences and $\beta$ is a shape function of the form (\ref{equ_beta}).
\end{itemize}
\end{theorem}




For instance, in the PRDS case, the procedure $R$ of Theorem \ref{th_gen_df}
with $\kappa=2$, $\alpha_0=\alpha/4$ and $\alpha_1=\alpha/2$, is the
adaptive linear step-up procedure at level $\alpha/2$ with the
following estimator for $\pi_0^{-1}$\,:
\[\frac{1}{1-\sqrt{(2|R_0|/m-1)_+}}\,,\] 
where $|R_0|$ is
the number of rejections of the LSU procedure at level $\alpha/4$ and
$(\cdot)_+$ denotes the positive part.

 Whether in the PRDS or arbitrary dependences case, with the above
 choice of parameters, we note that $R$ is less conservative than the
 non-adaptive step-up procedure with threshold collection $\Delta(i)=\alpha
 \beta(i)/m$ if $F_2(\abs{R_0}/\abs{\cH}) \geq 2$ or equivalently
 when $R_0$ rejects more than $F_2^{-1}(2) = 62,5\%$ of the null
 hypotheses. Conversely, $R$ is more conservative otherwise, and we
 can lose up to a factor $2$ in the threshold collection with respect
 to the standard one-stage version.  Therefore, here again this
 adaptive procedure is only useful in the cases where it is expected
 that a ``large'' proportion of null hypotheses can easily be
 rejected.  In particular, when we use Theorem \ref{th_gen_df} in
 the distribution-free case, it is relevant to choose the shape
 function $\beta$ from a prior distribution $\nu$ concentrated on the
 large numbers of $\{1,\dots,m\}$. Finally, note that it is not immediate
to see if this procedure will improve on the one of Theorem \ref{th_gen_holm}.
Namely, with the above choice parameters, it has to reject more hypotheses in
the first step than the procedure of Theorem \ref{th_gen_holm} in order
to beat the LSU, and the first step is performed at a smaller target level.
However, since the first step only controls the FDR, and not the FWER, it
can actually be much less conservative.



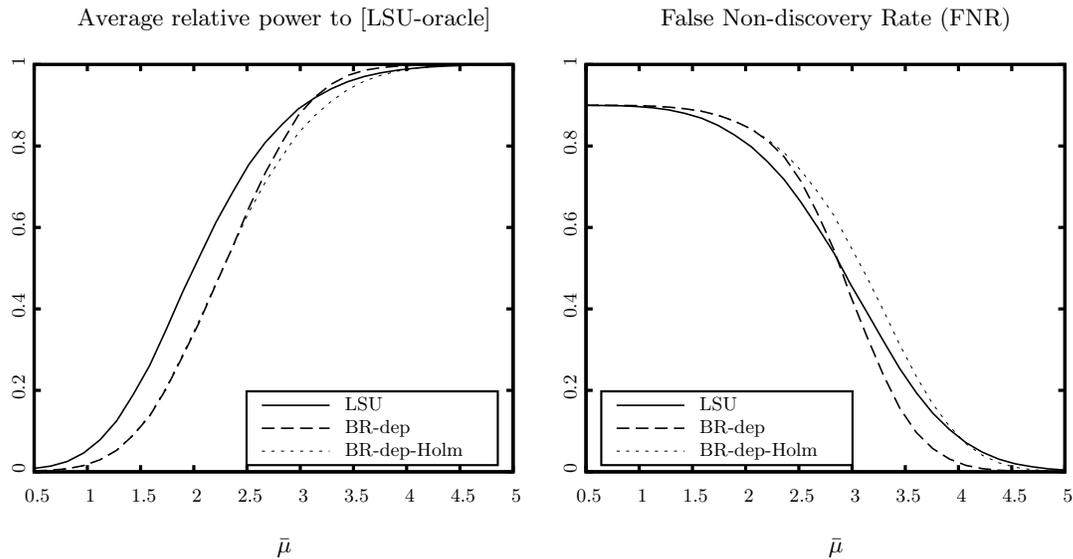
\begin{figure}[h!]
\begin{center}
\begin{tabular}{cc}
Average relative power to [LSU-oracle] & False Non-discovery Rate (FNR)\\
\hspace{-1cm} \scalebox{0.8}{\input{fig_dep_puis.tex}} & \hspace{-1cm}
\scalebox{0.8}{\input{fig_dep_fnr.tex}} \\
$\bar{\mu}$ & $\bar{\mu}$
\end{tabular}
\end{center}
\caption{\label{fig_adaptcorrpos} Relative power to oracle
and false non-discovery rate (FNR)
of the different procedures, as a function of the common alternative
hypothesis mean $\bar{\mu}$\,.
Parameters are $\alpha=5\%$\,, $m=1000$\,, $\pi_0=10\%$\,, $\rho=0.1$\,.
``BR08-dep-Holm'' corresponds to the procedure of Theorem \ref{th_gen_holm} using $\alpha_1=\alpha_0=\alpha/2$
and Holm's step-down
for the first step, and ``BR08-dep'' to the procedure of Theorem \ref{th_gen_df}
with $\kappa=2$,
$\alpha_0=\alpha/4$ and $\alpha_1=\alpha/2$\,. The shape function $\beta$ is the identity function.
Each point is an average over $10^4$ independent repetitions.
}
\end{figure}

To explore this issue, we performed limited experiments in a favorable situation to test the two above procedures,
i.e. with a small $\pi_0$\,.
Namely, we considered the simulation setting of Section~\ref{sec_simu} with $\rho=0.1$,
$m_0=100$ and $m=1000$ (hence $\pi_0=10\%$) and $\alpha=5\%$\,.
The common value $\bar{\mu}$ of the positive means varies in the range $[0,5]$\,.
Larger values of $\bar{\mu}$ correspond to a very large proportion of hypotheses that
are easy to reject, which favors the first stage of the two above procedures.
A positively correlated family of Gaussians satisfies the PRDS assumption (see \cite{BY2001})\,,
so that we use the identity shape function (linear step-up), and compare our procedures
against the standard LSU.
For the FWER-controlled  first stage of Theorem \ref{th_gen_holm}, we chose a standard
Holm procedure \cite{Holm1979}, which is a step-down procedure with threshold family $t(i) = \alpha m/(m-i+1)$\,.
In Figure \ref{fig_adaptcorrpos},  we report the average relative
power to the oracle LSU, and the False Non-discovery Rate (FNR), which is the converse
of the FDR for type II errors, i.e., the average of the ratio of non-rejected false hypotheses
over the total number of non-rejected hypotheses. Since we are in a situation where $\pi_0$ is
small, the FNR might actually be a more relevant criterion than the
raw power: in this situation, because of the small number of non
rejected hypotheses, two different procedures could have their power very similar
and close to 1, but noticeably different FNRs.

The conclusion is that there exists an (unfortunately relatively small) region where the adaptive procedures
improve over the standard LSU in terms of power. In terms of FNR, the improvement is more noticeable and
over a larger region. 
Finally, our two-step adaptive procedure of Theorem \ref{th_gen_df} appears
to outperform consistently the baseline of Theorem \ref{th_gen_holm}. These results are still unsatisfying to the extent
that the adaptive procedure improves over the non-adaptive one only in a region limited to
some quite particular cases, and underperforms otherwise. Nevertheless, this 
demonstrates
theoretically the possibility of provably adaptive procedures under dependence. Again, this theme appears
to have been theoretically studied in only a handful of previous works until now, and improving
significantly the theory in this setting is still an open challenge.

\begin{remark}
Some theoretical results for two-stage procedures under possible
dependences using a first stage with controlled FWER or controlled
FDR appeared earlier \citep{Far2007}. However, it appears that in this
reference, it is implicitly assumed that the two stages are actually independent,
because the proof relies on a conditioning argument wherein FDR control for
the second stage still holds conditionally to the first stage output. This
is the case for example if the two stages are performed on separate families of $p$-values
corresponding to a new independent observation. Here we specifically wanted to take into
account that we use the same collection of $p$-values for the two stages, and therefore
that the two stages cannot assumed to be independent. In this sense our results
are novel with respect to those of \cite{Far2007}.
\end{remark}

\begin{remark}
The theoretical problem of adaptive procedures under arbitrary
dependences was also considered by \cite{Sar2006} using two-stage
procedures. However, the procedures proposed there were reported not to
yield any significant improvement over non-adaptive procedures. In
fact, in the explicit procedures proposed by \cite{Sar2006}, it can be
seen that there exists a function $\beta$ of the form \eqref{equ_beta}
such that the second stage is always more conservative  (and sometimes
by a large factor) than the
non-adaptive step-up procedure with threshold collection $\Delta(i) =
\alpha \beta(i)/m$\,, which 
has FDR bounded by $\pi_0 \alpha$(see Theorem \ref{th_BF}).  
\end{remark}

\section{Conclusion and discussion}

We proposed several adaptive multiple testing procedures that provably
control the FDR under different hypotheses on the dependence of the $p$-values.
Firstly, we introduced the one- and two-stage procedures \textit{BR-1S} and \textit{BR-2S}
and we proved their theoretical validity
when the $p$-values are independent. The procedure
\textit{BR-2S} is less 
conservative in general (except in marginal situations) than the adaptive
procedure proposed by \cite{BKY2006}. Extensive simulations
showed that these new procedures appear to be
robustly controlling the FDR even in a positive dependence situation,
which is a very desirable property in practice.
This is an advantage with respect to the [Storey-$\frac{1}{2}$]
procedure,  which is less conservative
but breaks down under positive dependences. 
Moreover, our simulations 
showed that the choice of parameter 
$\lambda=\alpha$ instead of $\lambda=1/2$ in the Storey procedure 
resulted in  a much more
robust procedure under positive dependences, at the price of being
slightly more conservative. This fact is supported by a theoretical
investigation of the maximally dependent case. 
 These properties  do not appear to
have been reported before, and put forward Storey-$\alpha$ as
a procedure of considerable practical interest.

Secondly, we presented what we think is among the first
examples  
of adaptive multiple testing procedures with provable FDR
control in the PRDS and distribution-free cases. 
An  important
difference with respect to earlier works on this
topic is that the procedures
we introduced here are both theoretically founded and can be shown to
improve on non-adaptive procedures in certain (admittedly limited) circumstances.
Although  their interest at this point is mainly theoretical, this shows in principle that adaptivity can improve
performance in a theoretically rigorous way even without the independence assumption. 

The proofs of the results have been built upon the notion of {\em self-consistency} and other technical
tools introduced in a previous work \citep{BR2008EJS}. We believe
these tools allow for a more unified approach than in the classical adaptive multiple testing
literature, avoiding in particular to deal explicitly with the reordered $p$-values, which can
be somewhat cumbersome.

Another advantage of this approach is that it can be extended in a
relatively straightforward manner to the case of \textit{weighted FDR}, that is, the quantity
\eqref{equ_FDR} where the cardinality measure $|.|$ has been replaced by
a general measure $W(R)=\sum_{h\in R}w_h$ (with $W(\cH)=\sum_{h\in
\cH}w_h=m$). This allows in
particular to recover results very similar to those of \cite{BH2006} 
and can also be used to prove that a (generalized) Storey estimator
can be used to control the weighted FDR. 
The modifications needed to include this generalizations are relatively
minor; we omit the details here and refer the
reader to \citet{BR2008EJS} to see how the case of weighted FDR can be
handled using the same technical tools.


There remains a vast number of open issues concerning adaptive procedures.
We first want to underline once more that the theory for adaptive
procedures under dependence is still underdeveloped. It might actually
be too restrictive to look for procedures having theoretically controlled FDR
uniformly over arbitrary dependence situations such as what we studied in
Section \ref{newsu_df}. An interesting future theoretical direction
could be to prove that some of the adaptive procedures showing good
robustness in our simulations actually have controlled FDR under some
types of dependence, at least when the $p$-values are in some sense 
not too far from being independent.

\section{Proofs of the results}\label{sec_proofadaptproc}

\subsection{Proofs for Section \ref{section_indepadapt}}

The following proofs use the notations $\mbf{p}_{0,h}$ and $\mbf{p}_{-h}$ defined at the beginning of Section~\ref{sec:adaptplugin}.\\

\noindent{\bf Proof of Theorem \ref{th_main1}.}
Let $R$ denote a non-increasing self-consistent procedure
with respect to $\Delta$ defined in \eqref{equ_onestage_lambda}; by definition
$R$ satisfies
\[
R \subset \set{h\in\cH\telque p_h\leq
\min\paren{(1-\lambda)\frac{\alpha|R|}{m-|R|+1},\lambda}}.
\] 
Therefore, we have
\begin{eqnarray}
\FDR(R) &=& \sum_{h\in\cH_0}
  \e{\frac{\ind{h\in R(\mbf{p})}}{|R(\mbf{p})|}}\nonumber\\
&\leq& \sum_{h\in\cH_0}
  \e{\frac{\ind{p_h\leq (1-\lambda)\frac{\alpha|R(\mbf{p})|}{m-|R(\mbf{p})|+1}}}{|R(\mbf{p})|}}\nonumber\\
  &\leq& \sum_{h\in\cH_0}
 \e{\frac{\ind{p_h\leq (1-\lambda)\frac{\alpha|R(\mbf{p})|}{m-|R(\mbf{p}_{0,h})|+1}}}{|R(\mbf{p})|}}\nonumber\\
 &=& \sum_{h\in\cH_0}
 \e{\e{\frac{\ind{p_h\leq (1-\lambda)\frac{\alpha|R(\mbf{p})|}{m-|R(\mbf{p}_{0,h})|+1}}}{|R(\mbf{p})|}\bigg|\mbf{p}_{-h}}}\nonumber\\
&\leq& (1-\lambda)\alpha\sum_{h\in\cH_0}\e{\frac{1}{m-|R(\mbf{p}_{0,h})|+1}},\nonumber
\end{eqnarray}
The second inequality above comes from $|R(\mbf{p})|\leq |R(\mbf{p}_{0,h})|$, which itself holds because $|R|$ is
coordinate-wise non-increasing in each $p$-value.
The last inequality is obtained with Lemma~\ref{lemma_forindep} of
Section~\ref{proofs_selfbound} with $U=p_h$,
$g(U)=|R(\mbf{p}_{-h},U)|$ and
$c=\frac{(1-\lambda)\alpha}{m-|R(\mbf{p}_{0,h})|+1}$, because the
distribution of $p_h$ conditionnally to $\mbf{p}_{-h}$ is
stochastically lower bounded by a uniform distribution, $|R|$ is
coordinate-wise non-increasing and because
$\mbf{p}_{0,h}$ depends only on the $p$-values of $\mbf{p}_{-h}$. 
Finally, since the threshold collection of $R$ is upper bounded by
$\lambda$, we get $$(1-\lambda)\e{m/(m-|R(\mbf{p}_{0,h})|+1)}\leq \E G_1(\mbf{p}_{0,h}),$$ where
$G_1$ is the Storey estimator with parameter $\lambda$. We then use $\E
G_1(\mbf{p}_{0,h})\leq \pi_0^{-1}$ (see proof of Corollary
\ref{cor_stepupadapt}) to conclude.  
\hfill $\blacksquare$\par\noindent\smallskip\\

\noindent{\bf Proof of Lemma~\ref{lemma:probamarg}.}  Denote $G(t) =
\pi_0t + (1-\pi_0)F(t)$ the cdf of the $p$-values under the random
effects mixture model.
Let us denote by
$\hat{t}_m$ the threshold of the LSU procedure.
The proportion of rejected hypotheses from the initial pool is then exactly
$\wh{G}_m(\hat{t}_m)$\,, where $\wh{G}_m$ is the empirical cdf of the $p$-values.
It was  proved by \cite{GW2002} under the random effects model, that
as $m$ tends to infinity
the LSU threshold $\hat{t}_m$ converges in probability to
$t^{\star}$, which is the largest point $t\in[0,1]$ such that
$G(t)=\alpha^{-1}t$\,. Since $\wh{G}_m$ converges in probability
uniformly to $G$\,, we deduce that the proportion of rejected
hypotheses converges to $\alpha^{-1}t^*$ in probability; hence,
if $t^*>\alpha^2$\,, the probability that the proportion of rejected
hypotheses is less that $\alpha+1/m$ converges to zero; and conversely
converges to 1 if $t^*<\alpha^2$\,.

The definition of $t^*$ and the expression for $G$ in the 
Gaussian mean shift model imply the following relation whenever $t^*>0$\,:
\[
\mu=\overline{\Phi}^{-1}(t^{\star})-
\overline{\Phi}^{-1}\bigg(\frac{\alpha^{-1}-\pi_0}{1-\pi_0}t^{\star}\bigg).
\]
It is easily seen that if $\pi_0 < (1+\alpha)^{-1}$\,, the quantity $\mu^*$ in
the statement of the lemma is well defined and
 we have $t^*>\alpha^2$ for $\mu>\mu^*$. This gives the first part of the result.

Conversely, if $\pi_0 > (1+\alpha)^{-1}$
we have $t^*=0$\,, and if $\pi_0 < (1+\alpha)^{-1}$ but $\mu < \mu^*$\,, we have $t^*<\alpha^2$\,;
this leads to the second part of the result. 
\hfill
$\blacksquare$\par\noindent\smallskip\\

\noindent{\bf Proof of Theorem \ref{th_stepupadapt}.}
By definition of self-consistency, the procedure $R$ satisfies
\begin{equation}
R\subset\{h\in\cH\telque p_h\leq \alpha|R|G(\mbf{p})/m \}.\nonumber 
\end{equation}
Therefore, 
\[
\FDR(R) = \sum_{h\in\cH_0}
  \e{\frac{\ind{h\in R(\mbf{p})}}{|R(\mbf{p})|}}
\leq \sum_{h\in\cH_0} \e{\frac{\ind{p_h\leq \alpha|R(\mbf{p})|G(\mbf{p})/m}}{|R(\mbf{p})|}}\,.
\]
Since $G$ is non-increasing, we get:
\begin{eqnarray}
\FDR(R)  &\leq& \sum_{h\in\cH_0}
  \e{\frac{\ind{p_h\leq \alpha|R(\mbf{p})|G(\mbf{p}_{0,h})/m}}{|R(\mbf{p})|}}\nonumber\\
 &=& \sum_{h\in\cH_0}
  \e{\e{\frac{\ind{p_h\leq \alpha|R(\mbf{p})|G(\mbf{p}_{0,h})/m}}{|R(\mbf{p})|}\bigg| \mbf{p}_{-h}}}
  \leq \frac{\alpha }{m}\sum_{h\in\cH_0} \E G(\mbf{p}_{0,h}).\nonumber
\end{eqnarray}
The last step is obtained with Lemma \ref{lemma_forindep} of
Section~\ref{proofs_selfbound} with $U=p_h$,
$g(U)=|R(\mbf{p}_{-h},U)|$ and $c=\alpha G(\mbf{p}_{0,h})/m$, because
the distribution of $p_h$ conditionnally to $\mbf{p}_{-h}$ is
stochastically lower bounded by a uniform distribution, $|R|$ is
coordinate-wise non-increasing and $\mbf{p}_{0,h}$ depends only on the
$p$-values of $\mbf{p}_{-h}$. 
\hfill $\blacksquare$\par\noindent\smallskip\\

\noindent{\bf Proof of Corollary \ref{cor_stepupadapt}.}
We prove that the
  sufficient condition of  Theorem \ref{th_stepupadapt}
holds for the nonincreasing estimators $G_i$, $i=1,2,3,4$. 
To that end, we reproduce here without major changes the arguments  used by \cite{BKY2006}.
The bound for $G_1$ is obtained using Lemma \ref{lemma_estpi0indep} (see below) with
  $k=m_0$ and $q=1-\lambda$: for all $h\in\cH_0$,
$$\e{G_1({\mbf{p}}_{0,h})}\leq
m(1-\lambda)\E\bigg[\bigg(\sum_{h'\in\cH_{0}\backslash
  \{h\}}\ind{p_{h'}>\lambda}+1\bigg)^{-1}\bigg]\leq \pi_0^{-1}.$$ The proof
for $G_3$ and $G_4$ is deduced from the one of $G_1$ because $G_3\leq G_4\leq G_1$
pointwise.

Let us prove that $\E G_2(\mbf{p}_{0,h})\leq \pi_0^{-1}$, for any
$h\in\cH_0$ and any $k_0\in\{1,...,m\}$. 
If $k_0\leq m_1+1$, the result is
trivial. Suppose now $k_0>m_1+1$. Introduce the following auxiliary notation: 
for $\bp$ a family of $p$-values indexed by 
$\cH$\,, and a subset $B \subset \cH$\,, denote by $S(i,\bp,B)$
the $i$-th ordered $p$-value of the subfamily $(x_h')_{h'\in B}$\,.
Pointwise, $G_2$ can be rewritten as: 
\begin{align*}
 G_2(\mbf{p}_{0,h})&=\frac{m}{m+1-k_0}\bigg(1-S(k_0,\mbf{p}_{0,h},\cH)\bigg)\\
&=\frac{m}{m+1-k_0}\bigg(1-S(k_0-1, \bp, \cH\setminus\set{h}) \bigg)\\
&\leq \frac{m}{m+1-k_0}\bigg(1-S(k_0-1-m + m_0, \bp, \cH_0\setminus\set{h})\bigg),
\end{align*}
the latter coming from the relation $S(i,\bp,A) \geq
S(i-|A\setminus B|,\bp,B)$, for every finite sets $B\subsetneq A$ and integer $i> |A\setminus B|$\,.
Therefore, using that $m_0-1$ independent random
variables with marginal distributions stochastically lower bounded by
a uniform law have a $j$-largest value on average larger than $j/m_0$,
we obtain: 
\begin{align*}
 \E G_2(\mbf{p}_{0,h})&\leq \frac{m}{m+1-k_0}\bigg(1-\frac{k_0-1-m + m_0}{m_0}\bigg)= \pi_0^{-1}.
\end{align*}
\hfill $\blacksquare$\par\noindent\smallskip\\


\noindent{\bf Proof of Proposition~\ref{prop-rho1}.}
Let us first consider adaptive one-stage procedures: for any step-up procedure $R$ of threshold $\Delta(i)=\alpha\beta(i)/m$ we easily derive that the probability that $R$ makes any rejection is 
$$\prob{\exists i \telque p_i \leq \Delta(i)} =  \prob{\exists i \telque p_1 \leq \Delta(i)}=  \prob{p_1 \leq \Delta(m)}=\Delta(m),$$ which is  $\FDR(R)$ because $m_0=m$. The results for \textit{BR-1S-$\lambda$} and \textit{FDR09-$\eta$} follow.

With the same reasoning, we find that for any plug-in adaptive linear step-up procedure $R$ that uses an estimator $G(\mathbf{p})$,
\begin{equation}\label{equ_casrho1}
\FDR(R)=\prob{p_1\leq \alpha G(\mathbf{p})}.
\end{equation}
Next, for the Storey plug-in procedure, we have $G_1(p_1,...,p_1)=(1-\lambda) m/ (m\ind{p_1>\lambda}+1)$, so that applying \eqref{equ_casrho1}, we get
\begin{align*}
\FDR(\textit{Storey-$\lambda$})&=\prob{p_1\leq \alpha G_1(\mathbf{p})}\\
&= \prob{p_1\leq \lambda, \, p_1 \leq \alpha (1-\lambda) m }+\prob{p_1>\lambda, \, p_1 \leq \alpha (1-\lambda) m/(m+1) }\\
&= \min\bigg(\lambda,\alpha(1-\lambda)m\bigg) + \left(\frac{\alpha(1-\lambda)m}{m+1}-\lambda\right)_+.
\end{align*}
 For the quantile procedure, we have 
 $$\prob{p_1\leq \alpha (1-p_1) m / (m-k_0+1)} = \prob{p_1 ((1+\alpha)m-k_0+1)\leq \alpha m} = \frac{\alpha}{1+\alpha - (k_0-1)/m} .$$
 For the BKY06 procedure, we simply remark that since the linear step-up procedure of level $\lambda$ rejects all the hypotheses when $p_1\leq \lambda$ and rejects no hypothesis otherwise, the estimator $G_1$ and $G_3$ are equal in this case. The proof for \textit{BR-2S-$\lambda$} is similar.
\hfill $\blacksquare$\par\noindent

\subsection{Proofs for Section \ref{newsu_df}}

We begin with a technical lemma that will be useful for proving both Theorem \ref{th_gen_holm}
and \ref{th_gen_df}. It is related to techniques previously introduced by \cite{BR2008EJS}.
\begin{lemma}
\label{lem_dep}
Assume $R$ is a multiple testing procedure 
satisfying the self-consistency condition:
\[
R \subset \set{ h \in \cH | p_h \leq \alpha G(\bp) \beta(|R|) / m}\,,
\]
where $G(\bp)$ is a data-dependent factor. Then the following
inequality holds: 
\begin{equation}
FDR(R) \leq \alpha + \e{\frac{|R\cap \cH_0|}{|R|} \ind{|R|>0} \ind{ G(\bp) > \pi_0^{-1} }}\,, \label{lemma_gen_df}
\end{equation}
under either of the following conditions:
\begin{itemize}
\item the $p$-values $(p_h,h\in\cH)$ are PRDS on $\cH_0$\,, $R$ is non-increasing 
and $\beta$ is the identity function. 
\item the $p$-values have unspecified dependences and $\beta$ is a shape function of the form (\ref{equ_beta}).
\end{itemize}
\end{lemma}

{\bf Proof.}
We have
\begin{align*}
\FDR(R) & = \e{\frac{|R \cap \cH_0|}{|R|} \ind{|R|>0} }\\
& =  \e{\frac{|R \cap \cH_0|}{|R|} \ind{|R|>0} \ind{ G \leq \pi_0^{-1} }} 
+  \e{\frac{|R \cap \cH_0|}{|R|} \ind{|R|>0} \ind{ G > \pi_0^{-1} }}  \\
&\leq \sum_{h\in\cH_0} \e{\frac{\ind{p_h\leq \alpha \beta(|R|) /m_0}}{|R|} }
+  \e{\frac{|R \cap \cH_0|}{|R|}  \ind{|R|>0}\ind{ G > \pi_0^{-1} }}\,.
\end{align*}
The desired conclusion will therefore hold if we establish that for any $h \in \cH_0$\,,
and $c>0$\,:
\[
\e{\frac{\ind{p_h\leq c\beta(|R|)}}{|R|}}\leq c\,. \label{cond_th_gen_df}
\]
In the distribution-free case, this is a direct
consequence of Lemma \ref{occamshammer} of
Section~\ref{proofs_selfbound} with $U=p_h$ and $V=\beta(|R|)$. For
the PRDS case, we note that since $|R(\mathbf{p})|$ is coordinate-wise
nonincreasing in each $p$-value, for any $v>0$, $D=\{\mathbf{z}\in
[0,1]^{\cH}\telque |R(\mathbf{z})|<v\}$ is a measurable non-decreasing set, so
that the PRDS property implies that $u\mapsto \Proba(|R|<v \telque p_h
= u)$ is non-decreasing. This implies that $u\mapsto \Proba(|R|<v
\telque p_h \leq u)$ by the following argument (see also
\citealp{Lehm1966}, cited by \citealp{BY2001}, and \citealp{BR2008EJS}): 
putting $\gamma = \prob{p_h \leq u \telque p_h \leq u'}$\,,
\begin{align*}\
\prob{ \mbf{p} \in D \;|\; p_h \leq u'} 
& = \e{ \prob{ \mbf{p} \in D \;|\; p_h} \;|\; p_h \leq u'}\\
& = \gamma \e{ \prob{\mbf{p} \in D \;|\;p_h} \;|\; p_h \leq u} 
+ (1-\gamma) \e{ \prob{\mbf{p} \in D \;|\;p_h} \;|\; u < p_h \leq u'}\\
& \geq \e{ \prob{\mbf{p} \in D \;|\;p_h} \;|\; p_h \leq u} 
= \prob{ \mbf{p} \in D \;|\; p_h \leq u}\,. 
\end{align*}
We can then apply Lemma \ref{lemma_alt} of
Section~\ref{proofs_selfbound} with $U=p_h$ and $V=|R|$.
\hfill $\blacksquare$\par\noindent\smallskip\\

\noindent{\bf Proof of Theorem \ref{th_gen_holm}.}
By definition of a step-up procedure, the two-stage procedure $R$ 
satisfies the assumption of Lemma \ref{lem_dep} for
$G(\bp) = (1-\frac{|R_0|}{m})^{-1}$\,,
where $R_0$ is the first stage with FWER controlled at level
$\alpha_0$\,. Furthermore, it is easy to check that $|R|$ is
nonincreasing as a function of each $p$-value (since $|R_0|$ is).
Then, we can apply Lemma \ref{lem_dep}, and from inequality \eqref{lemma_gen_df} we deduce
\begin{align*}
FDR(R) & \leq \alpha_1 + \e{\frac{|R\cap \cH_0|}{|R|} \ind{ 1-\frac{|R_0|}{m} < \pi_0  }}\\
& \leq \alpha_1 + \prob{R_0 \cap \cH_0 \neq \emptyset} \\
& \leq \alpha_0 + \alpha_1\,.
\end{align*}
In the case where $R_0$ rejects all hypotheses, we
assumed implicitly that the second stage also does. 
\hfill $\blacksquare$\par\noindent\smallskip\\

\noindent{\bf Proof of Theorem \ref{th_gen_df}.}
Assume $\pi_0>0$ (otherwise
the result is trivial). By definition of a step-up procedure, the two-stage procedure $R$ 
satisfies the assumption of Lemma \ref{lem_dep}
for  $G(\bp) = F_\kappa(|R_0|/m)$\,, where $R_0$ is the first stage.
Furthermore, it is easy to check that $|R|$ is
nonincreasing as a function of each $p$-value (since $|R_0|$ is).
Then, we can apply Lemma \ref{lem_dep}, and from inequality \eqref{lemma_gen_df} we deduce
\begin{align*}
\FDR(R)
& \leq  \alpha_1 + \e{\frac{|R\cap \cH_0|}{|R|} \ind{F_{\kappa}(|R_0|/m) > \pi_0^{-1}  } }\\
& \leq \alpha_1 + m_0\e{\frac{\ind{F_{\kappa}(|R_0|/m)>\pi_0^{-1}}}{|R_0|}}
\end{align*}
For the
second inequality, we have used the two following facts:

(i) $F_{\kappa}(|R_0|/m) > \pi_0^{-1}$ implies $|R_0|>0$,

(ii) because of the assumption $\alpha_0 \leq \alpha_1$ and $F_{\kappa} \geq 1$\,, the output of the
second step is necessarily a set containing at least the output of the first step.
Hence $|R| \geq |R_0|$\,.

Let us now concentrate on further bounding this second term.
For this, first consider the generalized inverse of $F_{\kappa}$\,, 
$F_{\kappa}^{-1}(t) = \inf \set{ x\telque F_{\kappa}(x) > t}$\,.
Since $F_{\kappa}$ is a non-decreasing
left-continuous function, we have $F_{\kappa}(x) > t \Leftrightarrow x > F_{\kappa}^{-1}(t)$\,.
Furthermore, the expression of $F_{\kappa}^{-1}$ is given by: $\forall t\in[1,+\infty),
 F_{\kappa}^{-1}(t) = \kappa^{-1}t^{-2} -t^{-1} +1$ (providing in particular that $F_{\kappa}^{-1}(\pi_0^{-1})>1-\pi_0$).
Hence
\begin{align}
m_0 \e{\frac{\ind{F_{\kappa}(|R_0|/m) > \pi_0^{-1}}}{|R_0|}} & \leq
m_0 \e{\frac{\ind{|R_0|/m > F_{\kappa}^{-1}(\pi_0^{-1})}}{|R_0|}} \nonumber \\
& \leq \frac{\pi_0}{F_{\kappa}^{-1}(\pi_0^{-1})} \prob{|R_0|/m \geq F_{\kappa}^{-1}(\pi_0^{-1})}\,.
\label{inter_R}
\end{align}
Now, by assumption, the FDR of the first step $R_0$ is controlled at level $\pi_0\alpha_0$\,, so that
\begin{align*}
\pi_0\alpha_0 &\geq \e{\frac{|R_0 \cap \cH_0|}{|R_0|}\ind{|R_0|>0}}\\
  &\geq \e{\frac{|R_0|+m_0-m}{|R_0|}\ind{|R_0|>0}}\\
  &=\e{[1 + (\pi_0 -1) Z^{-1}]\ind{Z>0}}\,,
\end{align*}
where we denoted by $Z$ the random variable $|R_0|/m$\,. Hence by Markov's inequality, for all $t>1-\pi_0$,
\[
\prob{Z \geq t} \leq \Proba\bigg([1 + (\pi_0 -1) Z^{-1}]\ind{Z>0}\geq 1+(\pi_0 -1) t^{-1}\bigg) \leq \frac{\pi_0 \alpha_0}{1+(\pi_0 -1) t^{-1}}\,;
\]
choosing $t = F_{\kappa}^{-1}(\pi_0^{-1})$ and using this into \eqref{inter_R}, we obtain
\[
m_0 \e{\frac{\ind{F_{\kappa}(|R_0|/m) > \pi_0^{-1}}}{|R_0|}} \leq \alpha_0 \frac{\pi_0^2}{F_{\kappa}^{-1}(\pi_0^{-1}) - 1 + \pi_0}\,.
\]
If we want this last quantity to be less than $\kappa \alpha_0$\,, this yields the condition
$F_{\kappa}^{-1}(\pi_0^{-1}) \geq \kappa^{-1} \pi_0^2 -\pi_0 +1$\,, and this is true from the  expression of $F_{\kappa}^{-1}$
(note that this is how the formula for
$F_{\kappa}$ was determined in the first place). 
\hfill $\blacksquare$\par\noindent\smallskip\\

\section{Probabilistic lemmas}\label{proofs_selfbound}

The three following lemmas 
have been established in a previous work (see \citealp{BR2008EJS}, Lemma~3.2).

\begin{lemma}\label{lemma_forindep} Let $g:[0,1]\rightarrow
  (0,\infty)$ be a non-increasing 
function. Let $U$ be a random variable which has a distribution stochastically lower bounded by a uniform
  distribution, 
that is, $ \forall u \in [0,1], \; \Proba(U\leq u) \leq u$\,. Then, for any constant $c>0$, we have
$$\E\paren{\frac{\ind{U\leq cg(U)}}{g(U)}} \leq c\,.$$
\end{lemma}

\begin{lemma}\label{lemma_alt} Let $U,V$ be two non-negative real variables. Assume the following:
\begin{enumerate}
\item The distribution of $U$ is stochastically lower bounded by a uniform
  distribution, 
that is, $ \forall u \in [0,1], \; \Proba(U\leq u) \leq u$\,.
\item The conditional distribution of $V$ given $U\leq u$ is stochastically decreasing in $u$, that is,
\[
\forall\,  v\geq 0\,\qquad \forall\, 0\leq u \leq u'\,, \qquad \Proba(V < v\telque U \leq u) \leq \Proba(V < v\telque U \leq u')\,.
\]
\end{enumerate}
Then, for any constant $c>0$, we have
\[
\E\paren{\frac{\ind{U\leq cV}}{V}} \leq c\,.
\]
\end{lemma}

\begin{lemma}\label{occamshammer}
Let $U,V$ be two non-negative real variables and $\beta$ be a function of the form (\ref{equ_beta}). Assume that the distribution of $U$ is stochastically lower bounded by a uniform
  distribution, 
that is, $ \forall u \in [0,1], \; \Proba(U\leq u) \leq u$\,. 
Then, for any constant $c>0$, we have
$$\E\paren{\frac{\ind{U\leq c\beta(V)}}{V}} \leq c\,.$$
\end{lemma}

The following lemma was stated by \cite{BKY2006}. It is a major
point when we estimate $\pi_0^{-1}$ in the independent case. The proof
is left to the reader.
\begin{lemma}\label{lemma_estpi0indep}
For any $k\geq 2$, $q\in ]0,1]$\,, let $Y$ be a binomial 
  random variable with parameters $(k-1,q)$; then the following holds:
\begin{equation}
\E[1/(1+Y)]\leq 1/kq.\nonumber
\end{equation}
\end{lemma}



\bibliographystyle{abbrvnat}
\bibliography{BR08b}

\end{document}

%% file: macrosbis.tex
\ifx \isPCversion \thiscommandisnotdefinedtralala  

\newtheorem{theorem}{Theorem}[section]
\newtheorem{proposition}[theorem]{Proposition}
\newtheorem{corollary}[theorem]{Corollary}
\newtheorem{lemma}[theorem]{Lemma}

\theoremstyle{remark}
\newtheorem{remark}[theorem]{Remark}

\theoremstyle{definition}

\newtheorem{definition}[theorem]{Definition}

\fi

\newcommand{\Proba}{\mathbb{P}}
\newcommand{\E}{\mathbb{E}} 

\newcommand{\telque}{\;|\;}
\newcommand{\FDR}{\mathrm{FDR}}

\newcommand{\mbe}{\mathbb{E}}
\newcommand{\mbr}{\mathbb{R}}
\newcommand{\mbp}{\mathbb{P}}

\newcommand{\mtc}{\mathcal}
\newcommand{\mbf}{\mathbf}

\newcommand{\wt}[1]{{\widetilde{#1}}}
\newcommand{\wh}[1]{{\widehat{#1}}}

\newcommand{\ind}[1]{{\mbf{1}\set{#1}}}
\newcommand{\e}[1]{\mbe\brac{#1}}

\newcommand{\prob}[1]{\mbp\brac{#1}}

\newcommand{\paren}[1]{\left(#1\right)}
\newcommand{\brac}[1]{\left[#1\right]}

\newcommand{\set}[1]{\left\{#1\right\}}
\newcommand{\abs}[1]{\left| #1 \right|}

\newcommand{\eps}{\varepsilon}

\newcommand{\bp}{\mbf{p}}

\newcommand{\cH}{{\mtc{H}}}



%% file: compar_glob.tex
\begingroup
  \makeatletter
  \providecommand\color[2][]{%
    \GenericError{(gnuplot) \space\space\space\@spaces}{%
      Package color not loaded in conjunction with
      terminal option `colourtext'%
    }{See the gnuplot documentation for explanation.%
    }{Either use 'blacktext' in gnuplot or load the package
      color.sty in LaTeX.}%
    \renewcommand\color[2][]{}%
  }%
  \providecommand\includegraphics[2][]{%
    \GenericError{(gnuplot) \space\space\space\@spaces}{%
      Package graphicx or graphics not loaded%
    }{See the gnuplot documentation for explanation.%
    }{The gnuplot epslatex terminal needs graphicx.sty or graphics.sty.}%
    \renewcommand\includegraphics[2][]{}%
  }%
  \providecommand\rotatebox[2]{#2}%
  \@ifundefined{ifGPcolor}{%
    \newif\ifGPcolor
    \GPcolorfalse
  }{}%
  \@ifundefined{ifGPblacktext}{%
    \newif\ifGPblacktext
    \GPblacktexttrue
  }{}%
  \let\gplgaddtomacro\g@addto@macro
  \gdef\gplbacktext{}%
  \gdef\gplfronttext{}%
  \makeatother
  \ifGPblacktext
    \def\colorrgb#1{}%
    \def\colorgray#1{}%
  \else
    \ifGPcolor
      \def\colorrgb#1{\color[rgb]{#1}}%
      \def\colorgray#1{\color[gray]{#1}}%
      \expandafter\def\csname LTw\endcsname{\color{white}}%
      \expandafter\def\csname LTb\endcsname{\color{black}}%
      \expandafter\def\csname LTa\endcsname{\color{black}}%
      \expandafter\def\csname LT0\endcsname{\color[rgb]{1,0,0}}%
      \expandafter\def\csname LT1\endcsname{\color[rgb]{0,1,0}}%
      \expandafter\def\csname LT2\endcsname{\color[rgb]{0,0,1}}%
      \expandafter\def\csname LT3\endcsname{\color[rgb]{1,0,1}}%
      \expandafter\def\csname LT4\endcsname{\color[rgb]{0,1,1}}%
      \expandafter\def\csname LT5\endcsname{\color[rgb]{1,1,0}}%
      \expandafter\def\csname LT6\endcsname{\color[rgb]{0,0,0}}%
      \expandafter\def\csname LT7\endcsname{\color[rgb]{1,0.3,0}}%
      \expandafter\def\csname LT8\endcsname{\color[rgb]{0.5,0.5,0.5}}%
    \else
      \def\colorrgb#1{\color{black}}%
      \def\colorgray#1{\color[gray]{#1}}%
      \expandafter\def\csname LTw\endcsname{\color{white}}%
      \expandafter\def\csname LTb\endcsname{\color{black}}%
      \expandafter\def\csname LTa\endcsname{\color{black}}%
      \expandafter\def\csname LT0\endcsname{\color{black}}%
      \expandafter\def\csname LT1\endcsname{\color{black}}%
      \expandafter\def\csname LT2\endcsname{\color{black}}%
      \expandafter\def\csname LT3\endcsname{\color{black}}%
      \expandafter\def\csname LT4\endcsname{\color{black}}%
      \expandafter\def\csname LT5\endcsname{\color{black}}%
      \expandafter\def\csname LT6\endcsname{\color{black}}%
      \expandafter\def\csname LT7\endcsname{\color{black}}%
      \expandafter\def\csname LT8\endcsname{\color{black}}%
    \fi
  \fi
  \setlength{\unitlength}{0.0500bp}%
  \begin{picture}(7200.00,5040.00)%
    \gplgaddtomacro\gplbacktext{%
      \csname LTb\endcsname%
      \put(984,480){\makebox(0,0)[r]{\strut{} 0}}%
      \put(984,1451){\makebox(0,0)[r]{\strut{} 0.05}}%
      \put(984,2422){\makebox(0,0)[r]{\strut{} 0.1}}%
      \put(984,3393){\makebox(0,0)[r]{\strut{} 0.15}}%
      \put(984,4364){\makebox(0,0)[r]{\strut{} 0.2}}%
      \put(1122,240){\makebox(0,0){\strut{} 0}}%
      \put(2256,240){\makebox(0,0){\strut{} 200}}%
      \put(3390,240){\makebox(0,0){\strut{} 400}}%
      \put(4524,240){\makebox(0,0){\strut{} 600}}%
      \put(5658,240){\makebox(0,0){\strut{} 800}}%
      \put(6792,240){\makebox(0,0){\strut{} 1000}}%
      \put(240,2616){\rotatebox{90}{\makebox(0,0){\strut{}}}}%
      \put(3960,-120){\makebox(0,0){\strut{}}}%
    }%
    \gplgaddtomacro\gplfronttext{%
      \csname LTb\endcsname%
      \put(2199,4557){\makebox(0,0)[l]{\strut{}LSU}}%
      \csname LTb\endcsname%
      \put(2199,4341){\makebox(0,0)[l]{\strut{}AORC}}%
      \csname LTb\endcsname%
      \put(2199,4125){\makebox(0,0)[l]{\strut{}BR-1S, $\lambda=\alpha$ }}%
      \csname LTb\endcsname%
      \put(2199,3909){\makebox(0,0)[l]{\strut{}BR-1S, $\lambda=2\alpha$ }}%
      \csname LTb\endcsname%
      \put(2199,3693){\makebox(0,0)[l]{\strut{}BR-1S, $\lambda=3\alpha$ }}%
      \csname LTb\endcsname%
      \put(2199,3477){\makebox(0,0)[l]{\strut{}FDR09-1/2}}%
      \csname LTb\endcsname%
      \put(2199,3261){\makebox(0,0)[l]{\strut{}FDR09-1/3}}%
    }%
    \gplbacktext
    \put(0,0){\includegraphics{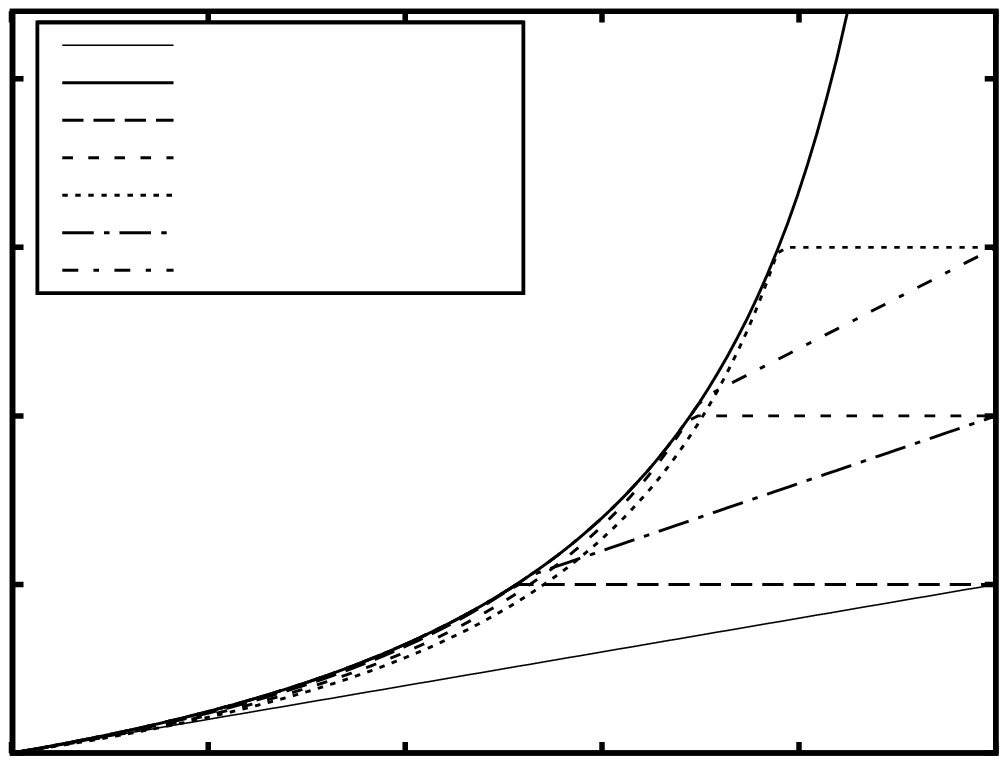}}%
    \gplfronttext
  \end{picture}%
\endgroup

%% file: fig_pi0_1_fdr.tex
\begingroup
  \makeatletter
  \providecommand\color[2][]{%
    \GenericError{(gnuplot) \space\space\space\@spaces}{%
      Package color not loaded in conjunction with
      terminal option `colourtext'%
    }{See the gnuplot documentation for explanation.%
    }{Either use 'blacktext' in gnuplot or load the package
      color.sty in LaTeX.}%
    \renewcommand\color[2][]{}%
  }%
  \providecommand\includegraphics[2][]{%
    \GenericError{(gnuplot) \space\space\space\@spaces}{%
      Package graphicx or graphics not loaded%
    }{See the gnuplot documentation for explanation.%
    }{The gnuplot epslatex terminal needs graphicx.sty or graphics.sty.}%
    \renewcommand\includegraphics[2][]{}%
  }%
  \providecommand\rotatebox[2]{#2}%
  \@ifundefined{ifGPcolor}{%
    \newif\ifGPcolor
    \GPcolorfalse
  }{}%
  \@ifundefined{ifGPblacktext}{%
    \newif\ifGPblacktext
    \GPblacktexttrue
  }{}%
  \let\gplgaddtomacro\g@addto@macro
  \gdef\gplbacktext{}%
  \gdef\gplfronttext{}%
  \makeatother
  \ifGPblacktext
    \def\colorrgb#1{}%
    \def\colorgray#1{}%
  \else
    \ifGPcolor
      \def\colorrgb#1{\color[rgb]{#1}}%
      \def\colorgray#1{\color[gray]{#1}}%
      \expandafter\def\csname LTw\endcsname{\color{white}}%
      \expandafter\def\csname LTb\endcsname{\color{black}}%
      \expandafter\def\csname LTa\endcsname{\color{black}}%
      \expandafter\def\csname LT0\endcsname{\color[rgb]{1,0,0}}%
      \expandafter\def\csname LT1\endcsname{\color[rgb]{0,1,0}}%
      \expandafter\def\csname LT2\endcsname{\color[rgb]{0,0,1}}%
      \expandafter\def\csname LT3\endcsname{\color[rgb]{1,0,1}}%
      \expandafter\def\csname LT4\endcsname{\color[rgb]{0,1,1}}%
      \expandafter\def\csname LT5\endcsname{\color[rgb]{1,1,0}}%
      \expandafter\def\csname LT6\endcsname{\color[rgb]{0,0,0}}%
      \expandafter\def\csname LT7\endcsname{\color[rgb]{1,0.3,0}}%
      \expandafter\def\csname LT8\endcsname{\color[rgb]{0.5,0.5,0.5}}%
    \else
      \def\colorrgb#1{\color{black}}%
      \def\colorgray#1{\color[gray]{#1}}%
      \expandafter\def\csname LTw\endcsname{\color{white}}%
      \expandafter\def\csname LTb\endcsname{\color{black}}%
      \expandafter\def\csname LTa\endcsname{\color{black}}%
      \expandafter\def\csname LT0\endcsname{\color{black}}%
      \expandafter\def\csname LT1\endcsname{\color{black}}%
      \expandafter\def\csname LT2\endcsname{\color{black}}%
      \expandafter\def\csname LT3\endcsname{\color{black}}%
      \expandafter\def\csname LT4\endcsname{\color{black}}%
      \expandafter\def\csname LT5\endcsname{\color{black}}%
      \expandafter\def\csname LT6\endcsname{\color{black}}%
      \expandafter\def\csname LT7\endcsname{\color{black}}%
      \expandafter\def\csname LT8\endcsname{\color{black}}%
    \fi
  \fi
  \setlength{\unitlength}{0.0500bp}%
  \begin{picture}(5760.00,4608.00)%
    \gplgaddtomacro\gplbacktext{%
      \csname LTb\endcsname%
      \put(986,480){\rotatebox{90}{\makebox(0,0){\strut{} 0}}}%
      \put(986,1760){\rotatebox{90}{\makebox(0,0){\strut{} 0.02}}}%
      \put(986,3040){\rotatebox{90}{\makebox(0,0){\strut{} 0.04}}}%
      \put(986,4320){\rotatebox{90}{\makebox(0,0){\strut{} 0.06}}}%
      \put(1128,240){\makebox(0,0){\strut{} 0}}%
      \put(1565,240){\makebox(0,0){\strut{} 0.1}}%
      \put(2002,240){\makebox(0,0){\strut{} 0.2}}%
      \put(2438,240){\makebox(0,0){\strut{} 0.3}}%
      \put(2875,240){\makebox(0,0){\strut{} 0.4}}%
      \put(3312,240){\makebox(0,0){\strut{} 0.5}}%
      \put(3749,240){\makebox(0,0){\strut{} 0.6}}%
      \put(4186,240){\makebox(0,0){\strut{} 0.7}}%
      \put(4622,240){\makebox(0,0){\strut{} 0.8}}%
      \put(5059,240){\makebox(0,0){\strut{} 0.9}}%
      \put(5496,240){\makebox(0,0){\strut{} 1}}%
      \put(528,2400){\rotatebox{90}{\makebox(0,0){\strut{}}}}%
      \put(3312,-72){\makebox(0,0){\strut{}}}%
    }%
    \gplgaddtomacro\gplfronttext{%
      \csname LTb\endcsname%
      \put(4200,2187){\makebox(0,0)[l]{\strut{}LSU-Oracle}}%
      \csname LTb\endcsname%
      \put(4200,1971){\makebox(0,0)[l]{\strut{}Storey-$\alpha$}}%
      \csname LTb\endcsname%
      \put(4200,1755){\makebox(0,0)[l]{\strut{}Storey-1/2}}%
      \csname LTb\endcsname%
      \put(4200,1539){\makebox(0,0)[l]{\strut{}Median LSU}}%
      \csname LTb\endcsname%
      \put(4200,1323){\makebox(0,0)[l]{\strut{}BKY06}}%
      \csname LTb\endcsname%
      \put(4200,1107){\makebox(0,0)[l]{\strut{}BR-1S-$\alpha$}}%
      \csname LTb\endcsname%
      \put(4200,891){\makebox(0,0)[l]{\strut{}BR-2S-$\alpha$}}%
      \csname LTb\endcsname%
      \put(4200,675){\makebox(0,0)[l]{\strut{}FDR09-1/2}}%
    }%
    \gplbacktext
    \put(0,0){\includegraphics{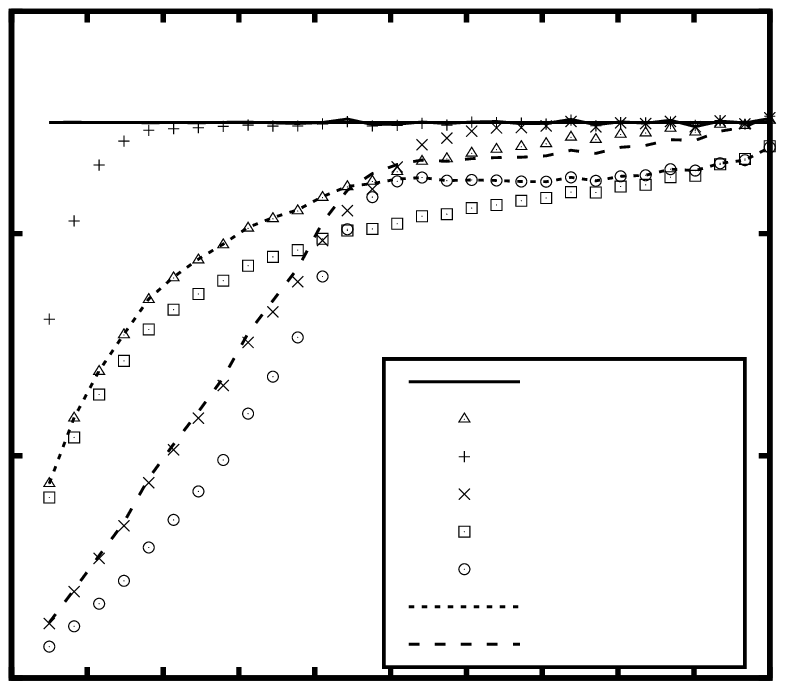}}%
    \gplfronttext
  \end{picture}%
\endgroup

%% file: fig_pi0_1_puis.tex
\begingroup
  \makeatletter
  \providecommand\color[2][]{%
    \GenericError{(gnuplot) \space\space\space\@spaces}{%
      Package color not loaded in conjunction with
      terminal option `colourtext'%
    }{See the gnuplot documentation for explanation.%
    }{Either use 'blacktext' in gnuplot or load the package
      color.sty in LaTeX.}%
    \renewcommand\color[2][]{}%
  }%
  \providecommand\includegraphics[2][]{%
    \GenericError{(gnuplot) \space\space\space\@spaces}{%
      Package graphicx or graphics not loaded%
    }{See the gnuplot documentation for explanation.%
    }{The gnuplot epslatex terminal needs graphicx.sty or graphics.sty.}%
    \renewcommand\includegraphics[2][]{}%
  }%
  \providecommand\rotatebox[2]{#2}%
  \@ifundefined{ifGPcolor}{%
    \newif\ifGPcolor
    \GPcolorfalse
  }{}%
  \@ifundefined{ifGPblacktext}{%
    \newif\ifGPblacktext
    \GPblacktexttrue
  }{}%
  \let\gplgaddtomacro\g@addto@macro
  \gdef\gplbacktext{}%
  \gdef\gplfronttext{}%
  \makeatother
  \ifGPblacktext
    \def\colorrgb#1{}%
    \def\colorgray#1{}%
  \else
    \ifGPcolor
      \def\colorrgb#1{\color[rgb]{#1}}%
      \def\colorgray#1{\color[gray]{#1}}%
      \expandafter\def\csname LTw\endcsname{\color{white}}%
      \expandafter\def\csname LTb\endcsname{\color{black}}%
      \expandafter\def\csname LTa\endcsname{\color{black}}%
      \expandafter\def\csname LT0\endcsname{\color[rgb]{1,0,0}}%
      \expandafter\def\csname LT1\endcsname{\color[rgb]{0,1,0}}%
      \expandafter\def\csname LT2\endcsname{\color[rgb]{0,0,1}}%
      \expandafter\def\csname LT3\endcsname{\color[rgb]{1,0,1}}%
      \expandafter\def\csname LT4\endcsname{\color[rgb]{0,1,1}}%
      \expandafter\def\csname LT5\endcsname{\color[rgb]{1,1,0}}%
      \expandafter\def\csname LT6\endcsname{\color[rgb]{0,0,0}}%
      \expandafter\def\csname LT7\endcsname{\color[rgb]{1,0.3,0}}%
      \expandafter\def\csname LT8\endcsname{\color[rgb]{0.5,0.5,0.5}}%
    \else
      \def\colorrgb#1{\color{black}}%
      \def\colorgray#1{\color[gray]{#1}}%
      \expandafter\def\csname LTw\endcsname{\color{white}}%
      \expandafter\def\csname LTb\endcsname{\color{black}}%
      \expandafter\def\csname LTa\endcsname{\color{black}}%
      \expandafter\def\csname LT0\endcsname{\color{black}}%
      \expandafter\def\csname LT1\endcsname{\color{black}}%
      \expandafter\def\csname LT2\endcsname{\color{black}}%
      \expandafter\def\csname LT3\endcsname{\color{black}}%
      \expandafter\def\csname LT4\endcsname{\color{black}}%
      \expandafter\def\csname LT5\endcsname{\color{black}}%
      \expandafter\def\csname LT6\endcsname{\color{black}}%
      \expandafter\def\csname LT7\endcsname{\color{black}}%
      \expandafter\def\csname LT8\endcsname{\color{black}}%
    \fi
  \fi
  \setlength{\unitlength}{0.0500bp}%
  \begin{picture}(5760.00,4608.00)%
    \gplgaddtomacro\gplbacktext{%
      \csname LTb\endcsname%
      \put(986,480){\rotatebox{90}{\makebox(0,0){\strut{} 0.9}}}%
      \put(986,1353){\rotatebox{90}{\makebox(0,0){\strut{} 0.925}}}%
      \put(986,2225){\rotatebox{90}{\makebox(0,0){\strut{} 0.95}}}%
      \put(986,3098){\rotatebox{90}{\makebox(0,0){\strut{} 0.975}}}%
      \put(986,3971){\rotatebox{90}{\makebox(0,0){\strut{} 1}}}%
      \put(1128,240){\makebox(0,0){\strut{} 0}}%
      \put(1565,240){\makebox(0,0){\strut{} 0.1}}%
      \put(2002,240){\makebox(0,0){\strut{} 0.2}}%
      \put(2438,240){\makebox(0,0){\strut{} 0.3}}%
      \put(2875,240){\makebox(0,0){\strut{} 0.4}}%
      \put(3312,240){\makebox(0,0){\strut{} 0.5}}%
      \put(3749,240){\makebox(0,0){\strut{} 0.6}}%
      \put(4186,240){\makebox(0,0){\strut{} 0.7}}%
      \put(4622,240){\makebox(0,0){\strut{} 0.8}}%
      \put(5059,240){\makebox(0,0){\strut{} 0.9}}%
      \put(5496,240){\makebox(0,0){\strut{} 1}}%
      \put(528,2400){\rotatebox{90}{\makebox(0,0){\strut{}}}}%
      \put(3312,-72){\makebox(0,0){\strut{}}}%
    }%
    \gplgaddtomacro\gplfronttext{%
      \csname LTb\endcsname%
      \put(4200,1971){\makebox(0,0)[l]{\strut{}Storey-$\alpha$}}%
      \csname LTb\endcsname%
      \put(4200,1755){\makebox(0,0)[l]{\strut{}Storey-1/2}}%
      \csname LTb\endcsname%
      \put(4200,1539){\makebox(0,0)[l]{\strut{}Median LSU}}%
      \csname LTb\endcsname%
      \put(4200,1323){\makebox(0,0)[l]{\strut{}BKY06}}%
      \csname LTb\endcsname%
      \put(4200,1107){\makebox(0,0)[l]{\strut{}BR-1S-$\alpha$}}%
      \csname LTb\endcsname%
      \put(4200,891){\makebox(0,0)[l]{\strut{}BR-2S-$\alpha$}}%
      \csname LTb\endcsname%
      \put(4200,675){\makebox(0,0)[l]{\strut{}FDR09-1/2}}%
    }%
    \gplbacktext
    \put(0,0){\includegraphics{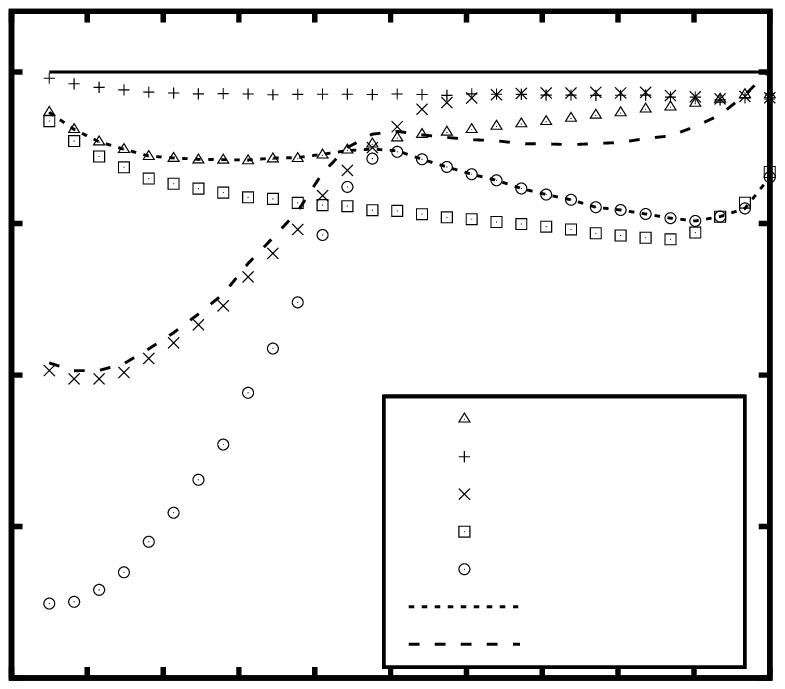}}%
    \gplfronttext
  \end{picture}%
\endgroup

%% file: fig_pi0_2_fdr.tex
\begingroup
  \makeatletter
  \providecommand\color[2][]{%
    \GenericError{(gnuplot) \space\space\space\@spaces}{%
      Package color not loaded in conjunction with
      terminal option `colourtext'%
    }{See the gnuplot documentation for explanation.%
    }{Either use 'blacktext' in gnuplot or load the package
      color.sty in LaTeX.}%
    \renewcommand\color[2][]{}%
  }%
  \providecommand\includegraphics[2][]{%
    \GenericError{(gnuplot) \space\space\space\@spaces}{%
      Package graphicx or graphics not loaded%
    }{See the gnuplot documentation for explanation.%
    }{The gnuplot epslatex terminal needs graphicx.sty or graphics.sty.}%
    \renewcommand\includegraphics[2][]{}%
  }%
  \providecommand\rotatebox[2]{#2}%
  \@ifundefined{ifGPcolor}{%
    \newif\ifGPcolor
    \GPcolorfalse
  }{}%
  \@ifundefined{ifGPblacktext}{%
    \newif\ifGPblacktext
    \GPblacktexttrue
  }{}%
  \let\gplgaddtomacro\g@addto@macro
  \gdef\gplbacktext{}%
  \gdef\gplfronttext{}%
  \makeatother
  \ifGPblacktext
    \def\colorrgb#1{}%
    \def\colorgray#1{}%
  \else
    \ifGPcolor
      \def\colorrgb#1{\color[rgb]{#1}}%
      \def\colorgray#1{\color[gray]{#1}}%
      \expandafter\def\csname LTw\endcsname{\color{white}}%
      \expandafter\def\csname LTb\endcsname{\color{black}}%
      \expandafter\def\csname LTa\endcsname{\color{black}}%
      \expandafter\def\csname LT0\endcsname{\color[rgb]{1,0,0}}%
      \expandafter\def\csname LT1\endcsname{\color[rgb]{0,1,0}}%
      \expandafter\def\csname LT2\endcsname{\color[rgb]{0,0,1}}%
      \expandafter\def\csname LT3\endcsname{\color[rgb]{1,0,1}}%
      \expandafter\def\csname LT4\endcsname{\color[rgb]{0,1,1}}%
      \expandafter\def\csname LT5\endcsname{\color[rgb]{1,1,0}}%
      \expandafter\def\csname LT6\endcsname{\color[rgb]{0,0,0}}%
      \expandafter\def\csname LT7\endcsname{\color[rgb]{1,0.3,0}}%
      \expandafter\def\csname LT8\endcsname{\color[rgb]{0.5,0.5,0.5}}%
    \else
      \def\colorrgb#1{\color{black}}%
      \def\colorgray#1{\color[gray]{#1}}%
      \expandafter\def\csname LTw\endcsname{\color{white}}%
      \expandafter\def\csname LTb\endcsname{\color{black}}%
      \expandafter\def\csname LTa\endcsname{\color{black}}%
      \expandafter\def\csname LT0\endcsname{\color{black}}%
      \expandafter\def\csname LT1\endcsname{\color{black}}%
      \expandafter\def\csname LT2\endcsname{\color{black}}%
      \expandafter\def\csname LT3\endcsname{\color{black}}%
      \expandafter\def\csname LT4\endcsname{\color{black}}%
      \expandafter\def\csname LT5\endcsname{\color{black}}%
      \expandafter\def\csname LT6\endcsname{\color{black}}%
      \expandafter\def\csname LT7\endcsname{\color{black}}%
      \expandafter\def\csname LT8\endcsname{\color{black}}%
    \fi
  \fi
  \setlength{\unitlength}{0.0500bp}%
  \begin{picture}(5760.00,4608.00)%
    \gplgaddtomacro\gplbacktext{%
      \csname LTb\endcsname%
      \put(986,480){\rotatebox{90}{\makebox(0,0){\strut{} 0}}}%
      \put(986,1440){\rotatebox{90}{\makebox(0,0){\strut{} 0.02}}}%
      \put(986,2400){\rotatebox{90}{\makebox(0,0){\strut{} 0.04}}}%
      \put(986,3360){\rotatebox{90}{\makebox(0,0){\strut{} 0.06}}}%
      \put(986,4320){\rotatebox{90}{\makebox(0,0){\strut{} 0.08}}}%
      \put(1128,240){\makebox(0,0){\strut{} 0}}%
      \put(1565,240){\makebox(0,0){\strut{} 0.1}}%
      \put(2002,240){\makebox(0,0){\strut{} 0.2}}%
      \put(2438,240){\makebox(0,0){\strut{} 0.3}}%
      \put(2875,240){\makebox(0,0){\strut{} 0.4}}%
      \put(3312,240){\makebox(0,0){\strut{} 0.5}}%
      \put(3749,240){\makebox(0,0){\strut{} 0.6}}%
      \put(4186,240){\makebox(0,0){\strut{} 0.7}}%
      \put(4622,240){\makebox(0,0){\strut{} 0.8}}%
      \put(5059,240){\makebox(0,0){\strut{} 0.9}}%
      \put(5496,240){\makebox(0,0){\strut{} 1}}%
      \put(528,2400){\rotatebox{90}{\makebox(0,0){\strut{}}}}%
      \put(3312,-72){\makebox(0,0){\strut{}}}%
    }%
    \gplgaddtomacro\gplfronttext{%
    }%
    \gplbacktext
    \put(0,0){\includegraphics{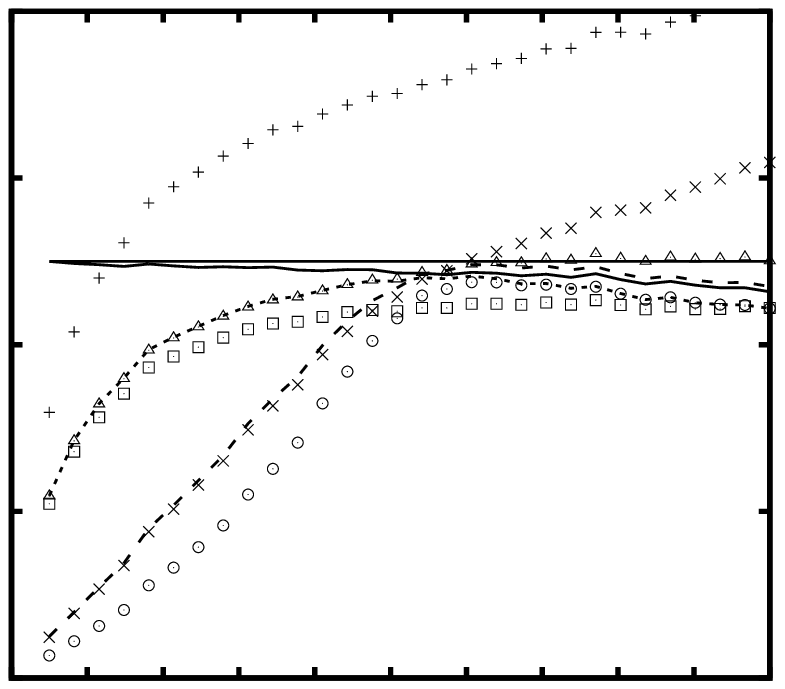}}%
    \gplfronttext
  \end{picture}%
\endgroup

%% file: fig_pi0_2_puis.tex
\begingroup
  \makeatletter
  \providecommand\color[2][]{%
    \GenericError{(gnuplot) \space\space\space\@spaces}{%
      Package color not loaded in conjunction with
      terminal option `colourtext'%
    }{See the gnuplot documentation for explanation.%
    }{Either use 'blacktext' in gnuplot or load the package
      color.sty in LaTeX.}%
    \renewcommand\color[2][]{}%
  }%
  \providecommand\includegraphics[2][]{%
    \GenericError{(gnuplot) \space\space\space\@spaces}{%
      Package graphicx or graphics not loaded%
    }{See the gnuplot documentation for explanation.%
    }{The gnuplot epslatex terminal needs graphicx.sty or graphics.sty.}%
    \renewcommand\includegraphics[2][]{}%
  }%
  \providecommand\rotatebox[2]{#2}%
  \@ifundefined{ifGPcolor}{%
    \newif\ifGPcolor
    \GPcolorfalse
  }{}%
  \@ifundefined{ifGPblacktext}{%
    \newif\ifGPblacktext
    \GPblacktexttrue
  }{}%
  \let\gplgaddtomacro\g@addto@macro
  \gdef\gplbacktext{}%
  \gdef\gplfronttext{}%
  \makeatother
  \ifGPblacktext
    \def\colorrgb#1{}%
    \def\colorgray#1{}%
  \else
    \ifGPcolor
      \def\colorrgb#1{\color[rgb]{#1}}%
      \def\colorgray#1{\color[gray]{#1}}%
      \expandafter\def\csname LTw\endcsname{\color{white}}%
      \expandafter\def\csname LTb\endcsname{\color{black}}%
      \expandafter\def\csname LTa\endcsname{\color{black}}%
      \expandafter\def\csname LT0\endcsname{\color[rgb]{1,0,0}}%
      \expandafter\def\csname LT1\endcsname{\color[rgb]{0,1,0}}%
      \expandafter\def\csname LT2\endcsname{\color[rgb]{0,0,1}}%
      \expandafter\def\csname LT3\endcsname{\color[rgb]{1,0,1}}%
      \expandafter\def\csname LT4\endcsname{\color[rgb]{0,1,1}}%
      \expandafter\def\csname LT5\endcsname{\color[rgb]{1,1,0}}%
      \expandafter\def\csname LT6\endcsname{\color[rgb]{0,0,0}}%
      \expandafter\def\csname LT7\endcsname{\color[rgb]{1,0.3,0}}%
      \expandafter\def\csname LT8\endcsname{\color[rgb]{0.5,0.5,0.5}}%
    \else
      \def\colorrgb#1{\color{black}}%
      \def\colorgray#1{\color[gray]{#1}}%
      \expandafter\def\csname LTw\endcsname{\color{white}}%
      \expandafter\def\csname LTb\endcsname{\color{black}}%
      \expandafter\def\csname LTa\endcsname{\color{black}}%
      \expandafter\def\csname LT0\endcsname{\color{black}}%
      \expandafter\def\csname LT1\endcsname{\color{black}}%
      \expandafter\def\csname LT2\endcsname{\color{black}}%
      \expandafter\def\csname LT3\endcsname{\color{black}}%
      \expandafter\def\csname LT4\endcsname{\color{black}}%
      \expandafter\def\csname LT5\endcsname{\color{black}}%
      \expandafter\def\csname LT6\endcsname{\color{black}}%
      \expandafter\def\csname LT7\endcsname{\color{black}}%
      \expandafter\def\csname LT8\endcsname{\color{black}}%
    \fi
  \fi
  \setlength{\unitlength}{0.0500bp}%
  \begin{picture}(5760.00,4608.00)%
    \gplgaddtomacro\gplbacktext{%
      \csname LTb\endcsname%
      \put(986,480){\rotatebox{90}{\makebox(0,0){\strut{} 0.9}}}%
      \put(986,1353){\rotatebox{90}{\makebox(0,0){\strut{} 0.925}}}%
      \put(986,2225){\rotatebox{90}{\makebox(0,0){\strut{} 0.95}}}%
      \put(986,3098){\rotatebox{90}{\makebox(0,0){\strut{} 0.975}}}%
      \put(986,3971){\rotatebox{90}{\makebox(0,0){\strut{} 1}}}%
      \put(1128,240){\makebox(0,0){\strut{} 0}}%
      \put(1565,240){\makebox(0,0){\strut{} 0.1}}%
      \put(2002,240){\makebox(0,0){\strut{} 0.2}}%
      \put(2438,240){\makebox(0,0){\strut{} 0.3}}%
      \put(2875,240){\makebox(0,0){\strut{} 0.4}}%
      \put(3312,240){\makebox(0,0){\strut{} 0.5}}%
      \put(3749,240){\makebox(0,0){\strut{} 0.6}}%
      \put(4186,240){\makebox(0,0){\strut{} 0.7}}%
      \put(4622,240){\makebox(0,0){\strut{} 0.8}}%
      \put(5059,240){\makebox(0,0){\strut{} 0.9}}%
      \put(5496,240){\makebox(0,0){\strut{} 1}}%
      \put(528,2400){\rotatebox{90}{\makebox(0,0){\strut{}}}}%
      \put(3312,-72){\makebox(0,0){\strut{}}}%
    }%
    \gplgaddtomacro\gplfronttext{%
    }%
    \gplbacktext
    \put(0,0){\includegraphics{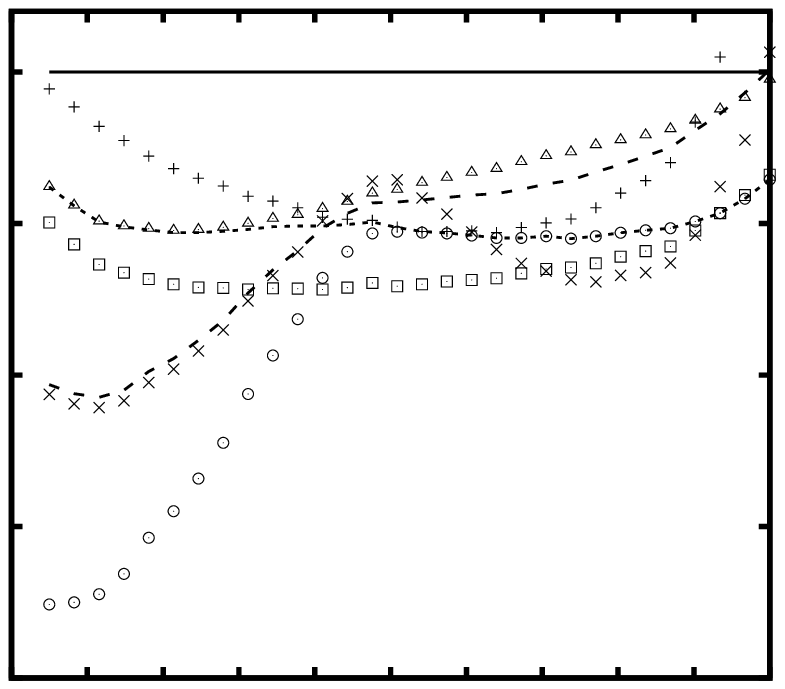}}%
    \gplfronttext
  \end{picture}%
\endgroup

%% file: fig_pi0_3_fdr.tex
\begingroup
  \makeatletter
  \providecommand\color[2][]{%
    \GenericError{(gnuplot) \space\space\space\@spaces}{%
      Package color not loaded in conjunction with
      terminal option `colourtext'%
    }{See the gnuplot documentation for explanation.%
    }{Either use 'blacktext' in gnuplot or load the package
      color.sty in LaTeX.}%
    \renewcommand\color[2][]{}%
  }%
  \providecommand\includegraphics[2][]{%
    \GenericError{(gnuplot) \space\space\space\@spaces}{%
      Package graphicx or graphics not loaded%
    }{See the gnuplot documentation for explanation.%
    }{The gnuplot epslatex terminal needs graphicx.sty or graphics.sty.}%
    \renewcommand\includegraphics[2][]{}%
  }%
  \providecommand\rotatebox[2]{#2}%
  \@ifundefined{ifGPcolor}{%
    \newif\ifGPcolor
    \GPcolorfalse
  }{}%
  \@ifundefined{ifGPblacktext}{%
    \newif\ifGPblacktext
    \GPblacktexttrue
  }{}%
  \let\gplgaddtomacro\g@addto@macro
  \gdef\gplbacktext{}%
  \gdef\gplfronttext{}%
  \makeatother
  \ifGPblacktext
    \def\colorrgb#1{}%
    \def\colorgray#1{}%
  \else
    \ifGPcolor
      \def\colorrgb#1{\color[rgb]{#1}}%
      \def\colorgray#1{\color[gray]{#1}}%
      \expandafter\def\csname LTw\endcsname{\color{white}}%
      \expandafter\def\csname LTb\endcsname{\color{black}}%
      \expandafter\def\csname LTa\endcsname{\color{black}}%
      \expandafter\def\csname LT0\endcsname{\color[rgb]{1,0,0}}%
      \expandafter\def\csname LT1\endcsname{\color[rgb]{0,1,0}}%
      \expandafter\def\csname LT2\endcsname{\color[rgb]{0,0,1}}%
      \expandafter\def\csname LT3\endcsname{\color[rgb]{1,0,1}}%
      \expandafter\def\csname LT4\endcsname{\color[rgb]{0,1,1}}%
      \expandafter\def\csname LT5\endcsname{\color[rgb]{1,1,0}}%
      \expandafter\def\csname LT6\endcsname{\color[rgb]{0,0,0}}%
      \expandafter\def\csname LT7\endcsname{\color[rgb]{1,0.3,0}}%
      \expandafter\def\csname LT8\endcsname{\color[rgb]{0.5,0.5,0.5}}%
    \else
      \def\colorrgb#1{\color{black}}%
      \def\colorgray#1{\color[gray]{#1}}%
      \expandafter\def\csname LTw\endcsname{\color{white}}%
      \expandafter\def\csname LTb\endcsname{\color{black}}%
      \expandafter\def\csname LTa\endcsname{\color{black}}%
      \expandafter\def\csname LT0\endcsname{\color{black}}%
      \expandafter\def\csname LT1\endcsname{\color{black}}%
      \expandafter\def\csname LT2\endcsname{\color{black}}%
      \expandafter\def\csname LT3\endcsname{\color{black}}%
      \expandafter\def\csname LT4\endcsname{\color{black}}%
      \expandafter\def\csname LT5\endcsname{\color{black}}%
      \expandafter\def\csname LT6\endcsname{\color{black}}%
      \expandafter\def\csname LT7\endcsname{\color{black}}%
      \expandafter\def\csname LT8\endcsname{\color{black}}%
    \fi
  \fi
  \setlength{\unitlength}{0.0500bp}%
  \begin{picture}(5760.00,4608.00)%
    \gplgaddtomacro\gplbacktext{%
      \csname LTb\endcsname%
      \put(986,480){\rotatebox{90}{\makebox(0,0){\strut{} 0}}}%
      \put(986,1440){\rotatebox{90}{\makebox(0,0){\strut{} 0.02}}}%
      \put(986,2400){\rotatebox{90}{\makebox(0,0){\strut{} 0.04}}}%
      \put(986,3360){\rotatebox{90}{\makebox(0,0){\strut{} 0.06}}}%
      \put(986,4320){\rotatebox{90}{\makebox(0,0){\strut{} 0.08}}}%
      \put(1128,240){\makebox(0,0){\strut{} 0}}%
      \put(1565,240){\makebox(0,0){\strut{} 0.1}}%
      \put(2002,240){\makebox(0,0){\strut{} 0.2}}%
      \put(2438,240){\makebox(0,0){\strut{} 0.3}}%
      \put(2875,240){\makebox(0,0){\strut{} 0.4}}%
      \put(3312,240){\makebox(0,0){\strut{} 0.5}}%
      \put(3749,240){\makebox(0,0){\strut{} 0.6}}%
      \put(4186,240){\makebox(0,0){\strut{} 0.7}}%
      \put(4622,240){\makebox(0,0){\strut{} 0.8}}%
      \put(5059,240){\makebox(0,0){\strut{} 0.9}}%
      \put(5496,240){\makebox(0,0){\strut{} 1}}%
      \put(528,2400){\rotatebox{90}{\makebox(0,0){\strut{}}}}%
      \put(3312,-72){\makebox(0,0){\strut{}}}%
    }%
    \gplgaddtomacro\gplfronttext{%
    }%
    \gplbacktext
    \put(0,0){\includegraphics{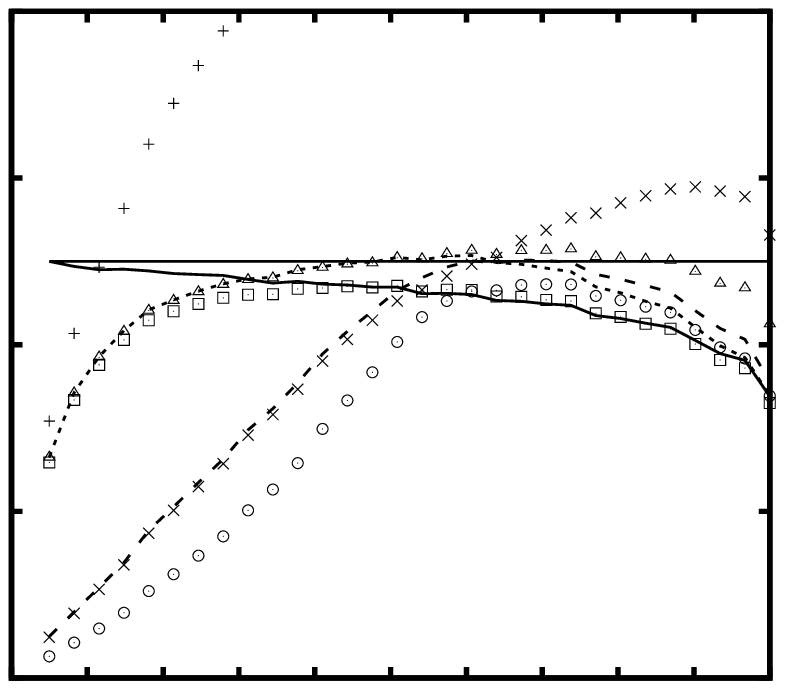}}%
    \gplfronttext
  \end{picture}%
\endgroup

%% file: fig_pi0_3_puis.tex
\begingroup
  \makeatletter
  \providecommand\color[2][]{%
    \GenericError{(gnuplot) \space\space\space\@spaces}{%
      Package color not loaded in conjunction with
      terminal option `colourtext'%
    }{See the gnuplot documentation for explanation.%
    }{Either use 'blacktext' in gnuplot or load the package
      color.sty in LaTeX.}%
    \renewcommand\color[2][]{}%
  }%
  \providecommand\includegraphics[2][]{%
    \GenericError{(gnuplot) \space\space\space\@spaces}{%
      Package graphicx or graphics not loaded%
    }{See the gnuplot documentation for explanation.%
    }{The gnuplot epslatex terminal needs graphicx.sty or graphics.sty.}%
    \renewcommand\includegraphics[2][]{}%
  }%
  \providecommand\rotatebox[2]{#2}%
  \@ifundefined{ifGPcolor}{%
    \newif\ifGPcolor
    \GPcolorfalse
  }{}%
  \@ifundefined{ifGPblacktext}{%
    \newif\ifGPblacktext
    \GPblacktexttrue
  }{}%
  \let\gplgaddtomacro\g@addto@macro
  \gdef\gplbacktext{}%
  \gdef\gplfronttext{}%
  \makeatother
  \ifGPblacktext
    \def\colorrgb#1{}%
    \def\colorgray#1{}%
  \else
    \ifGPcolor
      \def\colorrgb#1{\color[rgb]{#1}}%
      \def\colorgray#1{\color[gray]{#1}}%
      \expandafter\def\csname LTw\endcsname{\color{white}}%
      \expandafter\def\csname LTb\endcsname{\color{black}}%
      \expandafter\def\csname LTa\endcsname{\color{black}}%
      \expandafter\def\csname LT0\endcsname{\color[rgb]{1,0,0}}%
      \expandafter\def\csname LT1\endcsname{\color[rgb]{0,1,0}}%
      \expandafter\def\csname LT2\endcsname{\color[rgb]{0,0,1}}%
      \expandafter\def\csname LT3\endcsname{\color[rgb]{1,0,1}}%
      \expandafter\def\csname LT4\endcsname{\color[rgb]{0,1,1}}%
      \expandafter\def\csname LT5\endcsname{\color[rgb]{1,1,0}}%
      \expandafter\def\csname LT6\endcsname{\color[rgb]{0,0,0}}%
      \expandafter\def\csname LT7\endcsname{\color[rgb]{1,0.3,0}}%
      \expandafter\def\csname LT8\endcsname{\color[rgb]{0.5,0.5,0.5}}%
    \else
      \def\colorrgb#1{\color{black}}%
      \def\colorgray#1{\color[gray]{#1}}%
      \expandafter\def\csname LTw\endcsname{\color{white}}%
      \expandafter\def\csname LTb\endcsname{\color{black}}%
      \expandafter\def\csname LTa\endcsname{\color{black}}%
      \expandafter\def\csname LT0\endcsname{\color{black}}%
      \expandafter\def\csname LT1\endcsname{\color{black}}%
      \expandafter\def\csname LT2\endcsname{\color{black}}%
      \expandafter\def\csname LT3\endcsname{\color{black}}%
      \expandafter\def\csname LT4\endcsname{\color{black}}%
      \expandafter\def\csname LT5\endcsname{\color{black}}%
      \expandafter\def\csname LT6\endcsname{\color{black}}%
      \expandafter\def\csname LT7\endcsname{\color{black}}%
      \expandafter\def\csname LT8\endcsname{\color{black}}%
    \fi
  \fi
  \setlength{\unitlength}{0.0500bp}%
  \begin{picture}(5760.00,4608.00)%
    \gplgaddtomacro\gplbacktext{%
      \csname LTb\endcsname%
      \put(986,480){\rotatebox{90}{\makebox(0,0){\strut{} 0.9}}}%
      \put(986,2225){\rotatebox{90}{\makebox(0,0){\strut{} 0.95}}}%
      \put(986,3971){\rotatebox{90}{\makebox(0,0){\strut{} 1}}}%
      \put(1128,240){\makebox(0,0){\strut{} 0}}%
      \put(1565,240){\makebox(0,0){\strut{} 0.1}}%
      \put(2002,240){\makebox(0,0){\strut{} 0.2}}%
      \put(2438,240){\makebox(0,0){\strut{} 0.3}}%
      \put(2875,240){\makebox(0,0){\strut{} 0.4}}%
      \put(3312,240){\makebox(0,0){\strut{} 0.5}}%
      \put(3749,240){\makebox(0,0){\strut{} 0.6}}%
      \put(4186,240){\makebox(0,0){\strut{} 0.7}}%
      \put(4622,240){\makebox(0,0){\strut{} 0.8}}%
      \put(5059,240){\makebox(0,0){\strut{} 0.9}}%
      \put(5496,240){\makebox(0,0){\strut{} 1}}%
      \put(528,2400){\rotatebox{90}{\makebox(0,0){\strut{}}}}%
      \put(3312,-72){\makebox(0,0){\strut{}}}%
    }%
    \gplgaddtomacro\gplfronttext{%
    }%
    \gplbacktext
    \put(0,0){\includegraphics{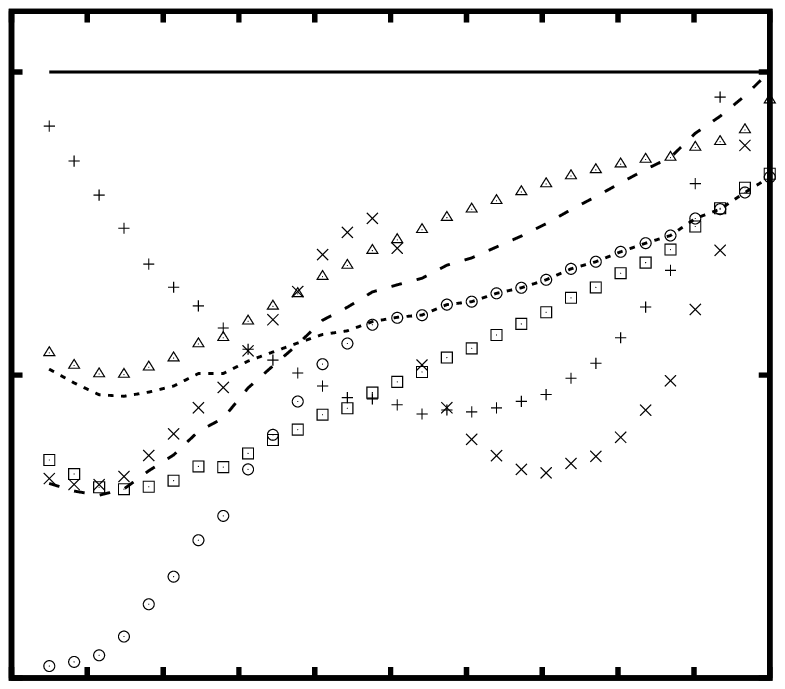}}%
    \gplfronttext
  \end{picture}%
\endgroup

%% file: fig_moy_1_fdr.tex
\begingroup
  \makeatletter
  \providecommand\color[2][]{%
    \GenericError{(gnuplot) \space\space\space\@spaces}{%
      Package color not loaded in conjunction with
      terminal option `colourtext'%
    }{See the gnuplot documentation for explanation.%
    }{Either use 'blacktext' in gnuplot or load the package
      color.sty in LaTeX.}%
    \renewcommand\color[2][]{}%
  }%
  \providecommand\includegraphics[2][]{%
    \GenericError{(gnuplot) \space\space\space\@spaces}{%
      Package graphicx or graphics not loaded%
    }{See the gnuplot documentation for explanation.%
    }{The gnuplot epslatex terminal needs graphicx.sty or graphics.sty.}%
    \renewcommand\includegraphics[2][]{}%
  }%
  \providecommand\rotatebox[2]{#2}%
  \@ifundefined{ifGPcolor}{%
    \newif\ifGPcolor
    \GPcolorfalse
  }{}%
  \@ifundefined{ifGPblacktext}{%
    \newif\ifGPblacktext
    \GPblacktexttrue
  }{}%
  \let\gplgaddtomacro\g@addto@macro
  \gdef\gplbacktext{}%
  \gdef\gplfronttext{}%
  \makeatother
  \ifGPblacktext
    \def\colorrgb#1{}%
    \def\colorgray#1{}%
  \else
    \ifGPcolor
      \def\colorrgb#1{\color[rgb]{#1}}%
      \def\colorgray#1{\color[gray]{#1}}%
      \expandafter\def\csname LTw\endcsname{\color{white}}%
      \expandafter\def\csname LTb\endcsname{\color{black}}%
      \expandafter\def\csname LTa\endcsname{\color{black}}%
      \expandafter\def\csname LT0\endcsname{\color[rgb]{1,0,0}}%
      \expandafter\def\csname LT1\endcsname{\color[rgb]{0,1,0}}%
      \expandafter\def\csname LT2\endcsname{\color[rgb]{0,0,1}}%
      \expandafter\def\csname LT3\endcsname{\color[rgb]{1,0,1}}%
      \expandafter\def\csname LT4\endcsname{\color[rgb]{0,1,1}}%
      \expandafter\def\csname LT5\endcsname{\color[rgb]{1,1,0}}%
      \expandafter\def\csname LT6\endcsname{\color[rgb]{0,0,0}}%
      \expandafter\def\csname LT7\endcsname{\color[rgb]{1,0.3,0}}%
      \expandafter\def\csname LT8\endcsname{\color[rgb]{0.5,0.5,0.5}}%
    \else
      \def\colorrgb#1{\color{black}}%
      \def\colorgray#1{\color[gray]{#1}}%
      \expandafter\def\csname LTw\endcsname{\color{white}}%
      \expandafter\def\csname LTb\endcsname{\color{black}}%
      \expandafter\def\csname LTa\endcsname{\color{black}}%
      \expandafter\def\csname LT0\endcsname{\color{black}}%
      \expandafter\def\csname LT1\endcsname{\color{black}}%
      \expandafter\def\csname LT2\endcsname{\color{black}}%
      \expandafter\def\csname LT3\endcsname{\color{black}}%
      \expandafter\def\csname LT4\endcsname{\color{black}}%
      \expandafter\def\csname LT5\endcsname{\color{black}}%
      \expandafter\def\csname LT6\endcsname{\color{black}}%
      \expandafter\def\csname LT7\endcsname{\color{black}}%
      \expandafter\def\csname LT8\endcsname{\color{black}}%
    \fi
  \fi
  \setlength{\unitlength}{0.0500bp}%
  \begin{picture}(5760.00,4608.00)%
    \gplgaddtomacro\gplbacktext{%
      \csname LTb\endcsname%
      \put(986,1248){\rotatebox{90}{\makebox(0,0){\strut{} 0.02}}}%
      \put(986,2784){\rotatebox{90}{\makebox(0,0){\strut{} 0.04}}}%
      \put(986,4320){\rotatebox{90}{\makebox(0,0){\strut{} 0.06}}}%
      \put(1175,240){\makebox(0,0){\strut{} 0.5}}%
      \put(1650,240){\makebox(0,0){\strut{} 1}}%
      \put(2125,240){\makebox(0,0){\strut{} 1.5}}%
      \put(2600,240){\makebox(0,0){\strut{} 2}}%
      \put(3075,240){\makebox(0,0){\strut{} 2.5}}%
      \put(3549,240){\makebox(0,0){\strut{} 3}}%
      \put(4024,240){\makebox(0,0){\strut{} 3.5}}%
      \put(4499,240){\makebox(0,0){\strut{} 4}}%
      \put(4974,240){\makebox(0,0){\strut{} 4.5}}%
      \put(5449,240){\makebox(0,0){\strut{} 5}}%
      \put(528,2400){\rotatebox{90}{\makebox(0,0){\strut{}}}}%
      \put(3312,-72){\makebox(0,0){\strut{}}}%
    }%
    \gplgaddtomacro\gplfronttext{%
      \csname LTb\endcsname%
      \put(4200,2187){\makebox(0,0)[l]{\strut{}LSU-Oracle}}%
      \csname LTb\endcsname%
      \put(4200,1971){\makebox(0,0)[l]{\strut{}Storey-$\alpha$}}%
      \csname LTb\endcsname%
      \put(4200,1755){\makebox(0,0)[l]{\strut{}Storey-1/2}}%
      \csname LTb\endcsname%
      \put(4200,1539){\makebox(0,0)[l]{\strut{}Median LSU}}%
      \csname LTb\endcsname%
      \put(4200,1323){\makebox(0,0)[l]{\strut{}BKY06}}%
      \csname LTb\endcsname%
      \put(4200,1107){\makebox(0,0)[l]{\strut{}BR-1S-$\alpha$}}%
      \csname LTb\endcsname%
      \put(4200,891){\makebox(0,0)[l]{\strut{}BR-2S-$\alpha$}}%
      \csname LTb\endcsname%
      \put(4200,675){\makebox(0,0)[l]{\strut{}FDR09-1/2}}%
    }%
    \gplbacktext
    \put(0,0){\includegraphics{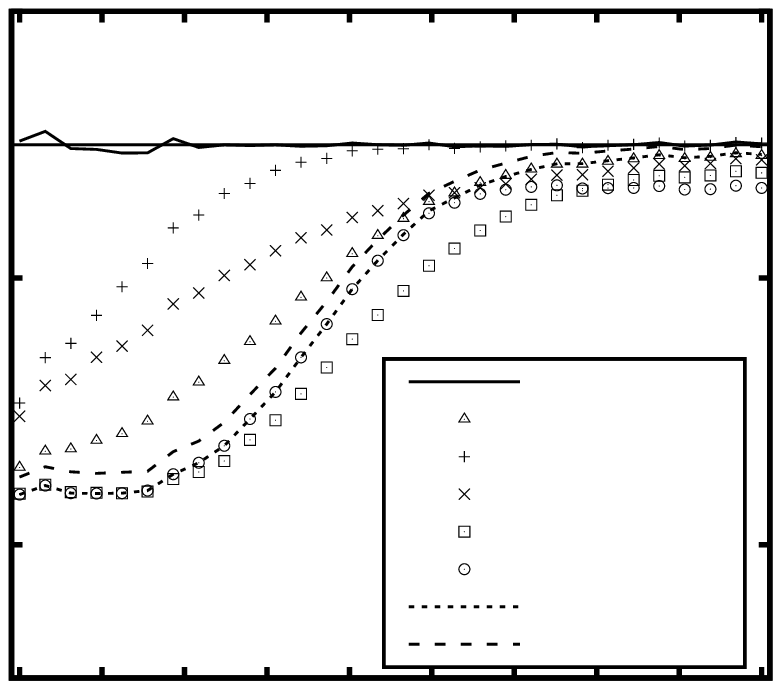}}%
    \gplfronttext
  \end{picture}%
\endgroup

%% file: fig_moy_1_puis.tex
\begingroup
  \makeatletter
  \providecommand\color[2][]{%
    \GenericError{(gnuplot) \space\space\space\@spaces}{%
      Package color not loaded in conjunction with
      terminal option `colourtext'%
    }{See the gnuplot documentation for explanation.%
    }{Either use 'blacktext' in gnuplot or load the package
      color.sty in LaTeX.}%
    \renewcommand\color[2][]{}%
  }%
  \providecommand\includegraphics[2][]{%
    \GenericError{(gnuplot) \space\space\space\@spaces}{%
      Package graphicx or graphics not loaded%
    }{See the gnuplot documentation for explanation.%
    }{The gnuplot epslatex terminal needs graphicx.sty or graphics.sty.}%
    \renewcommand\includegraphics[2][]{}%
  }%
  \providecommand\rotatebox[2]{#2}%
  \@ifundefined{ifGPcolor}{%
    \newif\ifGPcolor
    \GPcolorfalse
  }{}%
  \@ifundefined{ifGPblacktext}{%
    \newif\ifGPblacktext
    \GPblacktexttrue
  }{}%
  \let\gplgaddtomacro\g@addto@macro
  \gdef\gplbacktext{}%
  \gdef\gplfronttext{}%
  \makeatother
  \ifGPblacktext
    \def\colorrgb#1{}%
    \def\colorgray#1{}%
  \else
    \ifGPcolor
      \def\colorrgb#1{\color[rgb]{#1}}%
      \def\colorgray#1{\color[gray]{#1}}%
      \expandafter\def\csname LTw\endcsname{\color{white}}%
      \expandafter\def\csname LTb\endcsname{\color{black}}%
      \expandafter\def\csname LTa\endcsname{\color{black}}%
      \expandafter\def\csname LT0\endcsname{\color[rgb]{1,0,0}}%
      \expandafter\def\csname LT1\endcsname{\color[rgb]{0,1,0}}%
      \expandafter\def\csname LT2\endcsname{\color[rgb]{0,0,1}}%
      \expandafter\def\csname LT3\endcsname{\color[rgb]{1,0,1}}%
      \expandafter\def\csname LT4\endcsname{\color[rgb]{0,1,1}}%
      \expandafter\def\csname LT5\endcsname{\color[rgb]{1,1,0}}%
      \expandafter\def\csname LT6\endcsname{\color[rgb]{0,0,0}}%
      \expandafter\def\csname LT7\endcsname{\color[rgb]{1,0.3,0}}%
      \expandafter\def\csname LT8\endcsname{\color[rgb]{0.5,0.5,0.5}}%
    \else
      \def\colorrgb#1{\color{black}}%
      \def\colorgray#1{\color[gray]{#1}}%
      \expandafter\def\csname LTw\endcsname{\color{white}}%
      \expandafter\def\csname LTb\endcsname{\color{black}}%
      \expandafter\def\csname LTa\endcsname{\color{black}}%
      \expandafter\def\csname LT0\endcsname{\color{black}}%
      \expandafter\def\csname LT1\endcsname{\color{black}}%
      \expandafter\def\csname LT2\endcsname{\color{black}}%
      \expandafter\def\csname LT3\endcsname{\color{black}}%
      \expandafter\def\csname LT4\endcsname{\color{black}}%
      \expandafter\def\csname LT5\endcsname{\color{black}}%
      \expandafter\def\csname LT6\endcsname{\color{black}}%
      \expandafter\def\csname LT7\endcsname{\color{black}}%
      \expandafter\def\csname LT8\endcsname{\color{black}}%
    \fi
  \fi
  \setlength{\unitlength}{0.0500bp}%
  \begin{picture}(5760.00,4608.00)%
    \gplgaddtomacro\gplbacktext{%
      \csname LTb\endcsname%
      \put(986,480){\rotatebox{90}{\makebox(0,0){\strut{} 0.4}}}%
      \put(986,1739){\rotatebox{90}{\makebox(0,0){\strut{} 0.6}}}%
      \put(986,2998){\rotatebox{90}{\makebox(0,0){\strut{} 0.8}}}%
      \put(986,4257){\rotatebox{90}{\makebox(0,0){\strut{} 1}}}%
      \put(1175,240){\makebox(0,0){\strut{} 0.5}}%
      \put(1650,240){\makebox(0,0){\strut{} 1}}%
      \put(2125,240){\makebox(0,0){\strut{} 1.5}}%
      \put(2600,240){\makebox(0,0){\strut{} 2}}%
      \put(3075,240){\makebox(0,0){\strut{} 2.5}}%
      \put(3549,240){\makebox(0,0){\strut{} 3}}%
      \put(4024,240){\makebox(0,0){\strut{} 3.5}}%
      \put(4499,240){\makebox(0,0){\strut{} 4}}%
      \put(4974,240){\makebox(0,0){\strut{} 4.5}}%
      \put(5449,240){\makebox(0,0){\strut{} 5}}%
      \put(528,2400){\rotatebox{90}{\makebox(0,0){\strut{}}}}%
      \put(3312,-72){\makebox(0,0){\strut{}}}%
    }%
    \gplgaddtomacro\gplfronttext{%
      \csname LTb\endcsname%
      \put(4200,1971){\makebox(0,0)[l]{\strut{}Storey-$\alpha$}}%
      \csname LTb\endcsname%
      \put(4200,1755){\makebox(0,0)[l]{\strut{}Storey-1/2}}%
      \csname LTb\endcsname%
      \put(4200,1539){\makebox(0,0)[l]{\strut{}Median LSU}}%
      \csname LTb\endcsname%
      \put(4200,1323){\makebox(0,0)[l]{\strut{}BKY06}}%
      \csname LTb\endcsname%
      \put(4200,1107){\makebox(0,0)[l]{\strut{}BR-1S-$\alpha$}}%
      \csname LTb\endcsname%
      \put(4200,891){\makebox(0,0)[l]{\strut{}BR-2S-$\alpha$}}%
      \csname LTb\endcsname%
      \put(4200,675){\makebox(0,0)[l]{\strut{}FDR09-1/2}}%
    }%
    \gplbacktext
    \put(0,0){\includegraphics{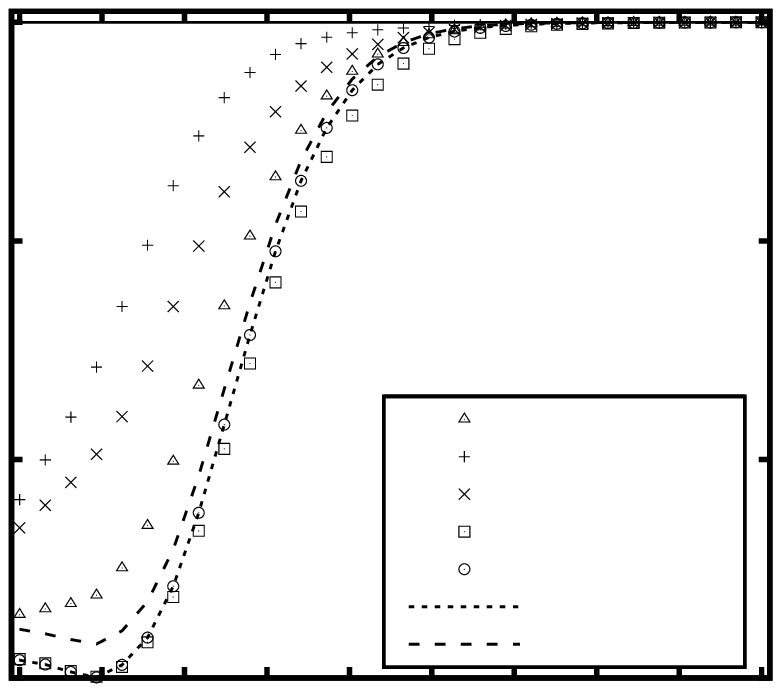}}%
    \gplfronttext
  \end{picture}%
\endgroup

%% file: fig_moy_2_fdr.tex
\begingroup
  \makeatletter
  \providecommand\color[2][]{%
    \GenericError{(gnuplot) \space\space\space\@spaces}{%
      Package color not loaded in conjunction with
      terminal option `colourtext'%
    }{See the gnuplot documentation for explanation.%
    }{Either use 'blacktext' in gnuplot or load the package
      color.sty in LaTeX.}%
    \renewcommand\color[2][]{}%
  }%
  \providecommand\includegraphics[2][]{%
    \GenericError{(gnuplot) \space\space\space\@spaces}{%
      Package graphicx or graphics not loaded%
    }{See the gnuplot documentation for explanation.%
    }{The gnuplot epslatex terminal needs graphicx.sty or graphics.sty.}%
    \renewcommand\includegraphics[2][]{}%
  }%
  \providecommand\rotatebox[2]{#2}%
  \@ifundefined{ifGPcolor}{%
    \newif\ifGPcolor
    \GPcolorfalse
  }{}%
  \@ifundefined{ifGPblacktext}{%
    \newif\ifGPblacktext
    \GPblacktexttrue
  }{}%
  \let\gplgaddtomacro\g@addto@macro
  \gdef\gplbacktext{}%
  \gdef\gplfronttext{}%
  \makeatother
  \ifGPblacktext
    \def\colorrgb#1{}%
    \def\colorgray#1{}%
  \else
    \ifGPcolor
      \def\colorrgb#1{\color[rgb]{#1}}%
      \def\colorgray#1{\color[gray]{#1}}%
      \expandafter\def\csname LTw\endcsname{\color{white}}%
      \expandafter\def\csname LTb\endcsname{\color{black}}%
      \expandafter\def\csname LTa\endcsname{\color{black}}%
      \expandafter\def\csname LT0\endcsname{\color[rgb]{1,0,0}}%
      \expandafter\def\csname LT1\endcsname{\color[rgb]{0,1,0}}%
      \expandafter\def\csname LT2\endcsname{\color[rgb]{0,0,1}}%
      \expandafter\def\csname LT3\endcsname{\color[rgb]{1,0,1}}%
      \expandafter\def\csname LT4\endcsname{\color[rgb]{0,1,1}}%
      \expandafter\def\csname LT5\endcsname{\color[rgb]{1,1,0}}%
      \expandafter\def\csname LT6\endcsname{\color[rgb]{0,0,0}}%
      \expandafter\def\csname LT7\endcsname{\color[rgb]{1,0.3,0}}%
      \expandafter\def\csname LT8\endcsname{\color[rgb]{0.5,0.5,0.5}}%
    \else
      \def\colorrgb#1{\color{black}}%
      \def\colorgray#1{\color[gray]{#1}}%
      \expandafter\def\csname LTw\endcsname{\color{white}}%
      \expandafter\def\csname LTb\endcsname{\color{black}}%
      \expandafter\def\csname LTa\endcsname{\color{black}}%
      \expandafter\def\csname LT0\endcsname{\color{black}}%
      \expandafter\def\csname LT1\endcsname{\color{black}}%
      \expandafter\def\csname LT2\endcsname{\color{black}}%
      \expandafter\def\csname LT3\endcsname{\color{black}}%
      \expandafter\def\csname LT4\endcsname{\color{black}}%
      \expandafter\def\csname LT5\endcsname{\color{black}}%
      \expandafter\def\csname LT6\endcsname{\color{black}}%
      \expandafter\def\csname LT7\endcsname{\color{black}}%
      \expandafter\def\csname LT8\endcsname{\color{black}}%
    \fi
  \fi
  \setlength{\unitlength}{0.0500bp}%
  \begin{picture}(5760.00,4608.00)%
    \gplgaddtomacro\gplbacktext{%
      \csname LTb\endcsname%
      \put(986,1248){\rotatebox{90}{\makebox(0,0){\strut{} 0.02}}}%
      \put(986,2784){\rotatebox{90}{\makebox(0,0){\strut{} 0.04}}}%
      \put(986,4320){\rotatebox{90}{\makebox(0,0){\strut{} 0.06}}}%
      \put(1175,240){\makebox(0,0){\strut{} 0.5}}%
      \put(1650,240){\makebox(0,0){\strut{} 1}}%
      \put(2125,240){\makebox(0,0){\strut{} 1.5}}%
      \put(2600,240){\makebox(0,0){\strut{} 2}}%
      \put(3075,240){\makebox(0,0){\strut{} 2.5}}%
      \put(3549,240){\makebox(0,0){\strut{} 3}}%
      \put(4024,240){\makebox(0,0){\strut{} 3.5}}%
      \put(4499,240){\makebox(0,0){\strut{} 4}}%
      \put(4974,240){\makebox(0,0){\strut{} 4.5}}%
      \put(5449,240){\makebox(0,0){\strut{} 5}}%
      \put(528,2400){\rotatebox{90}{\makebox(0,0){\strut{}}}}%
      \put(3312,-72){\makebox(0,0){\strut{}}}%
    }%
    \gplgaddtomacro\gplfronttext{%
    }%
    \gplbacktext
    \put(0,0){\includegraphics{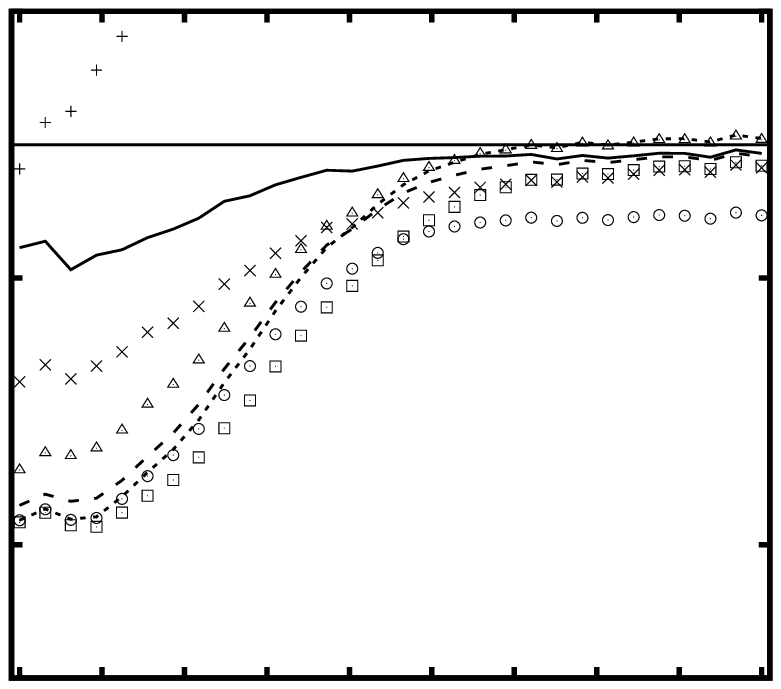}}%
    \gplfronttext
  \end{picture}%
\endgroup

%% file: fig_moy_2_puis.tex
\begingroup
  \makeatletter
  \providecommand\color[2][]{%
    \GenericError{(gnuplot) \space\space\space\@spaces}{%
      Package color not loaded in conjunction with
      terminal option `colourtext'%
    }{See the gnuplot documentation for explanation.%
    }{Either use 'blacktext' in gnuplot or load the package
      color.sty in LaTeX.}%
    \renewcommand\color[2][]{}%
  }%
  \providecommand\includegraphics[2][]{%
    \GenericError{(gnuplot) \space\space\space\@spaces}{%
      Package graphicx or graphics not loaded%
    }{See the gnuplot documentation for explanation.%
    }{The gnuplot epslatex terminal needs graphicx.sty or graphics.sty.}%
    \renewcommand\includegraphics[2][]{}%
  }%
  \providecommand\rotatebox[2]{#2}%
  \@ifundefined{ifGPcolor}{%
    \newif\ifGPcolor
    \GPcolorfalse
  }{}%
  \@ifundefined{ifGPblacktext}{%
    \newif\ifGPblacktext
    \GPblacktexttrue
  }{}%
  \let\gplgaddtomacro\g@addto@macro
  \gdef\gplbacktext{}%
  \gdef\gplfronttext{}%
  \makeatother
  \ifGPblacktext
    \def\colorrgb#1{}%
    \def\colorgray#1{}%
  \else
    \ifGPcolor
      \def\colorrgb#1{\color[rgb]{#1}}%
      \def\colorgray#1{\color[gray]{#1}}%
      \expandafter\def\csname LTw\endcsname{\color{white}}%
      \expandafter\def\csname LTb\endcsname{\color{black}}%
      \expandafter\def\csname LTa\endcsname{\color{black}}%
      \expandafter\def\csname LT0\endcsname{\color[rgb]{1,0,0}}%
      \expandafter\def\csname LT1\endcsname{\color[rgb]{0,1,0}}%
      \expandafter\def\csname LT2\endcsname{\color[rgb]{0,0,1}}%
      \expandafter\def\csname LT3\endcsname{\color[rgb]{1,0,1}}%
      \expandafter\def\csname LT4\endcsname{\color[rgb]{0,1,1}}%
      \expandafter\def\csname LT5\endcsname{\color[rgb]{1,1,0}}%
      \expandafter\def\csname LT6\endcsname{\color[rgb]{0,0,0}}%
      \expandafter\def\csname LT7\endcsname{\color[rgb]{1,0.3,0}}%
      \expandafter\def\csname LT8\endcsname{\color[rgb]{0.5,0.5,0.5}}%
    \else
      \def\colorrgb#1{\color{black}}%
      \def\colorgray#1{\color[gray]{#1}}%
      \expandafter\def\csname LTw\endcsname{\color{white}}%
      \expandafter\def\csname LTb\endcsname{\color{black}}%
      \expandafter\def\csname LTa\endcsname{\color{black}}%
      \expandafter\def\csname LT0\endcsname{\color{black}}%
      \expandafter\def\csname LT1\endcsname{\color{black}}%
      \expandafter\def\csname LT2\endcsname{\color{black}}%
      \expandafter\def\csname LT3\endcsname{\color{black}}%
      \expandafter\def\csname LT4\endcsname{\color{black}}%
      \expandafter\def\csname LT5\endcsname{\color{black}}%
      \expandafter\def\csname LT6\endcsname{\color{black}}%
      \expandafter\def\csname LT7\endcsname{\color{black}}%
      \expandafter\def\csname LT8\endcsname{\color{black}}%
    \fi
  \fi
  \setlength{\unitlength}{0.0500bp}%
  \begin{picture}(5760.00,4608.00)%
    \gplgaddtomacro\gplbacktext{%
      \csname LTb\endcsname%
      \put(986,907){\rotatebox{90}{\makebox(0,0){\strut{} 0.4}}}%
      \put(986,1760){\rotatebox{90}{\makebox(0,0){\strut{} 0.6}}}%
      \put(986,2613){\rotatebox{90}{\makebox(0,0){\strut{} 0.8}}}%
      \put(986,3467){\rotatebox{90}{\makebox(0,0){\strut{} 1}}}%
      \put(986,4320){\rotatebox{90}{\makebox(0,0){\strut{} 1.2}}}%
      \put(1175,240){\makebox(0,0){\strut{} 0.5}}%
      \put(1650,240){\makebox(0,0){\strut{} 1}}%
      \put(2125,240){\makebox(0,0){\strut{} 1.5}}%
      \put(2600,240){\makebox(0,0){\strut{} 2}}%
      \put(3075,240){\makebox(0,0){\strut{} 2.5}}%
      \put(3549,240){\makebox(0,0){\strut{} 3}}%
      \put(4024,240){\makebox(0,0){\strut{} 3.5}}%
      \put(4499,240){\makebox(0,0){\strut{} 4}}%
      \put(4974,240){\makebox(0,0){\strut{} 4.5}}%
      \put(5449,240){\makebox(0,0){\strut{} 5}}%
      \put(528,2400){\rotatebox{90}{\makebox(0,0){\strut{}}}}%
      \put(3312,-72){\makebox(0,0){\strut{}}}%
    }%
    \gplgaddtomacro\gplfronttext{%
    }%
    \gplbacktext
    \put(0,0){\includegraphics{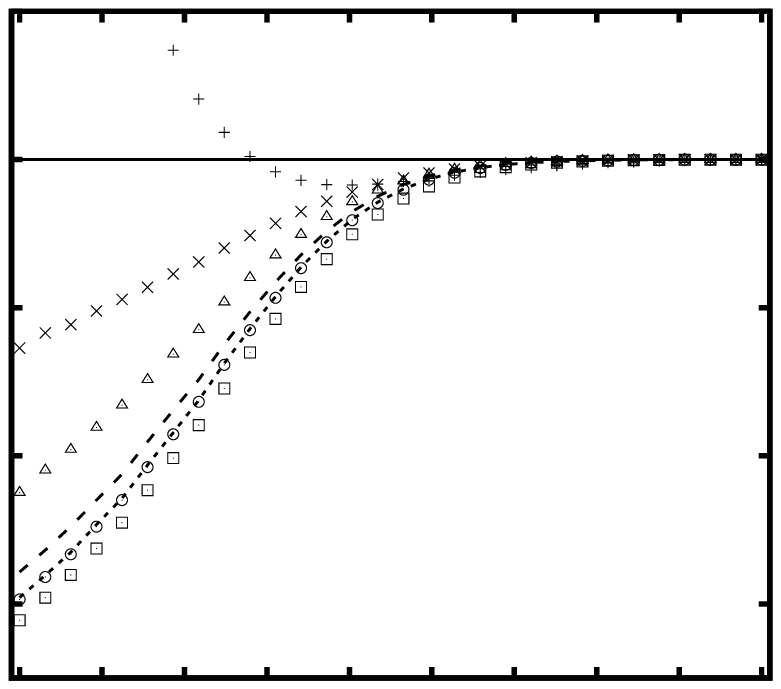}}%
    \gplfronttext
  \end{picture}%
\endgroup

%% file: fig_moy_3_fdr.tex
\begingroup
  \makeatletter
  \providecommand\color[2][]{%
    \GenericError{(gnuplot) \space\space\space\@spaces}{%
      Package color not loaded in conjunction with
      terminal option `colourtext'%
    }{See the gnuplot documentation for explanation.%
    }{Either use 'blacktext' in gnuplot or load the package
      color.sty in LaTeX.}%
    \renewcommand\color[2][]{}%
  }%
  \providecommand\includegraphics[2][]{%
    \GenericError{(gnuplot) \space\space\space\@spaces}{%
      Package graphicx or graphics not loaded%
    }{See the gnuplot documentation for explanation.%
    }{The gnuplot epslatex terminal needs graphicx.sty or graphics.sty.}%
    \renewcommand\includegraphics[2][]{}%
  }%
  \providecommand\rotatebox[2]{#2}%
  \@ifundefined{ifGPcolor}{%
    \newif\ifGPcolor
    \GPcolorfalse
  }{}%
  \@ifundefined{ifGPblacktext}{%
    \newif\ifGPblacktext
    \GPblacktexttrue
  }{}%
  \let\gplgaddtomacro\g@addto@macro
  \gdef\gplbacktext{}%
  \gdef\gplfronttext{}%
  \makeatother
  \ifGPblacktext
    \def\colorrgb#1{}%
    \def\colorgray#1{}%
  \else
    \ifGPcolor
      \def\colorrgb#1{\color[rgb]{#1}}%
      \def\colorgray#1{\color[gray]{#1}}%
      \expandafter\def\csname LTw\endcsname{\color{white}}%
      \expandafter\def\csname LTb\endcsname{\color{black}}%
      \expandafter\def\csname LTa\endcsname{\color{black}}%
      \expandafter\def\csname LT0\endcsname{\color[rgb]{1,0,0}}%
      \expandafter\def\csname LT1\endcsname{\color[rgb]{0,1,0}}%
      \expandafter\def\csname LT2\endcsname{\color[rgb]{0,0,1}}%
      \expandafter\def\csname LT3\endcsname{\color[rgb]{1,0,1}}%
      \expandafter\def\csname LT4\endcsname{\color[rgb]{0,1,1}}%
      \expandafter\def\csname LT5\endcsname{\color[rgb]{1,1,0}}%
      \expandafter\def\csname LT6\endcsname{\color[rgb]{0,0,0}}%
      \expandafter\def\csname LT7\endcsname{\color[rgb]{1,0.3,0}}%
      \expandafter\def\csname LT8\endcsname{\color[rgb]{0.5,0.5,0.5}}%
    \else
      \def\colorrgb#1{\color{black}}%
      \def\colorgray#1{\color[gray]{#1}}%
      \expandafter\def\csname LTw\endcsname{\color{white}}%
      \expandafter\def\csname LTb\endcsname{\color{black}}%
      \expandafter\def\csname LTa\endcsname{\color{black}}%
      \expandafter\def\csname LT0\endcsname{\color{black}}%
      \expandafter\def\csname LT1\endcsname{\color{black}}%
      \expandafter\def\csname LT2\endcsname{\color{black}}%
      \expandafter\def\csname LT3\endcsname{\color{black}}%
      \expandafter\def\csname LT4\endcsname{\color{black}}%
      \expandafter\def\csname LT5\endcsname{\color{black}}%
      \expandafter\def\csname LT6\endcsname{\color{black}}%
      \expandafter\def\csname LT7\endcsname{\color{black}}%
      \expandafter\def\csname LT8\endcsname{\color{black}}%
    \fi
  \fi
  \setlength{\unitlength}{0.0500bp}%
  \begin{picture}(5760.00,4608.00)%
    \gplgaddtomacro\gplbacktext{%
      \csname LTb\endcsname%
      \put(986,992){\rotatebox{90}{\makebox(0,0){\strut{} 0.02}}}%
      \put(986,2016){\rotatebox{90}{\makebox(0,0){\strut{} 0.04}}}%
      \put(986,3040){\rotatebox{90}{\makebox(0,0){\strut{} 0.06}}}%
      \put(986,4064){\rotatebox{90}{\makebox(0,0){\strut{} 0.08}}}%
      \put(1175,240){\makebox(0,0){\strut{} 0.5}}%
      \put(1650,240){\makebox(0,0){\strut{} 1}}%
      \put(2125,240){\makebox(0,0){\strut{} 1.5}}%
      \put(2600,240){\makebox(0,0){\strut{} 2}}%
      \put(3075,240){\makebox(0,0){\strut{} 2.5}}%
      \put(3549,240){\makebox(0,0){\strut{} 3}}%
      \put(4024,240){\makebox(0,0){\strut{} 3.5}}%
      \put(4499,240){\makebox(0,0){\strut{} 4}}%
      \put(4974,240){\makebox(0,0){\strut{} 4.5}}%
      \put(5449,240){\makebox(0,0){\strut{} 5}}%
      \put(528,2400){\rotatebox{90}{\makebox(0,0){\strut{}}}}%
      \put(3312,-72){\makebox(0,0){\strut{}}}%
    }%
    \gplgaddtomacro\gplfronttext{%
    }%
    \gplbacktext
    \put(0,0){\includegraphics{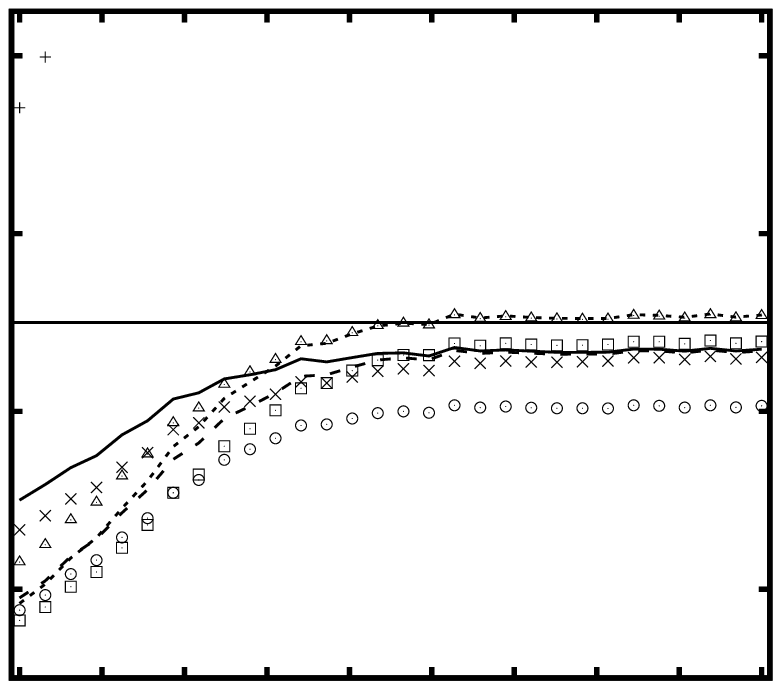}}%
    \gplfronttext
  \end{picture}%
\endgroup

%% file: fig_moy_3_puis.tex
\begingroup
  \makeatletter
  \providecommand\color[2][]{%
    \GenericError{(gnuplot) \space\space\space\@spaces}{%
      Package color not loaded in conjunction with
      terminal option `colourtext'%
    }{See the gnuplot documentation for explanation.%
    }{Either use 'blacktext' in gnuplot or load the package
      color.sty in LaTeX.}%
    \renewcommand\color[2][]{}%
  }%
  \providecommand\includegraphics[2][]{%
    \GenericError{(gnuplot) \space\space\space\@spaces}{%
      Package graphicx or graphics not loaded%
    }{See the gnuplot documentation for explanation.%
    }{The gnuplot epslatex terminal needs graphicx.sty or graphics.sty.}%
    \renewcommand\includegraphics[2][]{}%
  }%
  \providecommand\rotatebox[2]{#2}%
  \@ifundefined{ifGPcolor}{%
    \newif\ifGPcolor
    \GPcolorfalse
  }{}%
  \@ifundefined{ifGPblacktext}{%
    \newif\ifGPblacktext
    \GPblacktexttrue
  }{}%
  \let\gplgaddtomacro\g@addto@macro
  \gdef\gplbacktext{}%
  \gdef\gplfronttext{}%
  \makeatother
  \ifGPblacktext
    \def\colorrgb#1{}%
    \def\colorgray#1{}%
  \else
    \ifGPcolor
      \def\colorrgb#1{\color[rgb]{#1}}%
      \def\colorgray#1{\color[gray]{#1}}%
      \expandafter\def\csname LTw\endcsname{\color{white}}%
      \expandafter\def\csname LTb\endcsname{\color{black}}%
      \expandafter\def\csname LTa\endcsname{\color{black}}%
      \expandafter\def\csname LT0\endcsname{\color[rgb]{1,0,0}}%
      \expandafter\def\csname LT1\endcsname{\color[rgb]{0,1,0}}%
      \expandafter\def\csname LT2\endcsname{\color[rgb]{0,0,1}}%
      \expandafter\def\csname LT3\endcsname{\color[rgb]{1,0,1}}%
      \expandafter\def\csname LT4\endcsname{\color[rgb]{0,1,1}}%
      \expandafter\def\csname LT5\endcsname{\color[rgb]{1,1,0}}%
      \expandafter\def\csname LT6\endcsname{\color[rgb]{0,0,0}}%
      \expandafter\def\csname LT7\endcsname{\color[rgb]{1,0.3,0}}%
      \expandafter\def\csname LT8\endcsname{\color[rgb]{0.5,0.5,0.5}}%
    \else
      \def\colorrgb#1{\color{black}}%
      \def\colorgray#1{\color[gray]{#1}}%
      \expandafter\def\csname LTw\endcsname{\color{white}}%
      \expandafter\def\csname LTb\endcsname{\color{black}}%
      \expandafter\def\csname LTa\endcsname{\color{black}}%
      \expandafter\def\csname LT0\endcsname{\color{black}}%
      \expandafter\def\csname LT1\endcsname{\color{black}}%
      \expandafter\def\csname LT2\endcsname{\color{black}}%
      \expandafter\def\csname LT3\endcsname{\color{black}}%
      \expandafter\def\csname LT4\endcsname{\color{black}}%
      \expandafter\def\csname LT5\endcsname{\color{black}}%
      \expandafter\def\csname LT6\endcsname{\color{black}}%
      \expandafter\def\csname LT7\endcsname{\color{black}}%
      \expandafter\def\csname LT8\endcsname{\color{black}}%
    \fi
  \fi
  \setlength{\unitlength}{0.0500bp}%
  \begin{picture}(5760.00,4608.00)%
    \gplgaddtomacro\gplbacktext{%
      \csname LTb\endcsname%
      \put(986,1120){\rotatebox{90}{\makebox(0,0){\strut{} 0.6}}}%
      \put(986,2400){\rotatebox{90}{\makebox(0,0){\strut{} 0.8}}}%
      \put(986,3680){\rotatebox{90}{\makebox(0,0){\strut{} 1}}}%
      \put(1175,240){\makebox(0,0){\strut{} 0.5}}%
      \put(1650,240){\makebox(0,0){\strut{} 1}}%
      \put(2125,240){\makebox(0,0){\strut{} 1.5}}%
      \put(2600,240){\makebox(0,0){\strut{} 2}}%
      \put(3075,240){\makebox(0,0){\strut{} 2.5}}%
      \put(3549,240){\makebox(0,0){\strut{} 3}}%
      \put(4024,240){\makebox(0,0){\strut{} 3.5}}%
      \put(4499,240){\makebox(0,0){\strut{} 4}}%
      \put(4974,240){\makebox(0,0){\strut{} 4.5}}%
      \put(5449,240){\makebox(0,0){\strut{} 5}}%
      \put(528,2400){\rotatebox{90}{\makebox(0,0){\strut{}}}}%
      \put(3312,-72){\makebox(0,0){\strut{}}}%
    }%
    \gplgaddtomacro\gplfronttext{%
    }%
    \gplbacktext
    \put(0,0){\includegraphics{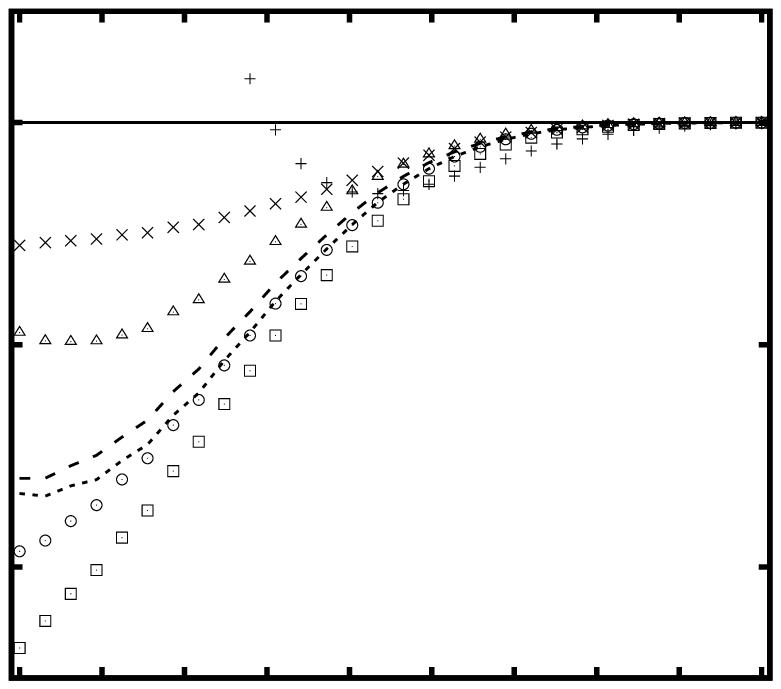}}%
    \gplfronttext
  \end{picture}%
\endgroup

%% file: fig_rho_1_fdr.tex
\begingroup
  \makeatletter
  \providecommand\color[2][]{%
    \GenericError{(gnuplot) \space\space\space\@spaces}{%
      Package color not loaded in conjunction with
      terminal option `colourtext'%
    }{See the gnuplot documentation for explanation.%
    }{Either use 'blacktext' in gnuplot or load the package
      color.sty in LaTeX.}%
    \renewcommand\color[2][]{}%
  }%
  \providecommand\includegraphics[2][]{%
    \GenericError{(gnuplot) \space\space\space\@spaces}{%
      Package graphicx or graphics not loaded%
    }{See the gnuplot documentation for explanation.%
    }{The gnuplot epslatex terminal needs graphicx.sty or graphics.sty.}%
    \renewcommand\includegraphics[2][]{}%
  }%
  \providecommand\rotatebox[2]{#2}%
  \@ifundefined{ifGPcolor}{%
    \newif\ifGPcolor
    \GPcolorfalse
  }{}%
  \@ifundefined{ifGPblacktext}{%
    \newif\ifGPblacktext
    \GPblacktexttrue
  }{}%
  \let\gplgaddtomacro\g@addto@macro
  \gdef\gplbacktext{}%
  \gdef\gplfronttext{}%
  \makeatother
  \ifGPblacktext
    \def\colorrgb#1{}%
    \def\colorgray#1{}%
  \else
    \ifGPcolor
      \def\colorrgb#1{\color[rgb]{#1}}%
      \def\colorgray#1{\color[gray]{#1}}%
      \expandafter\def\csname LTw\endcsname{\color{white}}%
      \expandafter\def\csname LTb\endcsname{\color{black}}%
      \expandafter\def\csname LTa\endcsname{\color{black}}%
      \expandafter\def\csname LT0\endcsname{\color[rgb]{1,0,0}}%
      \expandafter\def\csname LT1\endcsname{\color[rgb]{0,1,0}}%
      \expandafter\def\csname LT2\endcsname{\color[rgb]{0,0,1}}%
      \expandafter\def\csname LT3\endcsname{\color[rgb]{1,0,1}}%
      \expandafter\def\csname LT4\endcsname{\color[rgb]{0,1,1}}%
      \expandafter\def\csname LT5\endcsname{\color[rgb]{1,1,0}}%
      \expandafter\def\csname LT6\endcsname{\color[rgb]{0,0,0}}%
      \expandafter\def\csname LT7\endcsname{\color[rgb]{1,0.3,0}}%
      \expandafter\def\csname LT8\endcsname{\color[rgb]{0.5,0.5,0.5}}%
    \else
      \def\colorrgb#1{\color{black}}%
      \def\colorgray#1{\color[gray]{#1}}%
      \expandafter\def\csname LTw\endcsname{\color{white}}%
      \expandafter\def\csname LTb\endcsname{\color{black}}%
      \expandafter\def\csname LTa\endcsname{\color{black}}%
      \expandafter\def\csname LT0\endcsname{\color{black}}%
      \expandafter\def\csname LT1\endcsname{\color{black}}%
      \expandafter\def\csname LT2\endcsname{\color{black}}%
      \expandafter\def\csname LT3\endcsname{\color{black}}%
      \expandafter\def\csname LT4\endcsname{\color{black}}%
      \expandafter\def\csname LT5\endcsname{\color{black}}%
      \expandafter\def\csname LT6\endcsname{\color{black}}%
      \expandafter\def\csname LT7\endcsname{\color{black}}%
      \expandafter\def\csname LT8\endcsname{\color{black}}%
    \fi
  \fi
  \setlength{\unitlength}{0.0500bp}%
  \begin{picture}(5760.00,4608.00)%
    \gplgaddtomacro\gplbacktext{%
      \csname LTb\endcsname%
      \put(986,480){\rotatebox{90}{\makebox(0,0){\strut{} 0}}}%
      \put(986,1248){\rotatebox{90}{\makebox(0,0){\strut{} 0.02}}}%
      \put(986,2016){\rotatebox{90}{\makebox(0,0){\strut{} 0.04}}}%
      \put(986,2784){\rotatebox{90}{\makebox(0,0){\strut{} 0.06}}}%
      \put(986,3552){\rotatebox{90}{\makebox(0,0){\strut{} 0.08}}}%
      \put(986,4320){\rotatebox{90}{\makebox(0,0){\strut{} 0.1}}}%
      \put(1128,240){\makebox(0,0){\strut{} 0}}%
      \put(2002,240){\makebox(0,0){\strut{} 0.2}}%
      \put(2875,240){\makebox(0,0){\strut{} 0.4}}%
      \put(3749,240){\makebox(0,0){\strut{} 0.6}}%
      \put(4622,240){\makebox(0,0){\strut{} 0.8}}%
      \put(5496,240){\makebox(0,0){\strut{} 1}}%
      \put(528,2400){\rotatebox{90}{\makebox(0,0){\strut{}}}}%
      \put(3312,-72){\makebox(0,0){\strut{}}}%
    }%
    \gplgaddtomacro\gplfronttext{%
    }%
    \gplbacktext
    \put(0,0){\includegraphics{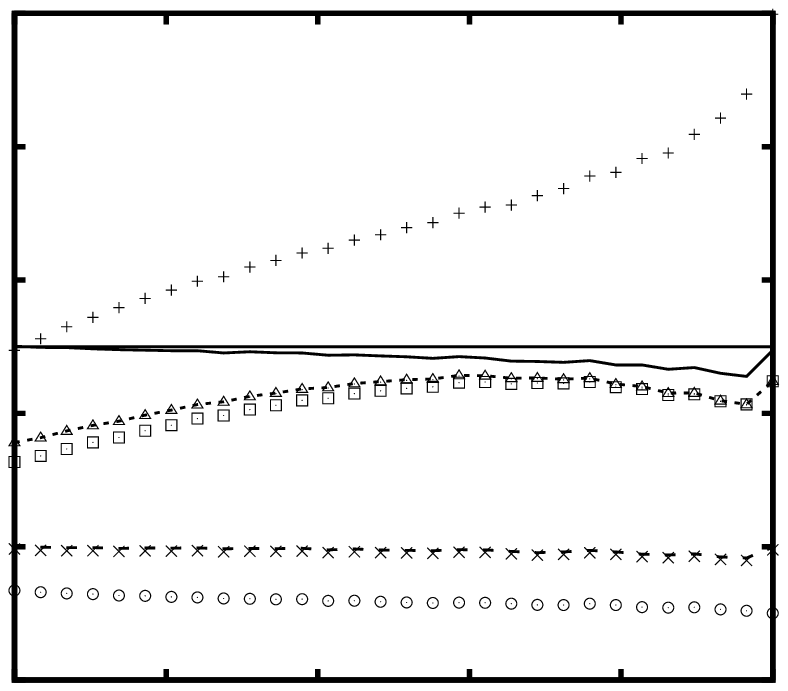}}%
    \gplfronttext
  \end{picture}%
\endgroup

%% file: fig_rho_1_puis.tex
\begingroup
  \makeatletter
  \providecommand\color[2][]{%
    \GenericError{(gnuplot) \space\space\space\@spaces}{%
      Package color not loaded in conjunction with
      terminal option `colourtext'%
    }{See the gnuplot documentation for explanation.%
    }{Either use 'blacktext' in gnuplot or load the package
      color.sty in LaTeX.}%
    \renewcommand\color[2][]{}%
  }%
  \providecommand\includegraphics[2][]{%
    \GenericError{(gnuplot) \space\space\space\@spaces}{%
      Package graphicx or graphics not loaded%
    }{See the gnuplot documentation for explanation.%
    }{The gnuplot epslatex terminal needs graphicx.sty or graphics.sty.}%
    \renewcommand\includegraphics[2][]{}%
  }%
  \providecommand\rotatebox[2]{#2}%
  \@ifundefined{ifGPcolor}{%
    \newif\ifGPcolor
    \GPcolorfalse
  }{}%
  \@ifundefined{ifGPblacktext}{%
    \newif\ifGPblacktext
    \GPblacktexttrue
  }{}%
  \let\gplgaddtomacro\g@addto@macro
  \gdef\gplbacktext{}%
  \gdef\gplfronttext{}%
  \makeatother
  \ifGPblacktext
    \def\colorrgb#1{}%
    \def\colorgray#1{}%
  \else
    \ifGPcolor
      \def\colorrgb#1{\color[rgb]{#1}}%
      \def\colorgray#1{\color[gray]{#1}}%
      \expandafter\def\csname LTw\endcsname{\color{white}}%
      \expandafter\def\csname LTb\endcsname{\color{black}}%
      \expandafter\def\csname LTa\endcsname{\color{black}}%
      \expandafter\def\csname LT0\endcsname{\color[rgb]{1,0,0}}%
      \expandafter\def\csname LT1\endcsname{\color[rgb]{0,1,0}}%
      \expandafter\def\csname LT2\endcsname{\color[rgb]{0,0,1}}%
      \expandafter\def\csname LT3\endcsname{\color[rgb]{1,0,1}}%
      \expandafter\def\csname LT4\endcsname{\color[rgb]{0,1,1}}%
      \expandafter\def\csname LT5\endcsname{\color[rgb]{1,1,0}}%
      \expandafter\def\csname LT6\endcsname{\color[rgb]{0,0,0}}%
      \expandafter\def\csname LT7\endcsname{\color[rgb]{1,0.3,0}}%
      \expandafter\def\csname LT8\endcsname{\color[rgb]{0.5,0.5,0.5}}%
    \else
      \def\colorrgb#1{\color{black}}%
      \def\colorgray#1{\color[gray]{#1}}%
      \expandafter\def\csname LTw\endcsname{\color{white}}%
      \expandafter\def\csname LTb\endcsname{\color{black}}%
      \expandafter\def\csname LTa\endcsname{\color{black}}%
      \expandafter\def\csname LT0\endcsname{\color{black}}%
      \expandafter\def\csname LT1\endcsname{\color{black}}%
      \expandafter\def\csname LT2\endcsname{\color{black}}%
      \expandafter\def\csname LT3\endcsname{\color{black}}%
      \expandafter\def\csname LT4\endcsname{\color{black}}%
      \expandafter\def\csname LT5\endcsname{\color{black}}%
      \expandafter\def\csname LT6\endcsname{\color{black}}%
      \expandafter\def\csname LT7\endcsname{\color{black}}%
      \expandafter\def\csname LT8\endcsname{\color{black}}%
    \fi
  \fi
  \setlength{\unitlength}{0.0500bp}%
  \begin{picture}(5760.00,4608.00)%
    \gplgaddtomacro\gplbacktext{%
      \csname LTb\endcsname%
      \put(986,480){\rotatebox{90}{\makebox(0,0){\strut{} 0.85}}}%
      \put(986,1680){\rotatebox{90}{\makebox(0,0){\strut{} 0.9}}}%
      \put(986,2880){\rotatebox{90}{\makebox(0,0){\strut{} 0.95}}}%
      \put(986,4080){\rotatebox{90}{\makebox(0,0){\strut{} 1}}}%
      \put(1128,240){\makebox(0,0){\strut{} 0}}%
      \put(2002,240){\makebox(0,0){\strut{} 0.2}}%
      \put(2875,240){\makebox(0,0){\strut{} 0.4}}%
      \put(3749,240){\makebox(0,0){\strut{} 0.6}}%
      \put(4622,240){\makebox(0,0){\strut{} 0.8}}%
      \put(5496,240){\makebox(0,0){\strut{} 1}}%
      \put(528,2400){\rotatebox{90}{\makebox(0,0){\strut{}}}}%
      \put(3312,-72){\makebox(0,0){\strut{}}}%
    }%
    \gplgaddtomacro\gplfronttext{%
    }%
    \gplbacktext
    \put(0,0){\includegraphics{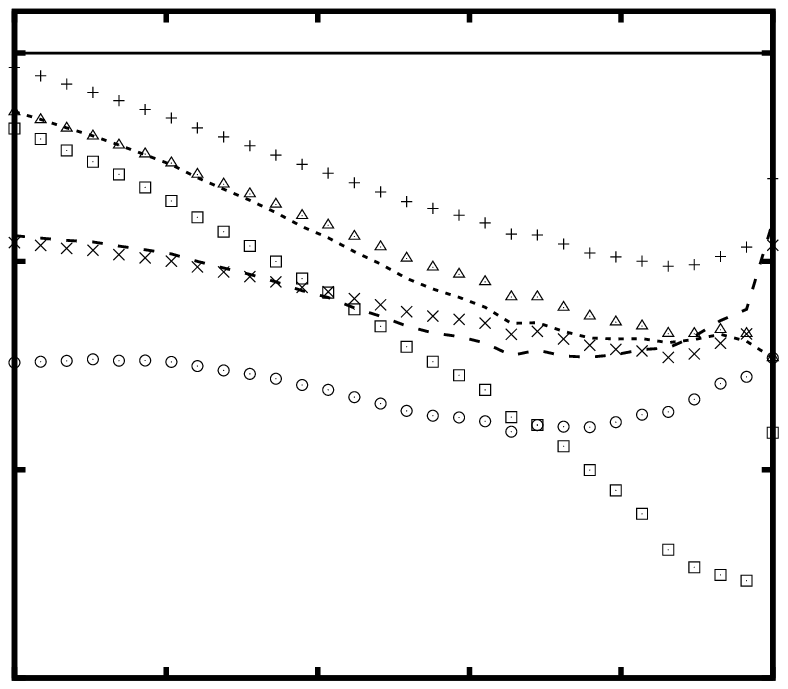}}%
    \gplfronttext
  \end{picture}%
\endgroup

%% file: fig_rho_2_fdr.tex
\begingroup
  \makeatletter
  \providecommand\color[2][]{%
    \GenericError{(gnuplot) \space\space\space\@spaces}{%
      Package color not loaded in conjunction with
      terminal option `colourtext'%
    }{See the gnuplot documentation for explanation.%
    }{Either use 'blacktext' in gnuplot or load the package
      color.sty in LaTeX.}%
    \renewcommand\color[2][]{}%
  }%
  \providecommand\includegraphics[2][]{%
    \GenericError{(gnuplot) \space\space\space\@spaces}{%
      Package graphicx or graphics not loaded%
    }{See the gnuplot documentation for explanation.%
    }{The gnuplot epslatex terminal needs graphicx.sty or graphics.sty.}%
    \renewcommand\includegraphics[2][]{}%
  }%
  \providecommand\rotatebox[2]{#2}%
  \@ifundefined{ifGPcolor}{%
    \newif\ifGPcolor
    \GPcolorfalse
  }{}%
  \@ifundefined{ifGPblacktext}{%
    \newif\ifGPblacktext
    \GPblacktexttrue
  }{}%
  \let\gplgaddtomacro\g@addto@macro
  \gdef\gplbacktext{}%
  \gdef\gplfronttext{}%
  \makeatother
  \ifGPblacktext
    \def\colorrgb#1{}%
    \def\colorgray#1{}%
  \else
    \ifGPcolor
      \def\colorrgb#1{\color[rgb]{#1}}%
      \def\colorgray#1{\color[gray]{#1}}%
      \expandafter\def\csname LTw\endcsname{\color{white}}%
      \expandafter\def\csname LTb\endcsname{\color{black}}%
      \expandafter\def\csname LTa\endcsname{\color{black}}%
      \expandafter\def\csname LT0\endcsname{\color[rgb]{1,0,0}}%
      \expandafter\def\csname LT1\endcsname{\color[rgb]{0,1,0}}%
      \expandafter\def\csname LT2\endcsname{\color[rgb]{0,0,1}}%
      \expandafter\def\csname LT3\endcsname{\color[rgb]{1,0,1}}%
      \expandafter\def\csname LT4\endcsname{\color[rgb]{0,1,1}}%
      \expandafter\def\csname LT5\endcsname{\color[rgb]{1,1,0}}%
      \expandafter\def\csname LT6\endcsname{\color[rgb]{0,0,0}}%
      \expandafter\def\csname LT7\endcsname{\color[rgb]{1,0.3,0}}%
      \expandafter\def\csname LT8\endcsname{\color[rgb]{0.5,0.5,0.5}}%
    \else
      \def\colorrgb#1{\color{black}}%
      \def\colorgray#1{\color[gray]{#1}}%
      \expandafter\def\csname LTw\endcsname{\color{white}}%
      \expandafter\def\csname LTb\endcsname{\color{black}}%
      \expandafter\def\csname LTa\endcsname{\color{black}}%
      \expandafter\def\csname LT0\endcsname{\color{black}}%
      \expandafter\def\csname LT1\endcsname{\color{black}}%
      \expandafter\def\csname LT2\endcsname{\color{black}}%
      \expandafter\def\csname LT3\endcsname{\color{black}}%
      \expandafter\def\csname LT4\endcsname{\color{black}}%
      \expandafter\def\csname LT5\endcsname{\color{black}}%
      \expandafter\def\csname LT6\endcsname{\color{black}}%
      \expandafter\def\csname LT7\endcsname{\color{black}}%
      \expandafter\def\csname LT8\endcsname{\color{black}}%
    \fi
  \fi
  \setlength{\unitlength}{0.0500bp}%
  \begin{picture}(5760.00,4608.00)%
    \gplgaddtomacro\gplbacktext{%
      \csname LTb\endcsname%
      \put(986,480){\rotatebox{90}{\makebox(0,0){\strut{} 0.02}}}%
      \put(986,1760){\rotatebox{90}{\makebox(0,0){\strut{} 0.04}}}%
      \put(986,3040){\rotatebox{90}{\makebox(0,0){\strut{} 0.06}}}%
      \put(986,4320){\rotatebox{90}{\makebox(0,0){\strut{} 0.08}}}%
      \put(1128,240){\makebox(0,0){\strut{} 0}}%
      \put(2002,240){\makebox(0,0){\strut{} 0.2}}%
      \put(2875,240){\makebox(0,0){\strut{} 0.4}}%
      \put(3749,240){\makebox(0,0){\strut{} 0.6}}%
      \put(4622,240){\makebox(0,0){\strut{} 0.8}}%
      \put(5496,240){\makebox(0,0){\strut{} 1}}%
      \put(528,2400){\rotatebox{90}{\makebox(0,0){\strut{}}}}%
      \put(3312,-72){\makebox(0,0){\strut{}}}%
    }%
    \gplgaddtomacro\gplfronttext{%
      \csname LTb\endcsname%
      \put(4200,4125){\makebox(0,0)[l]{\strut{}LSU-Oracle}}%
      \csname LTb\endcsname%
      \put(4200,3909){\makebox(0,0)[l]{\strut{}Storey-$\alpha$}}%
      \csname LTb\endcsname%
      \put(4200,3693){\makebox(0,0)[l]{\strut{}Storey-1/2}}%
      \csname LTb\endcsname%
      \put(4200,3477){\makebox(0,0)[l]{\strut{}Median LSU}}%
      \csname LTb\endcsname%
      \put(4200,3261){\makebox(0,0)[l]{\strut{}BKY06}}%
      \csname LTb\endcsname%
      \put(4200,3045){\makebox(0,0)[l]{\strut{}BR-1S-$\alpha$}}%
      \csname LTb\endcsname%
      \put(4200,2829){\makebox(0,0)[l]{\strut{}BR-2S-$\alpha$}}%
      \csname LTb\endcsname%
      \put(4200,2613){\makebox(0,0)[l]{\strut{}FDR09-1/2}}%
    }%
    \gplbacktext
    \put(0,0){\includegraphics{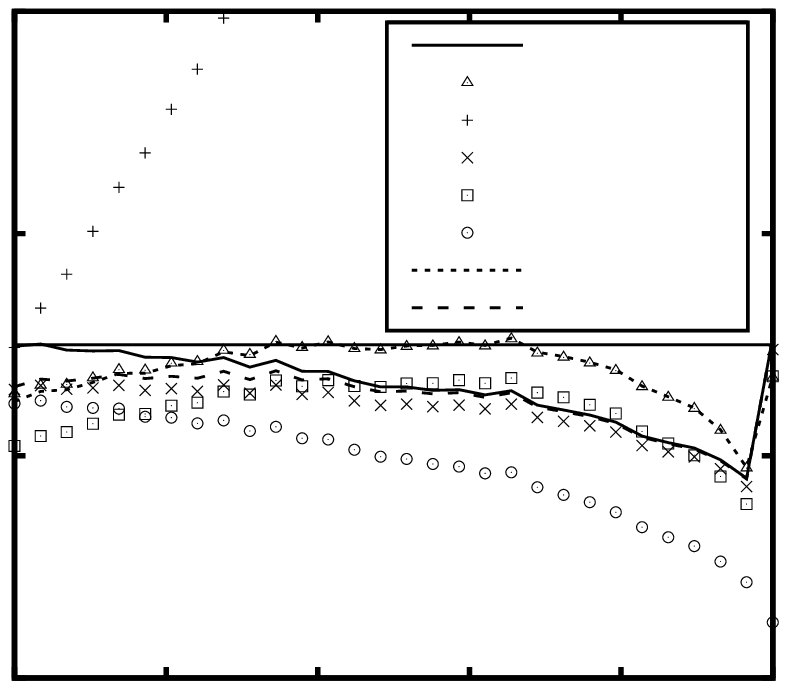}}%
    \gplfronttext
  \end{picture}%
\endgroup

%% file: fig_rho_2_puis.tex
\begingroup
  \makeatletter
  \providecommand\color[2][]{%
    \GenericError{(gnuplot) \space\space\space\@spaces}{%
      Package color not loaded in conjunction with
      terminal option `colourtext'%
    }{See the gnuplot documentation for explanation.%
    }{Either use 'blacktext' in gnuplot or load the package
      color.sty in LaTeX.}%
    \renewcommand\color[2][]{}%
  }%
  \providecommand\includegraphics[2][]{%
    \GenericError{(gnuplot) \space\space\space\@spaces}{%
      Package graphicx or graphics not loaded%
    }{See the gnuplot documentation for explanation.%
    }{The gnuplot epslatex terminal needs graphicx.sty or graphics.sty.}%
    \renewcommand\includegraphics[2][]{}%
  }%
  \providecommand\rotatebox[2]{#2}%
  \@ifundefined{ifGPcolor}{%
    \newif\ifGPcolor
    \GPcolorfalse
  }{}%
  \@ifundefined{ifGPblacktext}{%
    \newif\ifGPblacktext
    \GPblacktexttrue
  }{}%
  \let\gplgaddtomacro\g@addto@macro
  \gdef\gplbacktext{}%
  \gdef\gplfronttext{}%
  \makeatother
  \ifGPblacktext
    \def\colorrgb#1{}%
    \def\colorgray#1{}%
  \else
    \ifGPcolor
      \def\colorrgb#1{\color[rgb]{#1}}%
      \def\colorgray#1{\color[gray]{#1}}%
      \expandafter\def\csname LTw\endcsname{\color{white}}%
      \expandafter\def\csname LTb\endcsname{\color{black}}%
      \expandafter\def\csname LTa\endcsname{\color{black}}%
      \expandafter\def\csname LT0\endcsname{\color[rgb]{1,0,0}}%
      \expandafter\def\csname LT1\endcsname{\color[rgb]{0,1,0}}%
      \expandafter\def\csname LT2\endcsname{\color[rgb]{0,0,1}}%
      \expandafter\def\csname LT3\endcsname{\color[rgb]{1,0,1}}%
      \expandafter\def\csname LT4\endcsname{\color[rgb]{0,1,1}}%
      \expandafter\def\csname LT5\endcsname{\color[rgb]{1,1,0}}%
      \expandafter\def\csname LT6\endcsname{\color[rgb]{0,0,0}}%
      \expandafter\def\csname LT7\endcsname{\color[rgb]{1,0.3,0}}%
      \expandafter\def\csname LT8\endcsname{\color[rgb]{0.5,0.5,0.5}}%
    \else
      \def\colorrgb#1{\color{black}}%
      \def\colorgray#1{\color[gray]{#1}}%
      \expandafter\def\csname LTw\endcsname{\color{white}}%
      \expandafter\def\csname LTb\endcsname{\color{black}}%
      \expandafter\def\csname LTa\endcsname{\color{black}}%
      \expandafter\def\csname LT0\endcsname{\color{black}}%
      \expandafter\def\csname LT1\endcsname{\color{black}}%
      \expandafter\def\csname LT2\endcsname{\color{black}}%
      \expandafter\def\csname LT3\endcsname{\color{black}}%
      \expandafter\def\csname LT4\endcsname{\color{black}}%
      \expandafter\def\csname LT5\endcsname{\color{black}}%
      \expandafter\def\csname LT6\endcsname{\color{black}}%
      \expandafter\def\csname LT7\endcsname{\color{black}}%
      \expandafter\def\csname LT8\endcsname{\color{black}}%
    \fi
  \fi
  \setlength{\unitlength}{0.0500bp}%
  \begin{picture}(5760.00,4608.00)%
    \gplgaddtomacro\gplbacktext{%
      \csname LTb\endcsname%
      \put(986,480){\rotatebox{90}{\makebox(0,0){\strut{} 0.85}}}%
      \put(986,1680){\rotatebox{90}{\makebox(0,0){\strut{} 0.9}}}%
      \put(986,2880){\rotatebox{90}{\makebox(0,0){\strut{} 0.95}}}%
      \put(986,4080){\rotatebox{90}{\makebox(0,0){\strut{} 1}}}%
      \put(1128,240){\makebox(0,0){\strut{} 0}}%
      \put(2002,240){\makebox(0,0){\strut{} 0.2}}%
      \put(2875,240){\makebox(0,0){\strut{} 0.4}}%
      \put(3749,240){\makebox(0,0){\strut{} 0.6}}%
      \put(4622,240){\makebox(0,0){\strut{} 0.8}}%
      \put(5496,240){\makebox(0,0){\strut{} 1}}%
      \put(528,2400){\rotatebox{90}{\makebox(0,0){\strut{}}}}%
      \put(3312,-72){\makebox(0,0){\strut{}}}%
    }%
    \gplgaddtomacro\gplfronttext{%
      \csname LTb\endcsname%
      \put(2199,1971){\makebox(0,0)[l]{\strut{}Storey-$\alpha$}}%
      \csname LTb\endcsname%
      \put(2199,1755){\makebox(0,0)[l]{\strut{}Storey-1/2}}%
      \csname LTb\endcsname%
      \put(2199,1539){\makebox(0,0)[l]{\strut{}Median LSU}}%
      \csname LTb\endcsname%
      \put(2199,1323){\makebox(0,0)[l]{\strut{}BKY06}}%
      \csname LTb\endcsname%
      \put(2199,1107){\makebox(0,0)[l]{\strut{}BR-1S-$\alpha$}}%
      \csname LTb\endcsname%
      \put(2199,891){\makebox(0,0)[l]{\strut{}BR-2S-$\alpha$}}%
      \csname LTb\endcsname%
      \put(2199,675){\makebox(0,0)[l]{\strut{}FDR09-1/2}}%
    }%
    \gplbacktext
    \put(0,0){\includegraphics{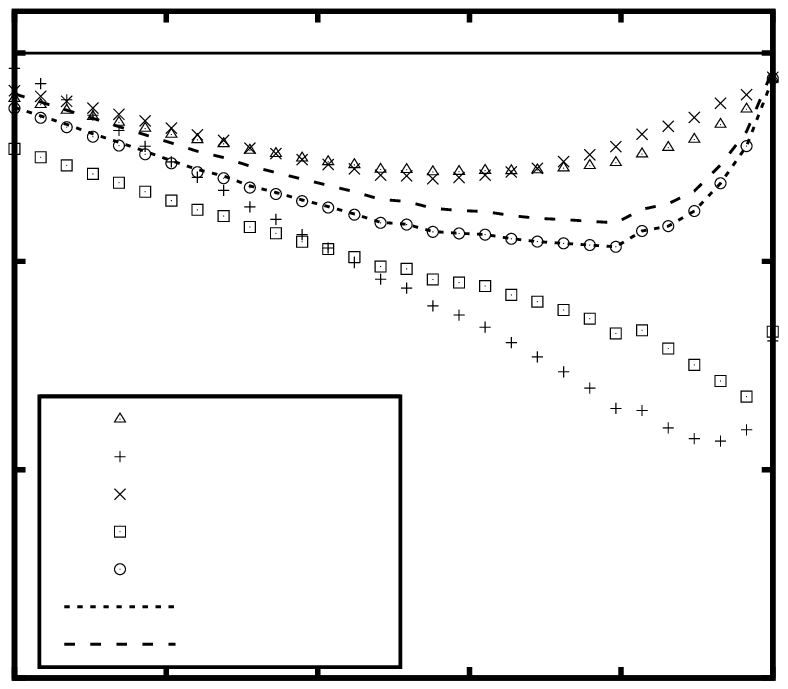}}%
    \gplfronttext
  \end{picture}%
\endgroup

%% file: fig_rho_3_fdr.tex
\begingroup
  \makeatletter
  \providecommand\color[2][]{%
    \GenericError{(gnuplot) \space\space\space\@spaces}{%
      Package color not loaded in conjunction with
      terminal option `colourtext'%
    }{See the gnuplot documentation for explanation.%
    }{Either use 'blacktext' in gnuplot or load the package
      color.sty in LaTeX.}%
    \renewcommand\color[2][]{}%
  }%
  \providecommand\includegraphics[2][]{%
    \GenericError{(gnuplot) \space\space\space\@spaces}{%
      Package graphicx or graphics not loaded%
    }{See the gnuplot documentation for explanation.%
    }{The gnuplot epslatex terminal needs graphicx.sty or graphics.sty.}%
    \renewcommand\includegraphics[2][]{}%
  }%
  \providecommand\rotatebox[2]{#2}%
  \@ifundefined{ifGPcolor}{%
    \newif\ifGPcolor
    \GPcolorfalse
  }{}%
  \@ifundefined{ifGPblacktext}{%
    \newif\ifGPblacktext
    \GPblacktexttrue
  }{}%
  \let\gplgaddtomacro\g@addto@macro
  \gdef\gplbacktext{}%
  \gdef\gplfronttext{}%
  \makeatother
  \ifGPblacktext
    \def\colorrgb#1{}%
    \def\colorgray#1{}%
  \else
    \ifGPcolor
      \def\colorrgb#1{\color[rgb]{#1}}%
      \def\colorgray#1{\color[gray]{#1}}%
      \expandafter\def\csname LTw\endcsname{\color{white}}%
      \expandafter\def\csname LTb\endcsname{\color{black}}%
      \expandafter\def\csname LTa\endcsname{\color{black}}%
      \expandafter\def\csname LT0\endcsname{\color[rgb]{1,0,0}}%
      \expandafter\def\csname LT1\endcsname{\color[rgb]{0,1,0}}%
      \expandafter\def\csname LT2\endcsname{\color[rgb]{0,0,1}}%
      \expandafter\def\csname LT3\endcsname{\color[rgb]{1,0,1}}%
      \expandafter\def\csname LT4\endcsname{\color[rgb]{0,1,1}}%
      \expandafter\def\csname LT5\endcsname{\color[rgb]{1,1,0}}%
      \expandafter\def\csname LT6\endcsname{\color[rgb]{0,0,0}}%
      \expandafter\def\csname LT7\endcsname{\color[rgb]{1,0.3,0}}%
      \expandafter\def\csname LT8\endcsname{\color[rgb]{0.5,0.5,0.5}}%
    \else
      \def\colorrgb#1{\color{black}}%
      \def\colorgray#1{\color[gray]{#1}}%
      \expandafter\def\csname LTw\endcsname{\color{white}}%
      \expandafter\def\csname LTb\endcsname{\color{black}}%
      \expandafter\def\csname LTa\endcsname{\color{black}}%
      \expandafter\def\csname LT0\endcsname{\color{black}}%
      \expandafter\def\csname LT1\endcsname{\color{black}}%
      \expandafter\def\csname LT2\endcsname{\color{black}}%
      \expandafter\def\csname LT3\endcsname{\color{black}}%
      \expandafter\def\csname LT4\endcsname{\color{black}}%
      \expandafter\def\csname LT5\endcsname{\color{black}}%
      \expandafter\def\csname LT6\endcsname{\color{black}}%
      \expandafter\def\csname LT7\endcsname{\color{black}}%
      \expandafter\def\csname LT8\endcsname{\color{black}}%
    \fi
  \fi
  \setlength{\unitlength}{0.0500bp}%
  \begin{picture}(5760.00,4608.00)%
    \gplgaddtomacro\gplbacktext{%
      \csname LTb\endcsname%
      \put(986,480){\rotatebox{90}{\makebox(0,0){\strut{} 0.02}}}%
      \put(986,1760){\rotatebox{90}{\makebox(0,0){\strut{} 0.04}}}%
      \put(986,3040){\rotatebox{90}{\makebox(0,0){\strut{} 0.06}}}%
      \put(986,4320){\rotatebox{90}{\makebox(0,0){\strut{} 0.08}}}%
      \put(1128,240){\makebox(0,0){\strut{} 0}}%
      \put(2002,240){\makebox(0,0){\strut{} 0.2}}%
      \put(2875,240){\makebox(0,0){\strut{} 0.4}}%
      \put(3749,240){\makebox(0,0){\strut{} 0.6}}%
      \put(4622,240){\makebox(0,0){\strut{} 0.8}}%
      \put(5496,240){\makebox(0,0){\strut{} 1}}%
      \put(528,2400){\rotatebox{90}{\makebox(0,0){\strut{}}}}%
      \put(3312,-72){\makebox(0,0){\strut{}}}%
    }%
    \gplgaddtomacro\gplfronttext{%
    }%
    \gplbacktext
    \put(0,0){\includegraphics{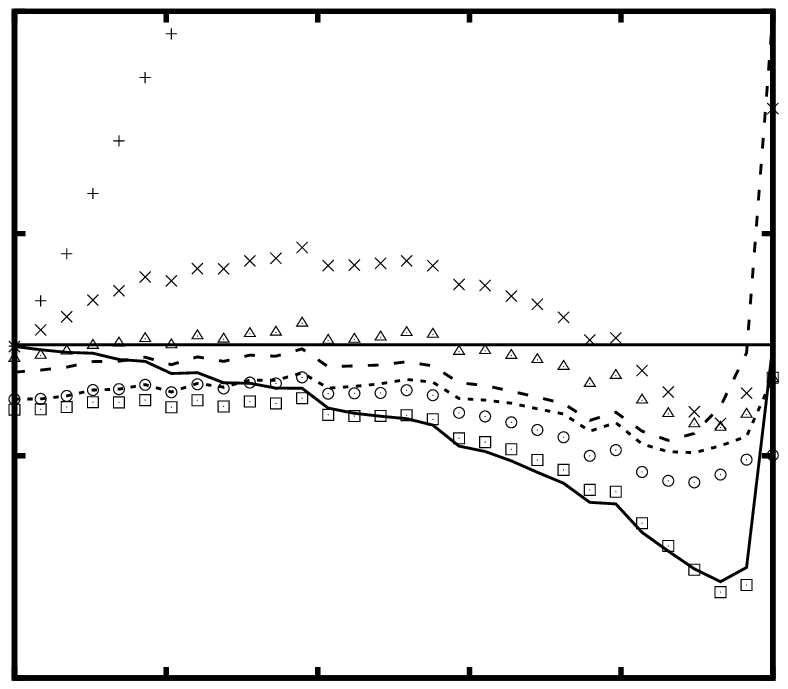}}%
    \gplfronttext
  \end{picture}%
\endgroup

%% file: fig_rho_3_puis.tex
\begingroup
  \makeatletter
  \providecommand\color[2][]{%
    \GenericError{(gnuplot) \space\space\space\@spaces}{%
      Package color not loaded in conjunction with
      terminal option `colourtext'%
    }{See the gnuplot documentation for explanation.%
    }{Either use 'blacktext' in gnuplot or load the package
      color.sty in LaTeX.}%
    \renewcommand\color[2][]{}%
  }%
  \providecommand\includegraphics[2][]{%
    \GenericError{(gnuplot) \space\space\space\@spaces}{%
      Package graphicx or graphics not loaded%
    }{See the gnuplot documentation for explanation.%
    }{The gnuplot epslatex terminal needs graphicx.sty or graphics.sty.}%
    \renewcommand\includegraphics[2][]{}%
  }%
  \providecommand\rotatebox[2]{#2}%
  \@ifundefined{ifGPcolor}{%
    \newif\ifGPcolor
    \GPcolorfalse
  }{}%
  \@ifundefined{ifGPblacktext}{%
    \newif\ifGPblacktext
    \GPblacktexttrue
  }{}%
  \let\gplgaddtomacro\g@addto@macro
  \gdef\gplbacktext{}%
  \gdef\gplfronttext{}%
  \makeatother
  \ifGPblacktext
    \def\colorrgb#1{}%
    \def\colorgray#1{}%
  \else
    \ifGPcolor
      \def\colorrgb#1{\color[rgb]{#1}}%
      \def\colorgray#1{\color[gray]{#1}}%
      \expandafter\def\csname LTw\endcsname{\color{white}}%
      \expandafter\def\csname LTb\endcsname{\color{black}}%
      \expandafter\def\csname LTa\endcsname{\color{black}}%
      \expandafter\def\csname LT0\endcsname{\color[rgb]{1,0,0}}%
      \expandafter\def\csname LT1\endcsname{\color[rgb]{0,1,0}}%
      \expandafter\def\csname LT2\endcsname{\color[rgb]{0,0,1}}%
      \expandafter\def\csname LT3\endcsname{\color[rgb]{1,0,1}}%
      \expandafter\def\csname LT4\endcsname{\color[rgb]{0,1,1}}%
      \expandafter\def\csname LT5\endcsname{\color[rgb]{1,1,0}}%
      \expandafter\def\csname LT6\endcsname{\color[rgb]{0,0,0}}%
      \expandafter\def\csname LT7\endcsname{\color[rgb]{1,0.3,0}}%
      \expandafter\def\csname LT8\endcsname{\color[rgb]{0.5,0.5,0.5}}%
    \else
      \def\colorrgb#1{\color{black}}%
      \def\colorgray#1{\color[gray]{#1}}%
      \expandafter\def\csname LTw\endcsname{\color{white}}%
      \expandafter\def\csname LTb\endcsname{\color{black}}%
      \expandafter\def\csname LTa\endcsname{\color{black}}%
      \expandafter\def\csname LT0\endcsname{\color{black}}%
      \expandafter\def\csname LT1\endcsname{\color{black}}%
      \expandafter\def\csname LT2\endcsname{\color{black}}%
      \expandafter\def\csname LT3\endcsname{\color{black}}%
      \expandafter\def\csname LT4\endcsname{\color{black}}%
      \expandafter\def\csname LT5\endcsname{\color{black}}%
      \expandafter\def\csname LT6\endcsname{\color{black}}%
      \expandafter\def\csname LT7\endcsname{\color{black}}%
      \expandafter\def\csname LT8\endcsname{\color{black}}%
    \fi
  \fi
  \setlength{\unitlength}{0.0500bp}%
  \begin{picture}(5760.00,4608.00)%
    \gplgaddtomacro\gplbacktext{%
      \csname LTb\endcsname%
      \put(986,480){\rotatebox{90}{\makebox(0,0){\strut{} 0.85}}}%
      \put(986,1680){\rotatebox{90}{\makebox(0,0){\strut{} 0.9}}}%
      \put(986,2880){\rotatebox{90}{\makebox(0,0){\strut{} 0.95}}}%
      \put(986,4080){\rotatebox{90}{\makebox(0,0){\strut{} 1}}}%
      \put(1128,240){\makebox(0,0){\strut{} 0}}%
      \put(2002,240){\makebox(0,0){\strut{} 0.2}}%
      \put(2875,240){\makebox(0,0){\strut{} 0.4}}%
      \put(3749,240){\makebox(0,0){\strut{} 0.6}}%
      \put(4622,240){\makebox(0,0){\strut{} 0.8}}%
      \put(5496,240){\makebox(0,0){\strut{} 1}}%
      \put(528,2400){\rotatebox{90}{\makebox(0,0){\strut{}}}}%
      \put(3312,-72){\makebox(0,0){\strut{}}}%
    }%
    \gplgaddtomacro\gplfronttext{%
    }%
    \gplbacktext
    \put(0,0){\includegraphics{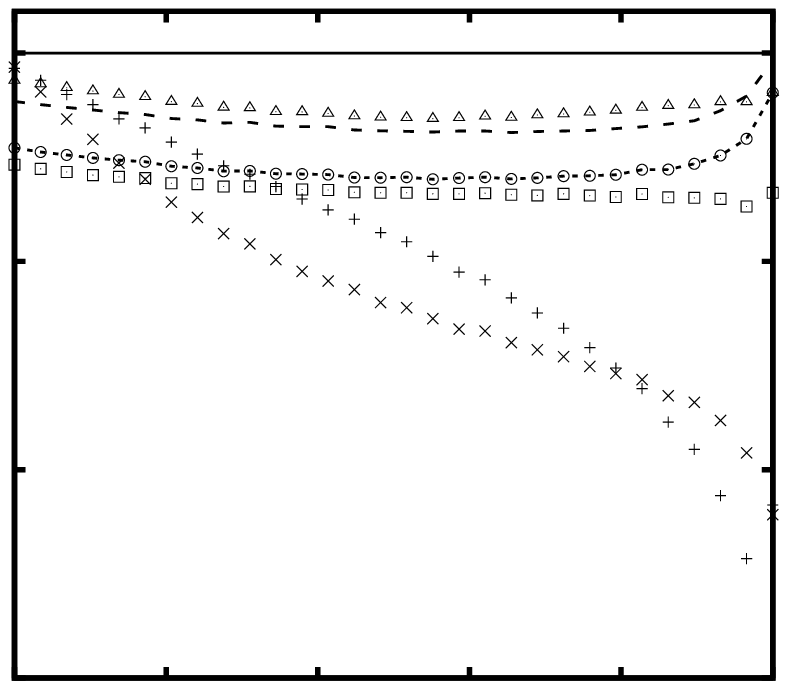}}%
    \gplfronttext
  \end{picture}%
\endgroup

%% file: fig_dep_puis.tex
\begingroup
  \makeatletter
  \providecommand\color[2][]{%
    \GenericError{(gnuplot) \space\space\space\@spaces}{%
      Package color not loaded in conjunction with
      terminal option `colourtext'%
    }{See the gnuplot documentation for explanation.%
    }{Either use 'blacktext' in gnuplot or load the package
      color.sty in LaTeX.}%
    \renewcommand\color[2][]{}%
  }%
  \providecommand\includegraphics[2][]{%
    \GenericError{(gnuplot) \space\space\space\@spaces}{%
      Package graphicx or graphics not loaded%
    }{See the gnuplot documentation for explanation.%
    }{The gnuplot epslatex terminal needs graphicx.sty or graphics.sty.}%
    \renewcommand\includegraphics[2][]{}%
  }%
  \providecommand\rotatebox[2]{#2}%
  \@ifundefined{ifGPcolor}{%
    \newif\ifGPcolor
    \GPcolorfalse
  }{}%
  \@ifundefined{ifGPblacktext}{%
    \newif\ifGPblacktext
    \GPblacktexttrue
  }{}%
  \let\gplgaddtomacro\g@addto@macro
  \gdef\gplbacktext{}%
  \gdef\gplfronttext{}%
  \makeatother
  \ifGPblacktext
    \def\colorrgb#1{}%
    \def\colorgray#1{}%
  \else
    \ifGPcolor
      \def\colorrgb#1{\color[rgb]{#1}}%
      \def\colorgray#1{\color[gray]{#1}}%
      \expandafter\def\csname LTw\endcsname{\color{white}}%
      \expandafter\def\csname LTb\endcsname{\color{black}}%
      \expandafter\def\csname LTa\endcsname{\color{black}}%
      \expandafter\def\csname LT0\endcsname{\color[rgb]{1,0,0}}%
      \expandafter\def\csname LT1\endcsname{\color[rgb]{0,1,0}}%
      \expandafter\def\csname LT2\endcsname{\color[rgb]{0,0,1}}%
      \expandafter\def\csname LT3\endcsname{\color[rgb]{1,0,1}}%
      \expandafter\def\csname LT4\endcsname{\color[rgb]{0,1,1}}%
      \expandafter\def\csname LT5\endcsname{\color[rgb]{1,1,0}}%
      \expandafter\def\csname LT6\endcsname{\color[rgb]{0,0,0}}%
      \expandafter\def\csname LT7\endcsname{\color[rgb]{1,0.3,0}}%
      \expandafter\def\csname LT8\endcsname{\color[rgb]{0.5,0.5,0.5}}%
    \else
      \def\colorrgb#1{\color{black}}%
      \def\colorgray#1{\color[gray]{#1}}%
      \expandafter\def\csname LTw\endcsname{\color{white}}%
      \expandafter\def\csname LTb\endcsname{\color{black}}%
      \expandafter\def\csname LTa\endcsname{\color{black}}%
      \expandafter\def\csname LT0\endcsname{\color{black}}%
      \expandafter\def\csname LT1\endcsname{\color{black}}%
      \expandafter\def\csname LT2\endcsname{\color{black}}%
      \expandafter\def\csname LT3\endcsname{\color{black}}%
      \expandafter\def\csname LT4\endcsname{\color{black}}%
      \expandafter\def\csname LT5\endcsname{\color{black}}%
      \expandafter\def\csname LT6\endcsname{\color{black}}%
      \expandafter\def\csname LT7\endcsname{\color{black}}%
      \expandafter\def\csname LT8\endcsname{\color{black}}%
    \fi
  \fi
  \setlength{\unitlength}{0.0500bp}%
  \begin{picture}(5760.00,4608.00)%
    \gplgaddtomacro\gplbacktext{%
      \csname LTb\endcsname%
      \put(705,480){\rotatebox{90}{\makebox(0,0){\strut{} 0}}}%
      \put(705,1248){\rotatebox{90}{\makebox(0,0){\strut{} 0.2}}}%
      \put(705,2016){\rotatebox{90}{\makebox(0,0){\strut{} 0.4}}}%
      \put(705,2784){\rotatebox{90}{\makebox(0,0){\strut{} 0.6}}}%
      \put(705,3552){\rotatebox{90}{\makebox(0,0){\strut{} 0.8}}}%
      \put(705,4320){\rotatebox{90}{\makebox(0,0){\strut{} 1}}}%
      \put(840,240){\makebox(0,0){\strut{} 0.5}}%
      \put(1341,240){\makebox(0,0){\strut{} 1}}%
      \put(1843,240){\makebox(0,0){\strut{} 1.5}}%
      \put(2344,240){\makebox(0,0){\strut{} 2}}%
      \put(2845,240){\makebox(0,0){\strut{} 2.5}}%
      \put(3347,240){\makebox(0,0){\strut{} 3}}%
      \put(3848,240){\makebox(0,0){\strut{} 3.5}}%
      \put(4349,240){\makebox(0,0){\strut{} 4}}%
      \put(4851,240){\makebox(0,0){\strut{} 4.5}}%
      \put(5352,240){\makebox(0,0){\strut{} 5}}%
      \put(240,2400){\rotatebox{90}{\makebox(0,0){\strut{}}}}%
      \put(3096,-120){\makebox(0,0){\strut{}}}%
    }%
    \gplgaddtomacro\gplfronttext{%
      \csname LTb\endcsname%
      \put(3768,1107){\makebox(0,0)[l]{\strut{}LSU}}%
      \csname LTb\endcsname%
      \put(3768,891){\makebox(0,0)[l]{\strut{}BR-dep}}%
      \csname LTb\endcsname%
      \put(3768,675){\makebox(0,0)[l]{\strut{}BR-dep-Holm}}%
    }%
    \gplbacktext
    \put(0,0){\includegraphics{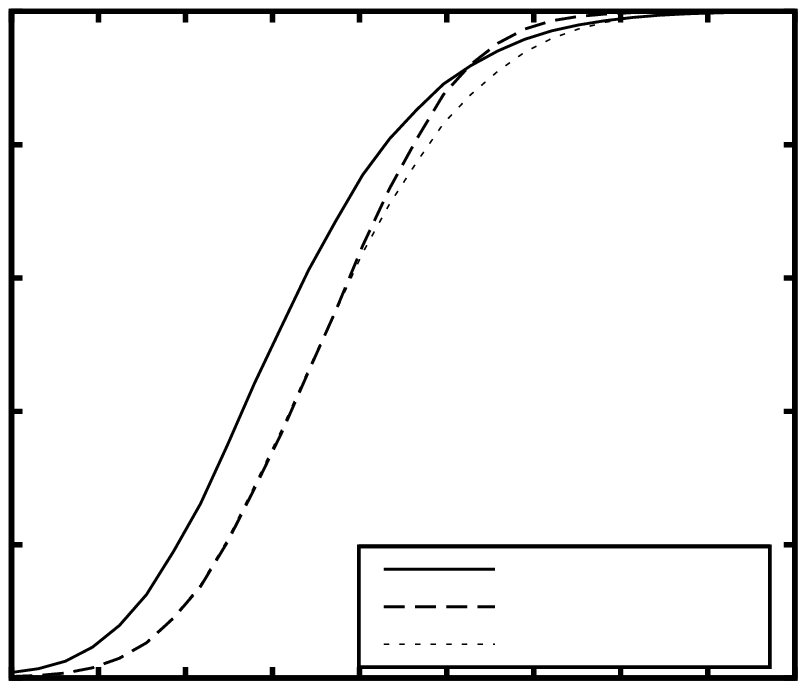}}%
    \gplfronttext
  \end{picture}%
\endgroup

%% file: fig_dep_fnr.tex
\begingroup
  \makeatletter
  \providecommand\color[2][]{%
    \GenericError{(gnuplot) \space\space\space\@spaces}{%
      Package color not loaded in conjunction with
      terminal option `colourtext'%
    }{See the gnuplot documentation for explanation.%
    }{Either use 'blacktext' in gnuplot or load the package
      color.sty in LaTeX.}%
    \renewcommand\color[2][]{}%
  }%
  \providecommand\includegraphics[2][]{%
    \GenericError{(gnuplot) \space\space\space\@spaces}{%
      Package graphicx or graphics not loaded%
    }{See the gnuplot documentation for explanation.%
    }{The gnuplot epslatex terminal needs graphicx.sty or graphics.sty.}%
    \renewcommand\includegraphics[2][]{}%
  }%
  \providecommand\rotatebox[2]{#2}%
  \@ifundefined{ifGPcolor}{%
    \newif\ifGPcolor
    \GPcolorfalse
  }{}%
  \@ifundefined{ifGPblacktext}{%
    \newif\ifGPblacktext
    \GPblacktexttrue
  }{}%
  \let\gplgaddtomacro\g@addto@macro
  \gdef\gplbacktext{}%
  \gdef\gplfronttext{}%
  \makeatother
  \ifGPblacktext
    \def\colorrgb#1{}%
    \def\colorgray#1{}%
  \else
    \ifGPcolor
      \def\colorrgb#1{\color[rgb]{#1}}%
      \def\colorgray#1{\color[gray]{#1}}%
      \expandafter\def\csname LTw\endcsname{\color{white}}%
      \expandafter\def\csname LTb\endcsname{\color{black}}%
      \expandafter\def\csname LTa\endcsname{\color{black}}%
      \expandafter\def\csname LT0\endcsname{\color[rgb]{1,0,0}}%
      \expandafter\def\csname LT1\endcsname{\color[rgb]{0,1,0}}%
      \expandafter\def\csname LT2\endcsname{\color[rgb]{0,0,1}}%
      \expandafter\def\csname LT3\endcsname{\color[rgb]{1,0,1}}%
      \expandafter\def\csname LT4\endcsname{\color[rgb]{0,1,1}}%
      \expandafter\def\csname LT5\endcsname{\color[rgb]{1,1,0}}%
      \expandafter\def\csname LT6\endcsname{\color[rgb]{0,0,0}}%
      \expandafter\def\csname LT7\endcsname{\color[rgb]{1,0.3,0}}%
      \expandafter\def\csname LT8\endcsname{\color[rgb]{0.5,0.5,0.5}}%
    \else
      \def\colorrgb#1{\color{black}}%
      \def\colorgray#1{\color[gray]{#1}}%
      \expandafter\def\csname LTw\endcsname{\color{white}}%
      \expandafter\def\csname LTb\endcsname{\color{black}}%
      \expandafter\def\csname LTa\endcsname{\color{black}}%
      \expandafter\def\csname LT0\endcsname{\color{black}}%
      \expandafter\def\csname LT1\endcsname{\color{black}}%
      \expandafter\def\csname LT2\endcsname{\color{black}}%
      \expandafter\def\csname LT3\endcsname{\color{black}}%
      \expandafter\def\csname LT4\endcsname{\color{black}}%
      \expandafter\def\csname LT5\endcsname{\color{black}}%
      \expandafter\def\csname LT6\endcsname{\color{black}}%
      \expandafter\def\csname LT7\endcsname{\color{black}}%
      \expandafter\def\csname LT8\endcsname{\color{black}}%
    \fi
  \fi
  \setlength{\unitlength}{0.0500bp}%
  \begin{picture}(5760.00,4608.00)%
    \gplgaddtomacro\gplbacktext{%
      \csname LTb\endcsname%
      \put(705,480){\rotatebox{90}{\makebox(0,0){\strut{} 0}}}%
      \put(705,1248){\rotatebox{90}{\makebox(0,0){\strut{} 0.2}}}%
      \put(705,2016){\rotatebox{90}{\makebox(0,0){\strut{} 0.4}}}%
      \put(705,2784){\rotatebox{90}{\makebox(0,0){\strut{} 0.6}}}%
      \put(705,3552){\rotatebox{90}{\makebox(0,0){\strut{} 0.8}}}%
      \put(705,4320){\rotatebox{90}{\makebox(0,0){\strut{} 1}}}%
      \put(840,240){\makebox(0,0){\strut{} 0.5}}%
      \put(1341,240){\makebox(0,0){\strut{} 1}}%
      \put(1843,240){\makebox(0,0){\strut{} 1.5}}%
      \put(2344,240){\makebox(0,0){\strut{} 2}}%
      \put(2845,240){\makebox(0,0){\strut{} 2.5}}%
      \put(3347,240){\makebox(0,0){\strut{} 3}}%
      \put(3848,240){\makebox(0,0){\strut{} 3.5}}%
      \put(4349,240){\makebox(0,0){\strut{} 4}}%
      \put(4851,240){\makebox(0,0){\strut{} 4.5}}%
      \put(5352,240){\makebox(0,0){\strut{} 5}}%
      \put(240,2400){\rotatebox{90}{\makebox(0,0){\strut{}}}}%
      \put(3096,-120){\makebox(0,0){\strut{}}}%
    }%
    \gplgaddtomacro\gplfronttext{%
      \csname LTb\endcsname%
      \put(1911,1107){\makebox(0,0)[l]{\strut{}LSU}}%
      \csname LTb\endcsname%
      \put(1911,891){\makebox(0,0)[l]{\strut{}BR-dep}}%
      \csname LTb\endcsname%
      \put(1911,675){\makebox(0,0)[l]{\strut{}BR-dep-Holm}}%
    }%
    \gplbacktext
    \put(0,0){\includegraphics{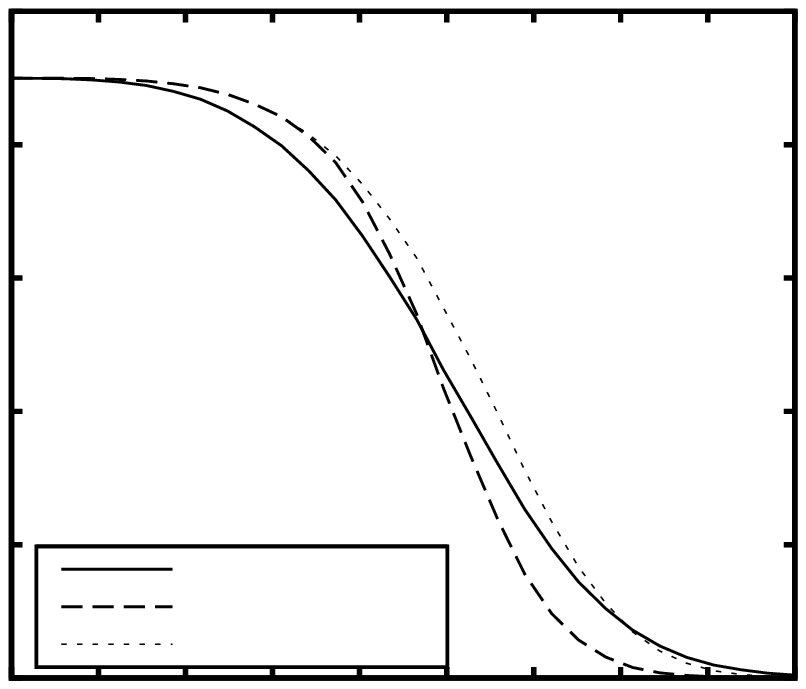}}%
    \gplfronttext
  \end{picture}%
\endgroup

%% file: BR08b.bbl
\begin{thebibliography}{34}
\providecommand{\natexlab}[1]{#1}
\providecommand{\url}[1]{\texttt{#1}}
\expandafter\ifx\csname urlstyle\endcsname\relax
  \providecommand{\doi}[1]{doi: #1}\else
  \providecommand{\doi}{doi: \begingroup \urlstyle{rm}\Url}\fi

\bibitem[Benjamini and Heller(2007)]{BH2006}
Y.~Benjamini and R.~Heller.
\newblock False discovery rates for spatial signals.
\newblock \emph{J. Amer. Statist. Assoc.}, 102\penalty0 (480):\penalty0
  1272--1281, 2007.

\bibitem[Benjamini and Hochberg(1995)]{BH1995}
Y.~Benjamini and Y.~Hochberg.
\newblock Controlling the false discovery rate: a practical and powerful
  approach to multiple testing.
\newblock \emph{J. Roy. Statist. Soc. Ser. B}, 57\penalty0 (1):\penalty0
  289--300, 1995.
\newblock ISSN 0035-9246.

\bibitem[Benjamini and Hochberg(2000)]{BH2000}
Y.~Benjamini and Y.~Hochberg.
\newblock On the adaptive control of the false discovery rate in multiple
  testing with independent statistics.
\newblock \emph{J. Behav. Educ. Statist.}, 25:\penalty0 60--83, 2000.

\bibitem[Benjamini and Yekutieli(2001)]{BY2001}
Y.~Benjamini and D.~Yekutieli.
\newblock The control of the false discovery rate in multiple testing under
  dependency.
\newblock \emph{Ann. Statist.}, 29\penalty0 (4):\penalty0 1165--1188, 2001.
\newblock ISSN 0090-5364.

\bibitem[Benjamini et~al.(2006)Benjamini, Krieger, and Yekutieli]{BKY2006}
Y.~Benjamini, A.~M. Krieger, and D.~Yekutieli.
\newblock Adaptive linear step-up procedures that control the false discovery
  rate.
\newblock \emph{Biometrika}, 93\penalty0 (3):\penalty0 491--507, 2006.
\newblock ISSN 0006-3444.

\bibitem[Black(2004)]{Black2004}
M.~A. Black.
\newblock A note on the adaptive control of false discovery rates.
\newblock \emph{J. R. Stat. Soc. Ser. B Stat. Methodol.}, 66\penalty0
  (2):\penalty0 297--304, 2004.
\newblock ISSN 1369-7412.

\bibitem[Blanchard and Fleuret(2007)]{BF2007}
G.~Blanchard and F.~Fleuret.
\newblock Occam's hammer.
\newblock In N.~Bshouty and C.~Gentile, editors, \emph{Proceedings of the 20th.
  conference on learning theory (COLT 2007)}, volume 4539 of \emph{Springer
  Lecture Notes on Computer Science}, pages 112--126, 2007.

\bibitem[Blanchard and Roquain(2008)]{BR2008EJS}
G.~Blanchard and E.~Roquain.
\newblock Two simple sufficient conditions for fdr control.
\newblock \emph{Electron. J. Stat.}, 2:\penalty0 963--992, 2008.

\bibitem[Dudoit et~al.(2003)Dudoit, Shaffer, and Boldrick]{DSP2003}
S.~Dudoit, J.~P. Shaffer, and J.~C. Boldrick.
\newblock Multiple hypothesis testing in microarray experiments.
\newblock \emph{Statist. Sci.}, 18\penalty0 (1):\penalty0 71--103, 2003.
\newblock ISSN 0883-4237.

\bibitem[Efron et~al.(2001)Efron, Tibshirani, Storey, and Tusher]{ETST2001}
B.~Efron, R.~Tibshirani, J.~D. Storey, and V.~Tusher.
\newblock Empirical {B}ayes analysis of a microarray experiment.
\newblock \emph{J. Amer. Statist. Assoc.}, 96\penalty0 (456):\penalty0
  1151--1160, 2001.
\newblock ISSN 0162-1459.

\bibitem[Farcomeni(2007)]{Far2007}
A.~Farcomeni.
\newblock Some results on the control of the false discovery rate under
  dependence.
\newblock \emph{Scandinavian Journal of Statistics}, 34\penalty0 (2):\penalty0
  275--297, 2007.

\bibitem[Finner and Roters(2001)]{FR2001}
H.~Finner and M.~Roters.
\newblock On the false discovery rate and expected type {I} errors.
\newblock \emph{Biom. J.}, 43\penalty0 (8):\penalty0 985--1005, 2001.
\newblock ISSN 0323-3847.

\bibitem[Finner et~al.(2007)Finner, Dickhaus, and Roters]{FDR2007}
H.~Finner, T.~Dickhaus, and M.~Roters.
\newblock Dependency and false discovery rate: Asymptotics.
\newblock \emph{Ann. Statist.}, 35\penalty0 (4):\penalty0 1432--1455, 2007.

\bibitem[Finner et~al.(2009)Finner, Dickhaus, and Roters]{FDR2008}
H.~Finner, T.~Dickhaus, and M.~Roters.
\newblock On the false discovery rate and an asymptotically optimal rejection
  curve.
\newblock \emph{Ann. Statist.}, 2009.
\newblock To appear.

\bibitem[Gavrilov et~al.(2009)Gavrilov, Benjamini, and Sarkar]{GBS2008}
Y.~Gavrilov, Y.~Benjamini, and S.~K. Sarkar.
\newblock An adaptive step-down procedure with proven {FDR} control under
  independence.
\newblock \emph{Ann. Statist.}, 2009.
\newblock To appear.

\bibitem[Genovese and Wasserman(2002)]{GW2002}
C.~Genovese and L.~Wasserman.
\newblock Operating characteristics and extensions of the false discovery rate
  procedure.
\newblock \emph{J. R. Stat. Soc. Ser. B Stat. Methodol.}, 64\penalty0
  (3):\penalty0 499--517, 2002.
\newblock ISSN 1369-7412.

\bibitem[Genovese and Wasserman(2004)]{GW2004}
C.~Genovese and L.~Wasserman.
\newblock A stochastic process approach to false discovery control.
\newblock \emph{Ann. Statist.}, 32\penalty0 (3):\penalty0 1035--1061, 2004.
\newblock ISSN 0090-5364.

\bibitem[Hoeffding(1963)]{Hoeff1963}
W.~Hoeffding.
\newblock Probability inequalities for sums of bounded random variables.
\newblock \emph{J. Amer. Statist. Assoc.}, 58:\penalty0 13--30, 1963.
\newblock ISSN 0162-1459.

\bibitem[Holm(1979)]{Holm1979}
S.~Holm.
\newblock A simple sequentially rejective multiple test procedure.
\newblock \emph{Scand. J. Statist.}, 6\penalty0 (2):\penalty0 65--70, 1979.
\newblock ISSN 0303-6898.

\bibitem[Jin(2008)]{Jin2008}
J.~Jin.
\newblock Proportion of non-zero normal means: universal oracle equivalences
  and uniformly consistent estimators.
\newblock \emph{J. Roy. Statist. Soc. Ser. B}, 70\penalty0 (3):\penalty0
  461--493, 2008.

\bibitem[Jin and Cai(2007)]{JC2007}
J.~Jin and T.~Cai.
\newblock Estimating the null and the proportion of nonnull effects in
  large-scale multiple comparisons.
\newblock \emph{J. Amer. Statist. Assoc.}, 102\penalty0 (478):\penalty0
  495--506, 2007.
\newblock ISSN 0162-1459.

\bibitem[Lehmann(1966)]{Lehm1966}
E.~L. Lehmann.
\newblock Some concepts of dependence.
\newblock \emph{Ann. Math. Statist.}, 37:\penalty0 1137--1153, 1966.

\bibitem[Meinshausen and Rice(2006)]{MR2006}
N.~Meinshausen and J.~Rice.
\newblock Estimating the proportion of false null hypotheses among a large
  number of independently tested hypotheses.
\newblock \emph{Ann. Statist.}, 34\penalty0 (1):\penalty0 373--393, 2006.
\newblock ISSN 0090-5364.

\bibitem[Neuvial(2008)]{Neu08}
P.~Neuvial.
\newblock Asymptotic properties of false discovery rate controlling procedures
  under independence.
\newblock \emph{Electron. J. Stat.}, 2:\penalty0 1065--1110, 2008.

\bibitem[Roquain(2007)]{Roq2007}
E.~Roquain.
\newblock \emph{Exceptional motifs in heterogeneous sequences. Contributions to
  theory and methodology of multiple testing}.
\newblock PhD thesis, Universit\'e Paris XI, 2007.

\bibitem[Sarkar(2008{\natexlab{a}})]{Sar2006}
S.~K. Sarkar.
\newblock Two-stage stepup procedures controlling {FDR}.
\newblock \emph{Journal of Statistical Planning and Inference}, 138\penalty0
  (4):\penalty0 1072--1084, 2008{\natexlab{a}}.

\bibitem[Sarkar(2008{\natexlab{b}})]{Sar2008}
S.~K. Sarkar.
\newblock On methods controlling the false discovery rate.
\newblock 2008{\natexlab{b}}.

\bibitem[Scott and Blanchard(2009)]{SB2009}
C.~Scott and G.~Blanchard.
\newblock Novelty detection: unlabaled data definitely help.
\newblock In \emph{AISTATS}, 2009.

\bibitem[Spj{\o}tvoll(1972)]{Spjo1972}
E.~Spj{\o}tvoll.
\newblock On the optimality of some multiple comparison procedures.
\newblock \emph{Ann. Math. Statist.}, 43\penalty0 (2):\penalty0 398--411, 1972.

\bibitem[Storey(2002)]{Storey2002}
J.~D. Storey.
\newblock A direct approach to false discovery rates.
\newblock \emph{J. R. Stat. Soc. Ser. B Stat. Methodol.}, 64\penalty0
  (3):\penalty0 479--498, 2002.
\newblock ISSN 1369-7412.

\bibitem[Storey(2007)]{Storey2007}
J.~D. Storey.
\newblock The optimal discovery procedure: a new approach to simultaneous
  significance testing.
\newblock \emph{J. R. Stat. Soc. Ser. B Stat. Methodol.}, 69\penalty0
  (3):\penalty0 347--368, 2007.
\newblock ISSN 1369-7412.

\bibitem[Storey and Tibshirani(2003)]{ST2003}
J.~D. Storey and R.~Tibshirani.
\newblock S{AM} thresholding and false discovery rates for detecting
  differential gene expression in {DNA} microarrays.
\newblock In \emph{The analysis of gene expression data}, Stat. Biol. Health,
  pages 272--290. Springer, New York, 2003.

\bibitem[Storey et~al.(2004)Storey, Taylor, and Siegmund]{STS2004}
J.~D. Storey, J.~E. Taylor, and D.~Siegmund.
\newblock Strong control, conservative point estimation and simultaneous
  conservative consistency of false discovery rates: a unified approach.
\newblock \emph{J. R. Stat. Soc. Ser. B Stat. Methodol.}, 66\penalty0
  (1):\penalty0 187--205, 2004.
\newblock ISSN 1369-7412.

\bibitem[Sun and Cai(2007)]{SC2007}
W.~Sun and T.~Cai.
\newblock Oracle and adaptive compound decision rules for false discovery rate
  control.
\newblock \emph{J. Amer. Statist. Assoc.}, 102\penalty0 (479):\penalty0
  901--912, 2007.
\newblock ISSN 0162-1459.

\end{thebibliography}
